\documentclass[secthm,seceqn,amsthm,ussrhead,12pt]{amsart}
\usepackage[utf8]{inputenc}
\usepackage[english]{babel}

\usepackage{amssymb,amsmath,amsthm,amsfonts,xcolor,enumerate,hyperref,comment,longtable,cleveref}

\usepackage{times}
\usepackage{cite}
\usepackage{pdflscape}
\usepackage{ulem}
\usepackage[mathcal]{euscript}
\usepackage{tikz}
\usepackage{hyperref}
\usepackage{cancel}
\usepackage{stmaryrd}

\usetikzlibrary{arrows}

\newtheorem*{theoremA}{Theorem A}
\newtheorem*{theoremB}{Theorem B}
\newtheorem{lemma}{Lemma}

\theoremstyle{definition}
\newtheorem*{definition}{Definition}

\theoremstyle{remark}
\newtheorem{remark}{Remark}

\usepackage{stmaryrd}
\usepackage{xcolor}

\newcommand{\aut}[1]{\operatorname{\mathrm{Aut}}{(#1)}}

\newcommand{\A}{\mathbf{A}}

\newcommand{\T}[2]{\mathbf{T}^{#1}_{#2}}

\newcommand{\la}{\langle}
\newcommand{\ra}{\rangle}

\newcommand{\Dt}[2]{\Delta_{#1#2}}
\newcommand{\Dl}[2]{[\Delta_{#1#2}]}
\newcommand{\nb}[1]{\nabla_{#1}}

\newcommand{\lb}{\lambda}
\newcommand{\0}{\theta}
\newcommand{\af}{\alpha}

\begin{document}
\noindent{\Large 
One-generated nilpotent terminal algebras}

\

   {\bf   Ivan Kaygorodov$^{a},$  Abror Khudoyberdiyev$^{b,c}$ \& Aloberdi Sattarov$^{c}$\\

    \medskip
}

{\tiny

$^{a}$ CMCC, Universidade Federal do ABC. Santo Andr\'e, Brazil

$^{b}$ National University of Uzbekistan, Tashkent, Uzbekistan

$^{c}$ Institute of Mathematics Academy of Science of Uzbekistan, Tashkent, Uzbekistan

\

\smallskip

   E-mail addresses:

\smallskip

    Ivan Kaygorodov (kaygorodov.ivan@gmail.com) 
    
    Abror Khudoyberdiyev (khabror@mail.ru)    
    
    Aloberdi Sattarov (saloberdi90@mail.ru)
    
}

\ 

\

\ 

\noindent {\bf Abstract.}
 We give an algebraic  classification of complex $5$-dimensional  one-generated nilpotent terminal  algebras. 

\

\noindent {\bf Keywords}: {\it Nilpotent algebra, terminal algebra,  algebraic classification,  central extension.}

\ 

\noindent {\bf MSC2010}: 17A30.

\ 

\section*{Introduction}

Algebraic classification (up to isomorphism) of algebras of small dimension from a certain variety defined by a family of polynomial identities is a classic problem in the theory of non-associative algebras. There are many results related to algebraic classification of small dimensional algebras in varieties of Jordan, Lie, Leibniz, Zinbiel and other algebras.
Another interesting approach of studying algebras of a fixed dimension is to study them from a geometric point of view (that is, to study degenerations and deformations of these algebras). The results in which the complete information about degenerations of a certain variety is obtained are generally referred to as the geometric classification of the algebras of these variety. There are many results related to geometric classification of Jordan, Lie, Leibniz, Zinbiel and other algebras \cite{ack, kv17,kppv,ikv17, kv16}.
Another interesting direction is a study of one-generated objects.
The description of one-generated finite groups is well-known: there is only one one-generated group of order $n$.
    In the case of algebras, there are some similar results,
    such that the description of $n$-dimensional one-generated nilpotent associative \cite{karel}, noncommutative Jordan \cite{jkk19},  Leibniz and Zinbiel algebras\cite{bakhrom}. 
    It was proven that there is only one $n$-dimensional one-generated nilpotent algebra in these varieties.
    But on the other side, as we can see in varieties of Novikov \cite{kkk18}, assosymmetric \cite{ikm19},  bicommutative \cite{kpv19}, commutative \cite{fkkv}, and terminal \cite{kkp19} algebras, there are more than one $4$-dimensional one-generated nilpotent algebra from these varieties.
    One-generated Novikov algebras in dimensions 5 and 6 were studied in \cite{ckkk19}.
In the present paper, we give the algebraic   classification of
$5$-dimensional one-generated nilpotent terminal   algebras, which were
introduced by Kantor   in \cite{Kantor72}.

Our main results are summarized below.

\begin{theoremA}
Up to isomorphism,  the variety of $4$-dimensional complex one-generated nilpotent terminal algebras has infinitely many isomorphism classes, described explicitly in Table A in terms of 
$1$ two-parametric family, $4$ one-parameter families,  and $5$ additional isomorphism classes.
\end{theoremA}

\begin{theoremB}
Up to isomorphism,  the variety of complex $5$-dimensional  one-generated nilpotent terminal algebras has infinitely many isomorphism classes, described explicitly in Table C in terms of 
$1$ three-parametric family, $18$ two-parametric families, $38$ one-parameter families,  and $22$ additional isomorphism classes.
\end{theoremB}

In 1972, Kantor introduced the notion a conservative algebra as a generalization of Jordan algebras \cite{Kantor72}. Unlike other classes of non-associative algebras, this class is not defined by a set of identities. To introduce the notion of a conservative algebra, we need some notation. 
 Let $\mathbb V$ be a vector space, let $A$ be a linear operator on $\mathbb V,$ and let $B$ and $C$ be bilinear operators on $\mathbb V.$ For all $x,y,z \in \mathbb V,$ put
\begin{longtable}{lcl}
$[B,x](y)    $&$=$&$ B(x,y),$\\
$[A,B](x,y)  $&$=$&$ A(B(x,y))- B(A(x),y)-B(x,A(y)),$\\
$[B,C](x,y,z)$&$=$&$ B(C(x,y),z)+B(x,C(y,z))+B(y,C(x,z))$\\
             &&$-C(B(x,y),z)-C(x,B(y,z))-C(y,B(x,z)).$ 
\end{longtable}

Consider an algebra as a vector space $\mathbb V$ over a field $\mathbb{C}$, together with an
 element $\mu$ of $\operatorname{Hom}(\mathbb V \otimes \mathbb V,  \mathbb V),$ so that $a \cdot b =\mu(a \otimes b).$ 
For an algebra $(\mathbb V, \mathcal P)$ with a multiplication $\mathcal P$ and $x\in \mathbb V$ we denote by $L_x^{\mathcal P}$ the operator of left multiplication by $x.$ %If the multiplication $P$ is fixed, we write $L_x$ instead of $L_x^P.$
Thus, Kantor defines conservative algebras as follows:

\begin{definition}
An algebra $(\mathbb V, \mathcal P),$ where $\mathbb V$ is the vector space and $\mathcal P$ is the multiplication, is called a conservative algebra if there is a new multiplication $\mathcal P^*: \mathbb V\times \mathbb V\rightarrow \mathbb V$ such that 
\begin{equation}\label{uno}
[L_b^{\mathcal P},[L_a^{\mathcal P}, {\mathcal P}]]=-[L_{{\mathcal P^*}(a,b)}^{\mathcal P},{\mathcal P}] 
\textrm{, for all $a, b \in \mathbb V.$}
\end{equation}
Simple calculations take us to the following identity with an additional multiplication $\mathcal P^*,$ which must  hold for all $a, b, x, y\in \mathbb V$:
\begin{equation}\label{dos}
\begin{split}
b(a(xy)-(ax)y-x(ay)) - a((bx)y) + (a(bx))y+(bx)(ay) -a(x(by))+(ax)(by)+x(a(by))= \\
= - \mathcal P^*(a,b)(xy)+(\mathcal P^*(a,b)x)y + x(\mathcal P^*(a,b)y).
\end{split}
\end{equation}
\end{definition}

The class of conservative algebras is very vast. It includes
 all associative algebras, 
 all quasi-associative algebras,
 all Jordan algebras, 
 all Lie algebras, 
 all (left) Leibniz algebras,
 all (left) Zinbiel algebras, 
 and many other classes of algebras.
On the other side, all conservative algebras are "rigid" (in sense of Cantarini and Kac \cite{kacan}).

However, this class is very hard to study and for now even the basic general questions about it remain unanswered, so it is a good idea to study its subclasses which are sufficiently wide but easier to deal with. In \cite{Kantor89} Kantor, studying the generalized TKK functor, introduced the class of terminal algebras which is a subclass of the class of conservative algebras. 

\begin{definition}
An algebra $(\mathbb V, \mathcal P),$ where $\mathbb V$ is a vector space and $\mathcal P$ is a multiplication, is called a terminal algebra if for all $a\in \mathbb V$ we have
\begin{equation}\label{tress}
[[[{\mathcal P},a],{\mathcal P}],{\mathcal P}]=0.
\end{equation}

\end{definition}

Note that we can expand the relation (\ref{tress}), obtaining an identity of degree 4. Therefore, the class of terminal algebras is a variety.
Also, about terminal and conservative algebras see, for example \cite{Kantor90, KPP18, cfk19, kv15, klp15, kkp19, p19}.
 
The following remark is obtained by straightforward calculations.

\begin{remark}
A commutative algebra satisfying (\ref{tress}) is a Jordan algebra. 
\end{remark}

Aside from Jordan algebras, the class of terminal algebras includes all Lie algebras, all (left) Leibniz algebras and some other types of algebras.

The following characterization of terminal algebras, proved by Kantor \cite[Theorem 2]{Kantor89}, provides a description of this class as a subclass of the class of conservative algebras.

\begin{remark}
An algebra $(\mathbb{V}, {\bf P})$ is terminal if and only if it is conservative
 and the multiplication in the associated superalgebra ${\bf P}^*$ can be defined by
 \begin{equation}
\label{terminal_associated} 
     {\bf P}^*(x,y)=\frac{2}{3}{\bf P}(x,y)+\frac{1}{3}{\bf P}(y,x).
 \end{equation}
\end{remark}

%Our method of classification of nilpotent terminal algebras is based on the calculation of central extensions of %smaller nilpotent algebras from the same variety.  The algebraic study of central extensions of Lie and non-Lie %algebras has a very big story \cite{omirov,hac16,ss78}.
%Skjelbred and Sund used central extensions of Lie algebras for a classification of nilpotent Lie algebras  %\cite{ss78}. After that, using the method described by Skjelbred and Sund were described  all non-Lie central %extensions of  all $4$-dimensional Malcev algebras \cite{hac16}, all non-associative central extensions of %$3$-dimensional Jordan algebras, all anticommutative central extensions of $3$-dimensional anticommutative %algebras,
%all central extensions of $2$-dimensional algebras.

The algebraic classification of nilpotent algebras will be obtained by the calculation of central extensions of algebras from the same variety which have a smaller dimension.
Central extensions of algebras from various varieties were studied, for example, in \cite{ss78,omirov,klp20}.
Skjelbred and Sund \cite{ss78} used central extensions of Lie algebras to classify nilpotent Lie algebras.
The same method was applied to describe,  
all non-Lie central extensions of  all $4$-dimensional Malcev algebras \cite{hac16},
all non-associative central extensions of all $3$-dimensional Jordan algebras \cite{ha17},
all anticommutative central extensions of $3$-dimensional anticommutative algebras,
all central extensions of $2$-dimensional algebras  
and some others were described.
One can also look at the classification of
$3$-dimensional nilpotent algebras \cite{fkkv},
$4$-dimensional nilpotent associative algebras \cite{degr1},
$4$-dimensional nilpotent Novikov algebras \cite{kkk18},
$4$-dimensional nilpotent bicommutative algebras \cite{kpv19},
$4$-dimensional nilpotent commutative algebras in \cite{fkkv},
$5$-dimensional nilpotent restricted Lie agebras \cite{usefi1},
$5$-dimensional nilpotent Jordan algebras \cite{ha16},
$5$-dimensional nilpotent anticommutative algebras \cite{fkkv},
$5$-dimensional nilpotent associative commutative algebras \cite{krs19},
$6$-dimensional nilpotent Lie algebras \cite{degr3, degr2},
$6$-dimensional nilpotent Malcev algebras \cite{hac18},
$6$-dimensional nilpotent Tortkara algebras \cite{gkk19},
$6$-dimensional nilpotent binary Lie algebras \cite{ack},
$6$-dimensional nilpotent anticommutative  algebras \cite{kkl20},
$8$-dimensional dual mock-Lie algebras \cite{ckls19} and so on.

\medskip 

\paragraph{\bf Motivation and contextualization} Given algebras ${\bf A}$ and ${\bf B}$ in the same variety, we write ${\bf A}\to {\bf B}$ and say that ${\bf A}$ {\it degenerates} to ${\bf B}$, or that ${\bf A}$ is a {\it deformation} of ${\bf B}$, if ${\bf B}$ is in the Zariski closure of the orbit of ${\bf A}$ (under the aforementioned base-change action of the general linear group). The study of degenerations of algebras is very  related to deformation theory, in the sense of Gerstenhaber. It offers an insightful geometric perspective on the subject and has been studied in various papers.
In particular, there are many results concerning degenerations of algebras of small dimensions in a  variety defined by a set of identities.
One of the main problems of the {\it geometric classification} of a variety of algebras is a description of its irreducible components. In the case studied in various papers finitely many orbits (i.e., isomorphism classes), the irreducible components are determined by the rigid algebras --- algebras whose orbit closure is an irreducible component of the variety under consideration. 

Observe that the algebraic classification of complex one-generated $5$-dimensional nilpotent terminal algebras presented in this paper is the first step in the full algebraic classification of complex one-generated $5$-dimensional nilpotent terminal algebras, 
which will halp to obtain the geometric classification of  all
complex $5$-dimensional nilpotent terminal algebras.

\section{Preliminaries and notations}
\subsection{Method of classification of nilpotent algebras}
Throughout this paper we  use the notation and methods described in \cite{hac16,cfk18}
and adapted for the terminal case with some modifications (see also \cite{Jac} for a discussion of extensions of algebras in an arbitrary nonassociative variety). Therefore, all statements in this subsection are given without proofs, which can be found in the papers cited above.

Let ${\bf A}$ be a terminal algebra over $\mathbb C$ and $\mathbb V$ a vector space of dimension $s$ over the same base field. Then the space ${\rm Z_T^2}(\bf A,\mathbb V )$ is defined as the set of all  maps $\theta :{\bf A} \times {\bf A} \to {\mathbb V}$ such that 
\[\theta(b,a(xy) - (ax)y - x(ay)) - \theta(a,(bx)y) + \theta(a(bx),y) + \theta(bx,ay)\]
\[- \theta(a,x(by)) + \theta(ax,by) + \theta(x,a(by)) = -\theta({\bf P}^*(a,b),xy) + \theta({\bf P}^*(a,b)x,y) + \theta(x,{\bf P}^*(a,b)y),\]
where ${\bf P}^*$ is given by (\ref{terminal_associated}). Its elements will be called \textit{cocycles}. For a linear map $f: \bf A \to \mathbb V$ define $\delta f: {\bf A} \times
{\bf A} \to {\mathbb V}$ by $\delta f (x,y ) =f(xy).$ One can check that $\delta f\in {\rm Z_T^2}({\bf A},{\mathbb V} ).$ Therefore, ${\rm B^2}({\bf A},{\mathbb V}) =\left\{ \theta =\delta f\ :f\in \operatorname{Hom}({\bf A},{\mathbb V}) \right\}$ is a subspace of ${\rm Z_T^2}({\bf A},{\mathbb V})$ whose elements are called \textit{coboundaries}. We define the \textit{second cohomology space} ${\rm H_T^2}({\bf A},{\mathbb V})$ as the quotient space ${\rm Z_T^2}({\bf A},{\mathbb V}) \big/{\rm B^2}({\bf A},{\mathbb V}).$
The equivalence class of $\theta \in {\rm Z^2}({\bf A},{\mathbb V})$  in ${\rm H^2}({\bf A},{\mathbb V})$ will be denoted by $[\theta].$

\bigskip

For a bilinear map $\theta :{\bf A} \times {\bf A} \to {\mathbb V}$ define on the linear space ${\bf A}_{\theta }:={\bf A}\oplus {\mathbb V}$ a bilinear product $[-,-]_{{\theta}}$ by 
\[[x+x^{\prime },y+y^{\prime }]_{{\theta}}= xy +\theta (x,y) \mbox{ for all } x,y\in {\bf A},x^{\prime },y^{\prime }\in {\mathbb V}.\]
Then ${\bf A}_{\theta }$ is an algebra called an $s${\it{-dimensional central extension}} of ${\bf A}$ by ${\mathbb V}.$ The following statement can be verified directly:

\begin{lemma}
The algebra ${\bf A_{\theta}}$ is terminal \textit{if and only if} $\theta \in {\rm Z_T^2}({\bf A}, {\mathbb V}).$
\end{lemma}

Recall that the {\it{annihilator}} of ${\bf A}$ is defined as the ideal $\operatorname{Ann}({\bf A} ) =\{ x\in {\bf A}:  x{\bf A}+{\bf A}x=0\}.$ Given $\theta \in {\rm Z_T^2}({\bf A}, {\mathbb V}),$ we call the set $\operatorname{Ann}(\theta)=\left\{ x\in {\bf A}:\theta (x, {\bf A} )+\theta ({\bf A},x ) = 0 \right\}$ the {\it{annihilator}} of $\theta.$ 

\begin{lemma}
\label{ann_of_ext}
$\operatorname{Ann}({\bf A}_{\theta }) = (\operatorname{Ann}(\theta) \cap \operatorname{Ann}({\bf A}))
 \oplus {\mathbb V}.$
\end{lemma}

\bigskip

Therefore, $0 \neq \mathbb{V} \subseteq \operatorname{Ann}({\bf A}_{\theta}).$ The following lemma shows that every algebra with a nonzero annihilator can be obtained in the way described above:

\begin{lemma}
Let ${\bf A}$ be an $n$-dimensional terminal algebra such that $\dim(\operatorname{Ann}({\bf A}))=m\neq0.$ Then there exists, up to isomorphism, a unique $(n-m)$-dimensional terminal algebra ${\bf A}^{\prime }$ and a bilinear map $\theta \in {\rm Z_T^2}({\bf A}, {\mathbb V})$ with $\operatorname{Ann}({\bf A}')\cap \operatorname{Ann}(\theta)=0,$ where $\mathbb V$ is a vector space of dimension $m,$ such that ${\bf A}\cong {\bf A}^{\prime }_{\theta}$ and ${\bf A}/\operatorname{Ann}({\bf A})\cong {\bf A}^{\prime }.$
\end{lemma}

In particular, any finite-dimensional nilpotent algebra is a central extension of another nilpotent algebra of strictly smaller dimension. Thus, to classify all nilpotent terminal algebras of a fixed dimension, we need to classify cocycles of nilpotent terminal algebras ${\bf A}'$ of smaller dimension (with an additional condition $\operatorname{Ann}({\bf A}')\cap \operatorname{Ann}(\theta)=0$) and central extensions that arise from them.

\bigskip

We can reduce the class of extensions that we need to consider.

\begin{definition}
Let ${\bf A}$ be an algebra and $I$ be a subspace of ${\rm Ann}({\bf A})$. If ${\bf A}={\bf A}_0 \oplus I$
then $I$ is called an {\it annihilator component} of ${\bf A}$.
 A central extension of an algebra $\bf A$ without an annihilator component is called a non-split central extension.
\end{definition} 

Clearly, we are only interested in non-split extensions (in the contrary case we can cut off annihilator components lying in $\mathbb{V}$ until we obtain a non-split extension). 

Let us fix a basis $e_{1},\ldots,e_{s}$ of ${\mathbb V},$ and $\theta \in {\rm Z_T^2}({\bf A},{\mathbb V}).$ Then $\theta$ can be uniquely written as $\theta (x,y) =\theta (x,y) =\sum_{i=1}^s\theta_{i}(x,y) e_{i},$ where $\theta_{i}\in {\rm Z_T^2}({\bf A},\mathbb C).$ Moreover, $\operatorname{Ann}(\theta)=\operatorname{Ann}(\theta_{1})\cap \operatorname{Ann}(\theta_{2})\cap\cdots \cap \operatorname{Ann}(\theta_{s}).$ Further, $\theta \in {\rm B^2}({\bf A},{\mathbb V}) $\ if and only if all $\theta_{i}\in {\rm B^2}({\bf A},\mathbb C).$ Using this presentation, one can determine whether the extension corresponding to a cocycle $\theta$ is split:

\begin{lemma} \cite[Lemma 13]{hac16}
\label{split_lindep}
Let $\theta (x,y) =\sum_{i=1}^s\theta_{i}(x,y) e_{i}\in {\rm Z_T^2}({\bf A},{\mathbb V})$ be such that $\operatorname{Ann}(\theta)\cap \operatorname{Ann}({\bf A}) = 0.$ Then ${\bf A}_{\theta }$ has an annihilator component if and only if $[\theta_{1}],[\theta_{2}],\ldots,[\theta_{s}]$ are linearly dependent in ${\rm H_T^2}({\bf A},\mathbb C).$
\end{lemma}  

\bigskip

Some cocycles give rise to isomorphic extensions:

\begin{lemma}
Let $\theta, \vartheta \in {\rm Z_T^2}({\bf A}, {\mathbb V})$ be such that $[\theta] = [\vartheta].$ Then ${\bf A}_\theta \cong {\bf A}_\vartheta.$  
\end{lemma}

%Thus, we need to compute the space ${\rm B^2}({\bf A}, \mathbb V).$ This is easily done with the aid of

%\begin{lemma}
%The map $\delta$ restricts to a linear isomorphism between $({\bf A}^2)^*$ and ${\rm B^2}({\bf A}, \mathbb C),$  where $({\bf A}^2)^* %\subseteq {\bf A}^*$ is the dual space of the square of ${\bf A}^2.$
%\end{lemma}

By above, the isomorphism classes of extensions correspond to certain equivalence classes on ${\rm H_T^2}({\bf A},{\mathbb V}).$ These classes can be given in terms of actions of certain groups on this space. In particular, let $\operatorname{Aut}({\bf A})$ be the automorphism group of ${\bf A},$ let $\phi \in \operatorname{Aut}({\bf A}),$ and let $\psi \in \operatorname{GL}(\mathbb V).$ For $\theta \in {\rm Z_T^2}({\bf A},{\mathbb V})$ define 
\[\phi \theta(x,y) =\theta (\phi(x), \phi(y)), \quad \psi\theta (x,y) = \psi(\theta(x,y)).\] 
Then $\phi\theta, \psi\theta \in {\rm Z_T^2}({\bf A},{\mathbb V}).$ Hence, $\operatorname{Aut}({\bf A})$ and $\operatorname{GL}(\mathbb V)$ act on ${\rm Z_T^2}({\bf A},{\mathbb V}).$ It is easy to verify that ${\rm B^2}({\bf A},{\mathbb V})$ is invariant under both actions. Therefore, we have induced actions on ${\rm H_T^2}({\bf A},{\mathbb V}).$

\begin{lemma}
Let $\theta, \vartheta \in {\rm Z_T^2}({\bf A},{\mathbb V})$ be such that $\operatorname{Ann}({\bf A}_\theta) = \operatorname{Ann}({\bf A}_\vartheta) = \mathbb{V}.$ Then ${\bf A}_\theta \cong {\bf A}_\vartheta$ if and only if there exist a $\phi \in \operatorname{Aut}({\bf A}), \psi \in \operatorname{GL}(\mathbb V)$ such that $[\phi\theta] = [\psi\vartheta].$
\end{lemma}

\bigskip

Now we rewrite the above lemma in a form more suitable for computations.

Let ${\mathbb U}$ be a finite-dimensional vector space over $\mathbb C.$ The {\it{Grassmannian}} $G_{k}({\mathbb U})$ is the set of all $k$-dimensional linear subspaces of ${\mathbb V}.$ Let $G_{s}({\rm H_T^2}({\bf A},\mathbb C) )$ be the Grassmannian of subspaces of dimension $s$ in ${\rm H_T^2}({\bf A},\mathbb C).$ There is a natural action of $\operatorname{Aut} ({\bf A})$ on $G_{s}({\rm H_T^2}({\bf A},\mathbb C))$ : for $\phi \in \operatorname{Aut} ({\bf A}), W=\langle [\theta_{1}],[\theta_{2}],\dots,[\theta_{s}] \rangle \in G_{s}({\rm H_T^2}({\bf A},\mathbb C)),$ define 
\[\phi W=\langle [\phi \theta_{1}], [\phi \theta_{2}],\dots,[\phi \theta_{s}]\rangle.\]
Note that this action is compatible with the action of $\operatorname{Aut} ({\bf A})$ on ${\rm H_T^2}({\bf A},{\mathbb V})$ and the above presentation of a cocycle as a collection of $s$ elements of ${\rm H_T^2}({\bf A},\mathbb C).$ Denote the orbit of $W$ under the action of $\operatorname{Aut} ({\bf A})$ by $\mathrm{Orb}(W).$ It is easy to check that given two bases of a subspace 
\begin{equation*}
W=\langle [\theta_{1}],[\theta_{2}],\dots,[\theta_{s}] \rangle =\langle [\vartheta_{1}],[\vartheta_{2}],\dots,[\vartheta_{s}]
\rangle \in G_{s}({\rm H^2}({\bf A},\mathbb C)),
\end{equation*}
we have $\cap_{i=1}^s\operatorname{Ann}(\theta_{i})\cap \operatorname{Ann}({\bf A}) =\cap_{i=1}^s\operatorname{Ann}(\vartheta_{i})\cap \operatorname{Ann}({\bf A}).$ Therefore, we can introduce the set

\begin{equation*}
T_{s}({\bf A}) =\left\{ W=\left\langle [\theta_{1}],%
[\theta_{2}],\dots,[\theta_{s}] \right\rangle \in
G_{s}({\rm H^2}({\bf A},\mathbb C) ) : \cap_{i=1}^s \operatorname{Ann}(\theta_{i})\cap \operatorname{Ann}({\bf A}) =0\right\},
\end{equation*}
which is stable under the action of $\operatorname{Aut}({\bf A}).$

\medskip

Let us denote by $E({\bf A},{\mathbb V})$ the set of all {\it non-split} central extensions of ${\bf A}$ by ${\mathbb V}.$ By Lemmas \ref{ann_of_ext} and \ref{split_lindep}, we can write
\begin{equation*}
E({\bf A},{\mathbb V}) =\left\{ {\bf A}_{\theta }:\theta (x,y) =\sum_{i=1}^s\theta_{i}(x,y) e_{i}\mbox{ and }\langle [\theta_{1}],[\theta_{2}],\dots,
[\theta_{s}] \rangle \in T_{s}({\bf A})\right\}.
\end{equation*}
Also, we have the next result, which can be proved in the same manner as \cite[Lemma 17]{hac16}.

\begin{lemma}
 Let $\theta(x,y) =\sum_{i=1}^s \theta_{i}(x,y) e_{i}$ and $\vartheta(x,y) = \sum_{i=1}^s\vartheta_{i}(x,y) e_{i}$ be such that ${\bf A}_{\theta},{\bf A}_{\vartheta }\in E({\bf A},{\mathbb V}).$  Then the algebras ${\bf A}_{\theta }$ and ${\bf A}_{\vartheta }$ are isomorphic if and only if 
\[\mathrm{Orb}\langle [\theta_{1}], [\theta_{2}],\dots,[\theta_{s}] \rangle = \mathrm{Orb}\langle [\vartheta_{1}],[\vartheta
_{2}],\dots,[\vartheta_{s}] \rangle. \]
\end{lemma}

Thus, there exists a one-to-one correspondence between the set of $\operatorname{Aut}({\bf A}) $-orbits on $T_{s}({\bf A})$ and the set of isomorphism classes of $E({\bf A},{\mathbb V}).$ Consequently, we have a procedure that allows us, given a terminal algebra ${\bf A}^{\prime }$ of dimension $n,$ to construct all non-split central extensions of ${\bf A}^{\prime }.$ This procedure is as follows:

\medskip

{\centerline{\it Procedure}}

\begin{enumerate}
\item For a given (nilpotent) terminal algebra $\bf{A}^{\prime }$
of dimension $n-s,$ determine ${T}_{s}(\bf{A}^{\prime })$
and $\operatorname{Aut}(\bf{A}^{\prime }).$

\item Determine the set of $\operatorname{Aut}(\bf{A}^{\prime })$-orbits on $%
{T}_{s}(\bf{A}^{\prime }).$

\item For each orbit, construct the terminal algebra corresponding to one of its
representatives.
\end{enumerate}

\subsection{Notations}
Let us introduce the following notations. Let ${\bf A}$ be a terminal algebra with
a basis $e_{1},e_{2}, \ldots, e_{n}.$ Then by $\Delta_{ij}$\ we will denote the
bilinear form
$\Delta_{ij}:{\bf A}\times {\bf A}\longrightarrow \mathbb C$
with $\Delta_{ij}(e_{l},e_{m}) = \delta_{il}\delta_{jm}.$
The set $\left\{ \Delta_{ij}:1\leq i, j\leq n\right\}$ is a basis for the linear space of 
bilinear forms on ${\bf A},$ so every $\theta \in
{\rm Z^2}({\bf A},\bf \mathbb V )$ can be uniquely written as $
\theta = \displaystyle \sum_{1\leq i,j\leq n} c_{ij}\Delta _{{i}{j}}$, where $
c_{ij}\in \mathbb C$.
Let us fix the following notations:
$$\begin{array}{lll}
\T {i}{j}& \mbox{---}& j\mbox{th }i\mbox{-dimensional one-generated terminal algebra.} \\
\end{array}$$

\subsection{The algebraic classification of $3$-dimensional one-generated nilpotent terminal algebras}

Observe that a 2-dimensional one-generated nilpotent algebra is isomorphic to
\begin{align*}%\label{alg-A_3}
    \begin{array}{llllll}
        \T {2}{01} &:& e_1 e_1 = e_2.
    \end{array}
\end{align*}
We have the following list of all 3-dimensional one-generated nilpotent terminal algebras with an annihilator of codimension 1 or 2:
\begin{align}\label{3-dim-term}
    \begin{array}{lllllllll}
        \T {3}{01}&:& e_1 e_1 = e_2 & e_1 e_2=e_3 \\
        \T 3{02}(\lb)&:& e_1 e_1 = e_2 & e_1 e_2=\lb e_3 & e_2 e_1=e_3.
    \end{array}
\end{align}

\subsection{Automorphism and cohomology groups of $3$-dimensional one-generated nilpotent terminal algebras}\label{aut-and-H^2} 
In the following table we give the description of the automorphism groups and the second cohomology spaces of  
$3$-dimensional one-generated nilpotent terminal algebras.

\begin{longtable}{|l|c|c|c|}
\hline
$\A$ & $\aut\A$ & ${\rm Z_T^2}(\A)$ & ${\rm H_T^2}(\A)$\\
\hline
$\T {3}{01}$
&
$\begin{pmatrix}
x &    0  &  0\\
y &  x^2  &  0\\
z &   xy  &  x^3
\end{pmatrix}$
&
$\Big\langle
\begin{array}{l}
     \Dt 11, \Dt 12, \Dt 13,\\
     \Dt 21, \Dt 22 - 3\Dt 31
\end{array}
\Big\rangle$
&
$\Big\langle
\begin{array}{l}
     \Dl 13, \Dl 21,\\
     \Dl 22 - 3\Dl 31
\end{array}
\Big\rangle$
\\
\hline
$\T 3{02}(\lambda)_{\lambda \neq 0}$
&
$\begin{pmatrix}
x &               0  &  0\\
y &             x^2  &  0\\
z &   (\lambda+1)xy  &  x^3
\end{pmatrix}$
&
$
\Big\langle
\begin{array}{l}
     \Dt 11, \Dt 12, \Dt 13 + \Dt 22,\\
     \Dt 21, \frac{1-\lb} 3 \Dt 22 + \Dt 31
\end{array} \Big\rangle
$
&
$
\Big\langle
\begin{array}{l}
     \Dl 12, \Dl 13 + \Dl 22,\\
     \frac{1-\lb} 3\Dl 22 + \Dl 31
\end{array}
\Big\rangle $
\\
\hline
$\T 3{02}(0)$
&
$\begin{pmatrix}
x &               0  &  0\\
y &             x^2  &  0\\
z &   xy  &  x^3
\end{pmatrix}$
&
$
\Big\langle
\begin{array}{l}
     \Dt 11, \Dt 12, \Dt 13 + \Dt 22,\\
     \Dt 21, \frac 1 3 \Dt 22 + \Dt 31, \Dt 23
\end{array}
\Big\rangle
$
&
$
\Big\langle
\begin{array}{l}
     \Dl 12, \Dl 13 + \Dl 22,\\
     \frac 1 3 \Dl 22 + \Dl 31, \Dl 23
\end{array}
\Big\rangle $
\\
\hline

\end{longtable}

\vskip0.5cm

\section{Classification of $4$-dimensional one-generated nilpotent terminal algebras}

Thanks to \cite{cfk18}, we have only one $4$-dimensional one-generated nilpotent terminal algebra with $2$-dimensional annihilator:
\begin{longtable}{lllllllllllll}
    $\T 4{01}$&:& $e_1 e_1 = e_2$ & $e_1 e_2=e_4$ & $e_2 e_1=e_3.$\end{longtable}
    
\subsection{Central extensions of $\T {3}{01}$}\label{ext-T_01^3}
	Let us use the following notations:
	\begin{align*}
	\nb 1 = \Dl 13, \quad \nb 2 = \Dl 21, \quad \nb 3 = \Dl 22 - 3\Dl 31.
	\end{align*}
	Take $\0=\sum\limits_{i=1}^3\af_i\nb i\in {\rm H_T^2}(\T {3}{01}).$
	For an element
	$$
	\phi=
	\begin{pmatrix}
	x &    0  &  0\\
	y &  x^2  &  0\\
	z &   xy  &  x^3
	\end{pmatrix}\in\aut{\T {3}{01}},
	$$
	we have 
	$$
	\phi^T\begin{pmatrix}
	0      &  0    & \af_1\\
	\af_2  & \af_3 & 0\\
	-3\af_3&  0    & 0
	\end{pmatrix} \phi=
	\begin{pmatrix}
	\af^*      &  \af^{**}    & \af^*_1\\
	\af^*_2    & \af^*_3      & 0\\
	-3\af^*_3  &  0           & 0
	\end{pmatrix},
	$$
	where
\[
\begin{array}{rclcrclcrcl}
	\af^*_1 &=& x^4\af_1,&&
	\af^*_2 &=& x^2(x\af_2 - 2y\af_3),&&
	\af^*_3 &=& x^4 \af_3.
	\end{array} \]

\subsubsection{$1$-dimensional central extensions of $\T {3}{01}$}

Since we are interested only in new algebras, we have the following cases:
\begin{enumerate}
\item if $\alpha_3\neq 0,$ then by choosing  $y=\frac {x\alpha_2}{2\af_3}$, we have the representative
$\langle \alpha \nabla_1+\nabla_3\rangle.$ 
\item if $\alpha_3= 0,$ then $\alpha_1\neq 0,$ and we have the representative $\langle  \nabla_1 \rangle$ and 
$\langle  \nabla_1+ \nabla_2\rangle$ depending on whether $\af_2=0$ or not.

\end{enumerate}

Hence, we have  the following $4$-dimensional algebras:
\begin{longtable}{lllllllllllllllllll}
\hline    $\T{4}{02}(\alpha)$&:& $e_1e_1 = e_2$& $e_1e_2 = e_3$& $e_1e_3 = \alpha e_4$& $e_2e_2 = e_4$& $e_3e_1 = -3e_4$\\
\hline    $\T 4{03}$  &: &  $e_1e_1 = e_2$ & $e_1e_2 = e_3$ &  $e_1e_3 = e_4$ \\
\hline    $\T{4}{04}$&:& $e_1e_1 = e_2$ & $e_1e_2 = e_3$& $e_1e_3 = e_4$& $e_2e_1 = e_4.$ \\
\hline \end{longtable}

%$$ \T 4{02}(\alpha), \quad \T 4{03}, \quad  \T 4{04}.$$

\subsubsection{$2$-dimensional central extensions of $\T {3}{01}$}
Consider the vector space generated by the following two cocycles:
$$\begin{array}{rcl}
\theta_1 &=& \alpha_1 \nabla_1+\alpha_2\nabla_2+\alpha_3\nabla_3  \\
\theta_2 &=& \beta_1 \nabla_1+\beta_2\nabla_2.
\end{array}$$

If $\alpha_3=0,$ then we have the representative $\langle \nabla_1, \nabla_2 \rangle.$
If $\alpha_3\neq 0,$ then:
\begin{enumerate}
    \item if $\beta_1=0,$ then we have the family of representatives
    $\langle \nabla_2, \alpha\nabla_1+ \nabla_3 \rangle.$
    \item if $\beta_1\neq 0,$  then we have the representatives
    $\langle \nabla_1, \nabla_3 \rangle,$ and 
    $\langle \nabla_1+\nabla_2, \nabla_3 \rangle$ depending on whether $\beta_2=0$ or not.
\end{enumerate}

Hence, we have the following $5$-dimensional algebras:
\begin{longtable}{llllllllllllllllllllllll}
\hline     $\T{5}{01} $ &  $e_1 e_1 = e_2$ & $e_1 e_2=e_3$ & $e_1 e_3=e_4$ & $e_2 e_1=e_5$ &\\
\hline     $\T{5}{02}(\alpha) $ &  $e_1 e_1 = e_2$ & $e_1 e_2=e_3$ & $e_1 e_3=\alpha e_5$    & $e_2 e_1=e_4$ & $e_2 e_2=e_5$  & $e_3 e_1=-3e_5$ &&\\
\hline     $\T{5}{03} $ &  $e_1 e_1 = e_2$ & $e_1 e_2=e_3$ & $e_1 e_3=e_4$ &  $e_2 e_2=e_5$ & $e_3 e_1=-3e_5$\\
\hline     $\T{5}{04} $ &  $e_1 e_1 = e_2$ & $e_1 e_2=e_3$ & $e_1 e_3=e_4$ & $e_2 e_1=e_4$ & $e_2 e_2=e_5$   & $e_3 e_1=-3e_5.$ &&\\
\hline \end{longtable}

%\[\T 5{01}, \quad \T 5{02}(\alpha), \quad  \T 5{03}, \quad   \T 5{04}\]

\subsection{Central extensions of $\T {3}{02}(\lambda)_{\lambda \neq 0}$}

Let us use the following notations:
	\begin{align*}
	\nb 1 = \Dl 12, \quad  \nb 2 = \Dl 13 + \Dl 22, \quad \nb 3 = \frac{1-\lb} 3\Dl 22 +  \Dl 31.
	\end{align*}
	Take $\0=\sum\limits_{i=1}^3\af_i\nb i\in {\rm H_T^2}(\T {3}{02}(\lambda)).$
	For an element
	$$
	\phi=
	\begin{pmatrix}
x &               0  &  0\\
y &             x^2  &  0\\
z &   (\lambda+1)xy  &  x^3
\end{pmatrix}
\in\aut{\T {3}{01}},
	$$
	 we have
	$$
	\phi^T\begin{pmatrix}
	0      &  \af_1    & \af_2\\
	0  & \af_2+ \frac{1-\lb} 3\af_3 & 0\\
	\af_3&  0    & 0
	\end{pmatrix} \phi=
	\begin{pmatrix}
	\alpha^*      &  \af_1^* + \lambda\alpha^{**}    & \af_2^*\\
	\alpha^{**}  & \af_2^*+\frac{1-\lb} 3 \af_3^* & 0\\
	\af_3^*&  0    & 0
	\end{pmatrix},
	$$
	where
\[
\begin{array}{rclcrclcrcl}
	\af^*_1 &=& x^3\af_1+x^2y\Big( 2\af_2 - \frac{2\lb^2  + 5\lb  - 1} 3 \af_3\Big) , &&
	\af^*_2 &=& x^4\af_2,&&
	\af^*_3 &=& x^4\af_3.
	\end{array} \]

\subsubsection{$1$-dimensional central extensions of $\T {3}{02}(\lambda)_{\lambda \neq 0}$}
	
Since we are interested only in new algebras, we have the following cases:
\begin{enumerate}
\item if $\alpha_3\neq 0,$ then:
\begin{enumerate}
\item if $\alpha_2\neq \frac{2\lb^2  + 5\lb  - 1} 6 \af_3,$ then by choosing  $y=\frac {3x\alpha_1}{6\af_2 - (2\lb^2  + 5\lb  - 1) \af_3}$, we have the representative
$\langle \alpha \nabla_2 + \nabla_3\rangle,$ with $\alpha \neq \frac{2\lb^2  + 5\lb  - 1} 6.$

\item if $\alpha_2= \frac{2\lb^2  + 5\lb  - 1} 6 \af_3,$ then we have the representatives
\[\langle  \frac{2\lb^2  + 5\lb  - 1} 6 \nabla_2 + \nabla_3\rangle\mbox{ and }\langle \nabla_1+ \frac{2\lb^2  + 5\lb  - 1} 6 \nabla_2 + \nabla_3\rangle,\] depending on whether $\alpha_1=0$ or not. 
\end{enumerate}

\item if $\alpha_3= 0,$ then $\alpha_2\neq 0,$ and by choosing $y=- \frac {x\af_1}{2\af_2},$ we have the representative $\langle  \nabla_2 \rangle.$ 

\end{enumerate}

	Hence, we have the following distinct orbits:
	\[\langle \alpha \nabla_2 + \nabla_3\rangle, \quad \langle \nabla_1+ \frac{2\lb^2  + 5\lb  - 1} 6 \nabla_2 + \nabla_3\rangle ,\quad \langle  \nabla_2 \rangle \]
	
\subsubsection{$2$-dimensional central extensions of $\T {3}{02}(\lambda)_{\lambda \neq 0}$}
	Consider the vector space generated by the following two cocycles:
$$\begin{array}{rcl}
\theta_1 &=& \alpha_1 \nabla_1+\alpha_2\nabla_2+\alpha_3\nabla_3  \\
\theta_2 &=& \beta_1 \nabla_1+\beta_2\nabla_2.
\end{array}$$

If $\alpha_3=0,$ then we have the representative $\langle \nabla_1, \nabla_2 \rangle.$
If $\alpha_3\neq 0,$ then:

\begin{enumerate}
    \item if $\beta_2=0,$ then we have the family of representatives
    $\langle \nabla_1, \alpha\nabla_2+ \nabla_3 \rangle.$

    \item if $\beta_2\neq 0,$  then choosing $y=-\frac{x\beta_1}{2\beta_2}$, we have the representatives
                $\langle \nabla_2, \nabla_3 \rangle$ and $\langle \nabla_2, \nabla_1+\nabla_3 \rangle$ 
                depending on whether $\alpha_1^*=0$ or not.
    
\end{enumerate}

	Hence, we have the following distinct orbits:
\[ \langle \nabla_1, \nabla_2 \rangle, \quad  \langle \nabla_1, \alpha\nabla_2+ \nabla_3 \rangle, \quad \langle \nabla_2, \nabla_3 \rangle,   \quad  \langle \nabla_2, \nabla_1+\nabla_3 \rangle.\]

\subsection{Central extensions of $\T {3}{02}(0)$}
	Let us use the following notations:
	\begin{align*}
	\nb 1 = \Dl 12, \quad \nb 2 = \Dl 13 + \Dl 22, \quad \nb 3 = \frac 1 3 \Dl 22 + \Dl 31, \quad \nb 4 = \Dl 23.
	\end{align*}
	Take $\0=\sum\limits_{i=1}^4\af_i\nb i\in {\rm H_T^2}(\T {3}{01}(0)).$
For an element
	$$
	\phi=
	\begin{pmatrix}
    x &               0  &  0\\
    y &             x^2  &  0\\
    z &   xy  &  x^3
    \end{pmatrix}
\in\aut{\T {3}{01}(0)},
	$$
	we have
	$$
	\phi^T\begin{pmatrix}
	0      &  \af_1    & \af_2\\
	0  & \af_2+\frac 1 3\af_3 & \af_4\\
	\af_3&  0    & 0
	\end{pmatrix} \phi=
	\begin{pmatrix}
	\alpha^*      &  \af^*_1    & \af^*_2\\
	\alpha^{**}  & \af^*_2+\frac 1 3 \af^*_3 & \af^*_4\\
	\af^*_3&  0    & 0
	\end{pmatrix},
	$$
	where
\[
\begin{array}{rclcrclcrcl}
	\af^*_1 &=& x^3 \af_1 + x^2 y (2\af_2 + \frac 1 3 \af_3) + x y^2 \af_4, &&
    \af^*_2 &=& x^4 \af_2 + x^3 y \af_4,\\
    \af^*_3 &=& x^4\af_3, &&
	\af^*_4 &=& x^5 \af_4.
	\end{array} \]
	
\subsubsection{$1$-dimensional central extensions of $\T {3}{02}(0)$}

	\begin{enumerate}
\item if $\alpha_4 =  0,$ then:
\begin{enumerate}
\item if $\alpha_3\neq 0, $ $\alpha_2\neq - \frac {1} 6 \alpha_3,$
then by choosing  $y=\frac {3x\alpha_1}{6\af_2 + \af_3}$, we have the representative
$\langle \alpha \nabla_2 + \nabla_3\rangle,$ with $\alpha \neq - \frac {1} 6.$ 

\item if $\alpha_3\neq 0, $ $\alpha_2= - \frac {1} 6 \alpha_3,$
then we have the representatives
$\langle  -\frac{1} 6 \nabla_2 + \nabla_3\rangle$ and $\langle \nabla_1 - \frac{1} 6 \nabla_2 + \nabla_3\rangle,$ depending on whether $\alpha_1=0$ or not. 

\item if $\alpha_3= 0,$ then $\alpha_2\neq 0,$ and by choosing $y= - \frac {x\af_1}{2\af_2},$ we have the representative $\langle  \nabla_2 \rangle.$ 

\end{enumerate}

\item if $\alpha_4 \neq 0,$ then by choosing $y= - \frac {x\af_2}{\af_4},$ 
we have the representatives 
\[\langle  \nabla_4 \rangle, \ \langle  \nabla_1 + \nabla_4 \rangle \mbox{ 
and }\langle \alpha \nabla_1 + \nabla_3 +\nabla_4 \rangle\] 
depending on whether $\alpha_3=0$ and 
$3\alpha_1\alpha_4 - 3\alpha_2^2 - \alpha_2\alpha_3=0$ or not.

\end{enumerate}

	Hence, we have the following distinct orbits:
	\[\langle \alpha \nabla_2 + \nabla_3\rangle, \quad \langle \nabla_1- \frac{ 1} 6 \nabla_2 + \nabla_3\rangle ,\quad \langle  \nabla_2 \rangle,  \quad \langle  \nabla_4 \rangle, \quad \langle  \nabla_1 + \nabla_4 \rangle, \quad
\langle \alpha \nabla_1 + \nabla_3 +\nabla_4 \rangle.\]

Summarizing the cases $\lambda\neq0$ and $\lambda=0,$ we have  the following new $4$-dimensional algebras:

\begin{longtable}{llllllllllllllllllllll}
\hline    $\T{4}{05}(\lambda, \alpha)$ & 
$e_1 e_1 = e_2$ & $e_1 e_2=\lambda e_3$ & $e_1 e_3=\alpha e_4$ \\
& $e_2 e_1=e_3$ & $e_2 e_2=\big(\alpha+\frac{1-\lambda}3\big) e_4$ & $e_3 e_1= e_4$\\
\hline    $\T{4}{06}(\lambda)$ & $e_1 e_1 = e_2$ & $e_1 e_2=\lambda e_3+e_4$ & \multicolumn{2}{l}{$e_1 e_3= \frac{2\lb  + 5\lb  - 1}{6}e_4$} \\ & $e_2 e_1=e_3$ & $e_2 e_2= \frac{(2\lb  + 1)(\lb  + 1)}{6}e_4$ & $e_3 e_1= e_4$\\
\hline    $\T{4}{07}(\lambda)$ & $e_1 e_1 = e_2$ & $e_1 e_2=\lambda e_3$ & $e_1 e_3= e_4$ & $e_2 e_1=e_3$ & $e_2 e_2= e_4$\\
\hline    $\T 4{08}$   & $e_1 e_1 = e_2$ & $e_2 e_1=e_3$    & $e_2e_3 = e_4$\\
\hline    $\T 4{09}$   & $e_1 e_1 = e_2$ & $e_1 e_2=e_4$ & $e_2 e_1=e_3$ & $e_2e_3 = e_4$\\
\hline    $\T 4{10}(\af)$   & $e_1 e_1 = e_2$ & $e_1 e_2=\af e_4$ &  $e_2 e_1=e_3$ & $e_2e_2=\frac 1 3 e_4$ & $e_2e_3 = e_4$ & $e_3e_1=e_4$\\
\hline    \end{longtable}

%$$ \T 4{05}(\alpha, \lambda), \quad \T 4{06}(\lambda), \quad \T 4{07}(\lambda), \quad  \T 4{08}, \quad \T 4{09}, \quad  \T 4{10}(\alpha).$$

\subsubsection{$2$-dimensional central extensions of $\T {3}{02}(0)$}
Consider the vector space generated by the following two cocycles
$$\begin{array}{rcl}
\theta_1 &=& \alpha_1 \nabla_1+\alpha_2\nabla_2+\alpha_3\nabla_3 + \alpha_4\nabla_4 \\
\theta_2 &=& \beta_1 \nabla_1+\beta_2\nabla_2 +\beta_3\nabla_3.
\end{array}$$
	
	Consider the following cases:
    \begin{enumerate}
    \item if $\alpha_4=0,$ then similarly to the case of $\lambda \neq 0$, we have the representatives
    \[ \langle \nabla_1, \nabla_2 \rangle, \ \langle \nabla_1, \alpha\nabla_2+ \nabla_3 \rangle, \
    \langle \nabla_2, \nabla_3 \rangle \mbox{ and }\langle \nabla_2, \nabla_1 + \nabla_3 \rangle.\]

    \item if $\alpha_4\neq 0,$ then:
    \begin{enumerate}
		\item if $\beta_3=\beta_2=0,$ then choosing $y=- \frac{x\alpha_2}{\alpha_4},$ we have the representatives
            $\langle \nabla_1, \nabla_4 \rangle,$ $\langle \nabla_1, 3 \nabla_3 +\nabla_4 \rangle$  depending on whether  $\alpha_3=0$  or not.
        \item  $\beta_3=0, \beta_2\neq0,$ 
        then choosing $y=-\frac{x\beta_1}{2\beta_2}$, we have the representatives
        \[\langle \nabla_2, \nabla_4 \rangle, \ \langle \nabla_2, \nabla_1 + \nabla_4 \rangle \mbox{  and }\langle \nabla_2, 3 \nabla_3 +\nabla_4 \rangle\] 
        %and \langle \nabla_2, \alpha \nabla_1 +\nabla_3 +\nabla_4 \rangle,$ 
        depending on whether $\alpha_1^*=0,$ $\alpha_3^*=0$ or not.
        \item  $\beta_3\neq 0,$ then we may suppose $\alpha_3=0$ and choosing 
        $y=- \frac{x\alpha_2}{\alpha_4},$ we have the representatives 
        \[\langle \alpha\nabla_2+3 \nabla_3, \nabla_4 \rangle, \ 
        \langle \nabla_1 + \alpha\nabla_2+3 \nabla_3, \nabla_4 \rangle, \ \langle\beta\nabla_1+ \alpha\nabla_2+3 \nabla_3, \nabla_1 + \nabla_4 \rangle.\]

    \end{enumerate}

\end{enumerate}

Summarizing the cases $\lambda\neq0$ and $\lambda=0,$ we have  the following  $5$-dimensional algebras:

\begin{longtable}{llllllllllllllllll}
\hline     $\T{5}{05}(\lambda) $ &  $e_1 e_1 = e_2$ & $e_1 e_2=  \lambda e_3+e_4$ & $e_1 e_3= e_5$ & $e_2 e_1=e_3$  & $e_2 e_2= e_5$\\
\hline     $\T{5}{06}(\lambda,\alpha) $ &  $e_1 e_1 = e_2$& $e_1 e_2= \lambda e_3+e_4$ & $e_1 e_3= \alpha e_5$ \\ & $e_2 e_1=e_3$ & $e_2 e_2= \left(\alpha+\frac{1-\lambda} 3 \right) e_5$& $e_3 e_1=  e_5$ &&\\
\hline     $\T{5}{07}(\lambda) $ &  $e_1 e_1 = e_2$ & $e_1 e_2= \lambda e_3$ & $e_1 e_3=  e_4$ \\ & $e_2 e_1=e_3$ & $e_2 e_2=  e_4 + \frac{1-\lambda} 3e_5$& $e_3 e_1=  e_5$ &&\\
\hline     $\T{5}{08}(\lambda) $ &  $e_1 e_1 = e_2$& $e_1 e_2= \lambda e_3+e_5$ & $e_1 e_3=  e_4$ \\ & $e_2 e_1=e_3$ & $e_2 e_2=  e_4 + \frac{1-\lambda} {3} e_5$& $e_3 e_1=  e_5$ &&\\
\hline     $\T{5}{09} $ &  $e_1 e_1 = e_2$ & $e_1 e_2= e_4$ & $e_2 e_1=e_3$ & $e_2 e_3=e_5$ &\\
\hline     $\T{5}{10} $ &  $e_1 e_1 = e_2$ & $e_1 e_2= e_4$ & $e_2 e_1=e_3$ & $e_2 e_2= e_5$ \\ & $e_2 e_3=e_5$ & $e_3 e_1=3e_5$ &&\\
\hline     $\T{5}{11} $ &  $e_1 e_1 = e_2$ & $e_1 e_3= e_4$ & $e_2 e_1=e_3$ & $e_2 e_2=e_4$  & $e_2 e_3=e_5$\\
\hline     $\T{5}{12} $ &  $e_1 e_1 = e_2$ & $e_1 e_2= e_5$ & $e_1 e_3=e_4$ & $e_2 e_1=e_3$ \\ & $e_2 e_2=e_4$  & $e_2 e_3=e_5$ &&\\
\hline     $\T{5}{13}(\alpha) $ &  $e_1 e_1 = e_2$ & $e_1 e_2=  \alpha e_5$ & $e_1 e_3= e_4$ & $e_2 e_1=e_3$ \\ & $e_2 e_2=e_4 + e_5$& $e_2 e_3=e_5$& $e_3 e_1=3e_4$&&  \\
\hline     $\T{5}{14}(\alpha) $ &  $e_1 e_1 = e_2$ & $e_1 e_3= \alpha e_4$ & $e_2 e_1=e_3$\\ & $e_2 e_2= (\alpha +1)e_4$ & $e_2 e_3=e_5$& $e_3 e_1=3e_4$ && \\
\hline     $\T{5}{15}(\alpha) $ &  $e_1 e_1 = e_2$ & $e_1 e_2=  e_4$  & $e_1 e_3= \alpha e_4$  & $e_2 e_1=e_3$ \\& $e_2 e_2= (\alpha+1)e_4$ & $e_2 e_3=e_5$& $e_3 e_1=3 e_4$ &&\\
\hline     $\T{5}{16}(\alpha, \beta) $ &  $e_1 e_1 = e_2$& $e_1 e_2= \beta e_4 + e_5$ & $e_1 e_3= \alpha e_4$  \\& $e_2 e_1=e_3$ &      $e_2 e_2= (\alpha+1)e_4$ & $e_2 e_3=e_5$& $e_3 e_1=3e_4$ && \\
\hline     \end{longtable}

%$$ \T 5{05}, \dots, \T 5{16}.$$

\subsection{Classification theorem}
Summarizing results of the present sections,
 we have the following theorem.

\begin{theoremA}\label{main-alg}
Let $\mathbf A$ be a complex $4$-dimensional one-generated nilpotent terminal algebra. 
Then $\mathbf A$ is isomorphic to one of the algebras $\T 4{01}-\T 4{10}$  found in Table A (See, Appendix).
\end{theoremA}

\section{The algebraic classification of $5$-dimensional one-generated nilpotent terminal algebras}

To complete the algebraic classification of $5$-dimensional one-generated nilpotent terminal algebras, we need to describe all 1-dimensional extensions of the algebras from the table \ref{3-dim-term}. This is done in Subsections~\ref{ext-T_01^4}--\ref{ext-T_08^4} and summarized in the main theorem of the paper.

\subsection{$1$-dimensional  central extensions of $\T {4}{01}$}\label{ext-T_01^4}
	Let us use the following notations:
	\begin{align*}
	\nb 1 = \Dl 13 + \Dl 22, \quad \nb 2 = \Dl 14, \quad \nb 3 = \Dl 22 + 3\Dl 31, \quad \nb 4 = \Dl 23, \quad \nb 5 = \Dl 31+\Dl 41.
	\end{align*}
	
Take $\0=\sum\limits_{i=1}^5\af_i\nb i\in {\rm H_T^2}(\T {4}{01}).$
For an element
	$$
	\phi=
	\begin{pmatrix}
	x &    0  &  0&  0\\
	y &  x^2  &  0&  0\\
	z &   xy  &  x^3 &  0\\
    t &   xy  &  0 &  x^3
	\end{pmatrix}\in\aut{\T {4}{01}},
	$$
	we have
	$$
	\phi^T\begin{pmatrix}
	0      &  0     & \af_1  & \af_2\\
	0      & \af_1+\af_3  & \af_4  & 0    \\
    3\af_3+\af_5  & 0  & 0  & 0    \\
    \af_5      & 0  & 0  & 0
	\end{pmatrix} \phi=
	\begin{pmatrix}
    \af^*      &  \af^{**}     & \af^*_1  & \af^*_2\\
	\af^{***}      & \af^*_1 +\af^*_3  & \af^*_4  & 0    \\
    3\af^*_3 + \af^*_5  & 0  & 0  & 0    \\
    \af^*_5     & 0  & 0  & 0
	\end{pmatrix},
	$$
where
\[
\begin{array}{rclcrclcrcl}
	\af^*_1 &=& x^4\af_1+x^3y\af_4,&&
	\af^*_2 &=& x^4\af_2,&&
	\af^*_3 &=& x^4\af_3,\\
    \af^*_4 &=& x^5\af_4,&&
    \af^*_5 &=& x^4\af_5.
	\end{array} \]

%Hence, $\phi\langle\0\rangle=\langle\0^*\rangle,$ where $\0^*=\sum\limits_{i=1}^5 \af_i^*  \nb i.$
	
Since we are interested only in new algebras, we have $(\af_2, \af_5) \neq (0,0)$ and the following cases:
			\begin{enumerate}
			\item if $\af_4\ne 0$, then choosing $y=-\frac{x \af_1}{\af_4}$, we have the family of representatives 
			$\la\af_2\nb 2+\af_3\nb 3+\af_4\nb 4+\af_5\nb 5\ra.$ It two three distinct orbits with representatives:
\[\la\nb 2+\af\nb 3+\nb 4+\beta\nb 5\ra \mbox{ and  }\la\alpha \nb 3+\nb 4+\nb 5\ra.\]
			\item if $\af_4= 0$, then we have the family of representatives $\la\af_1\nb 1+\af_2\nb 2+\af_3\nb 3+\af_5\nb 5\ra.$
It gives the following representatives of distinct orbits:
\[ \la \af \nb 1 + \beta \nb 2 + \gamma \nb 3+\nb 5\ra_{(\alpha; \gamma)\neq(0 ; - \frac 1 3 )}
\mbox{ and }\la \af\nb 1+\nb 2 + \beta \nb 3\ra_{(\alpha; \beta)\neq(0 ; 0 )}.\]
		\end{enumerate}

Hence, we have the following $5$-dimensional algebras:
\begin{longtable}{lllllllllllll}
\hline    $\T{5}{17}(\alpha, \beta) $ &  $e_1 e_1 = e_2$ & $e_1 e_2=e_4$ & $e_1 e_4=e_5$ & $e_2 e_1=e_3$ &\\
        & $e_2 e_2= \alpha e_5$ & $e_2 e_3 = e_5$ & $e_3 e_1=(3\alpha+\beta) e_5$ & $e_4 e_1= \beta e_5$ &\\
     \hline$\T{5}{18}(\alpha) $ &  $e_1 e_1 = e_2$ & $e_1 e_2=e_4$  & $e_2 e_1=e_3$ & $e_2 e_2= \alpha e_5$ & \\
        & $e_2 e_3 = e_5$  &  $e_3 e_1=(3\alpha+1) e_5$ & $e_4 e_1=  e_5$ && \\
     \hline$\T{5}{19}(\alpha, \beta, \gamma)_{(\alpha; \gamma)\neq(0 ; - \frac 1 3 )} $ &  $e_1 e_1 = e_2$ & $e_1 e_2=e_4$ & $e_1 e_3=\alpha e_5$ & $e_1 e_4=\beta e_5$ & \\
      & $e_2 e_1=e_3$ & $e_2 e_2= (\alpha+\gamma) e_5$ & $e_3 e_1=(3\gamma+1) e_5$ & $e_4 e_1= e_5$ &\\ \hline
    
     $\T{5}{20}(\alpha, \beta)_{(\alpha; \beta)\neq(0 ; 0 )} $ &  $e_1 e_1 = e_2$ & $e_1 e_2=e_4$ & $e_1 e_3=\alpha e_5$ & $e_1 e_4= e_5$ &  \\
     &  $e_2 e_1=e_3$ & $e_2 e_2= (\alpha+\beta) e_5$ & $e_3 e_1=3\beta e_5$ &&\\ \hline
  \end{longtable}

%\[\T 5{17}(\alpha, \beta), \quad  \T 5{18}(\alpha), \quad  \T 5{19}(\alpha, \beta, \gamma), \quad \T 5{20}(\alpha, \beta).\]

\subsection{$1$-dimensional  central extensions of $\T {4}{02}(\alpha)$}\label{ext-T_02^4}
	
\begin{enumerate}

\item $\alpha\neq 6.$ Let us use the following notations:
	\begin{align*}
	\nb 1 = \Dl 13, \quad \nb 2 = \Dl 14 + 2\Dl 23-3\Dl 32, \quad \nb 3 = \Dl 21 .
	\end{align*}
	
Take $\0=\sum\limits_{i=1}^3\af_i\nb i\in {\rm H_T^2}(\T {4}{02}(\alpha)).$
	For the element
	$$
	\phi=
	\begin{pmatrix}
x &    0  &  0 &  0\\
0 &  x^2  &  0&  0\\
y &   0  &  x^3&  0\\
z &  xy(\alpha-3)  &   0 & x^4
\end{pmatrix}\in\aut{\T {4}{02}(\alpha)},
	$$
we have
	$$
	\phi^T\begin{pmatrix}
	0      &  0     & \af_1  & \af_2\\
	\af_3      & 0  & 2\af_2  & 0    \\
    0 & -3\af_2 & 0  & 0    \\
    0      & 0  & 0  & 0
	\end{pmatrix} \phi=
	\begin{pmatrix}
    \af^{**}      &  \af^{*}     & \af^*_1  & \af^*_2\\
	\af^{*}_3      & 0  & 2\af^*_2  & 0    \\
    0  & -3\af^*_2  & 0  & 0    \\
    0     & 0  & 0  & 0
	\end{pmatrix},
	$$
	where
\[
\begin{array}{rclcrclcrcl}
	\af^*_1 &=& x^4\af_1,&&
	\af^*_2 &=& x^5\af_2,&&
	\af^*_3 &=& x^3\af_3+2x^2y\af_2.
	\end{array} \]

Since $\af_2\ne 0,$ then choosing $y=-\frac{x\af_3}{2\af_2},$ we have the representatives $\la \nb 2\ra_{\alpha\neq 6}$ and
$\la \nb 1 + \nb 2\ra_{\alpha\neq 6}$ depending on whether $\alpha_1=0$ or not.

\item $\alpha= 6.$ Let us use the following notations:
	\begin{align*}
	\nb 1 = \Dl 13, \quad \nb 2 = \Dl 14 + 2\Dl 23-3\Dl 32, \quad \nb 3 = \Dl 21, \quad \nb 4 =\Dl 41 + \Dl 23-3\Dl 32 .
	\end{align*}
Take $\0=\sum\limits_{i=1}^4\af_i\nb i\in {\rm H_T^2}(\T {4}{02}(6)).$
	For the element
	$$
	\phi=
	\begin{pmatrix}
x &    0  &  0 &  0\\
0 &  x^2  &  0&  0\\
y &   0  &  x^3&  0\\
z &  3xy  &   0 & x^4
\end{pmatrix}\in\aut{\T {4}{02}(6)},
	$$
we have
	$$
	\phi^T\begin{pmatrix}
	0      &  0     & \af_1  & \af_2\\
	\af_3      & 0  & 2\af_2+\af_4  & 0    \\
    0 & -3(\af_2+\af_4) & 0  & 0    \\
    \af_4      & 0  & 0  & 0
	\end{pmatrix} \phi=
	\begin{pmatrix}
    \af^{**}      &  \af^{*}     & \af^*_1  & \af^*_2\\
	\af^{*}_3      & 0  & 2\af^*_2+\af^{*}_4  & 0    \\
    0  & -3(\af^*_2+\af^{*}_4)  & 0  & 0    \\
    \af^{*}_4     & 0  & 0  & 0
	\end{pmatrix},
	$$
where
\[
\begin{array}{rclcrclcrcl}
	\af^*_1 &=& x^4\af_1, &&
	\af^*_2 &=& x^5\af_2,\\
	\af^*_3 &=& x^3\af_3+2x^2y(\af_2+2\af_4), &&
    \af^*_4 &=& x^5\af_4.
	\end{array} \]

Since we are interested only in new algebras, we have $(\af_2, \af_4) \neq (0,0)$ and the following cases:

\begin{enumerate}
\item if $\alpha_4 =0,$ then $\af_2\ne 0,$ and choosing $y=-\frac{x\af_3}{2\af_2},$ we have the representatives $\la \nb 2\ra$ and
$\la \nb 1 + \nb 2\ra$ depending on whether $\alpha_1=0$ or not.

\item if $\alpha_4 \neq 0,$ then:
\begin{enumerate}
\item if $\alpha_2 \neq -2\alpha_4,$ then choosing $y=-\frac{x\af_3}{2(\af_2+2\af_4)},$ we have the representatives 
\[\la \beta\nb 2+ \nb 4\ra_{\beta\neq -2}\mbox{ and }\la \nb 1 + \beta \nb 2 + \nb 4\ra_{\beta\neq -2} \] depending on whether $\alpha_1=0$ or not.

\item if $\alpha_2 = -2\alpha_4,$ then we have the family of representatives \[\la \alpha_1\nb 1 - 2\alpha_4\nb 2 +\alpha_3\nb 3+\alpha_4\nb 4 \ra,\] which gives three distinct orbits with representatives
    \[\la - 2\nb 2 +\nb 4 \ra, \ \la \nb 1 - 2\nb 2 +\nb 4 \ra \mbox{ and }\la \beta\nb 1 - 2\nb 2 +\nb 3+\nb 4 \ra.\]
    
\end{enumerate}

\end{enumerate}

\end{enumerate}
		
Summarizing, all cases, we have the following representatives of distinct orbits:
 \[\la \nb 2\ra, \quad \la \nb 1 + \nb 2\ra, \quad \la \beta\nb 2+ \nb 4\ra_{\alpha=6},
 \quad \la \nb 1 + \beta \nb 2 + \nb 4\ra_{\alpha=6}, \quad \la \beta\nb 1 - 2\nb 2 +\nb 3+\nb 4 \ra_{\alpha=6}.\]

Hence, we have following $5$-dimensional algebras:

\begin{longtable}{llllllllllllll}
  \hline $\T{5}{21}(\alpha) $ &  $e_1 e_1 = e_2$ & $e_1 e_2=e_3$ & $e_1 e_3=\alpha e_4$ & $e_1 e_4=e_5$  &\\
        & $e_2 e_2= e_4$  & $e_2 e_3 = 2 e_5$ & $e_3 e_1 = -3 e_4$ & $e_3 e_2=-3e_5$ &\\
    \hline $\T{5}{22}(\alpha) $ &  $e_1 e_1 = e_2$ & $e_1 e_2=e_3$ & $e_1 e_3=\alpha e_4 +e_5$ & $e_1 e_4=e_5$ &\\
        & $e_2 e_2= e_4$  & $e_2 e_3 = 2 e_5$ & $e_3 e_1 = -3 e_4$ & $e_3 e_2=-3e_5$ &\\
    \hline $\T{5}{23}(\beta) $ &  $e_1 e_1 = e_2$ & $e_1 e_2=e_3$ & $e_1 e_3=6 e_4$ & $e_1 e_4=\beta e_5$ & $e_2 e_2= e_4$\\
        &  $e_2 e_3 = (2\beta+1) e_5$ & $e_3 e_1 = -3 e_4$  & $e_3 e_2 = -3(\beta+1) e_5$  & $e_4 e_1=e_5$&\\
    \hline $\T{5}{24}(\beta) $ &  $e_1 e_1 = e_2$ & $e_1 e_2=e_3$ & $e_1 e_3=6 e_4 + e_5$ & $e_1 e_4=\beta e_5$ & $e_2 e_2= e_4$\\
        &  $e_2 e_3 = (2\beta+1) e_5$ & $e_3 e_1 = -3 e_4$  & $e_3 e_2 = -3(\beta+1) e_5$  & $e_4 e_1=e_5$&\\
    \hline $\T{5}{25}(\beta) $ &  $e_1 e_1 = e_2$ & $e_1 e_2=e_3$ & $e_1 e_3=6 e_4 + \beta e_5$ & $e_1 e_4=-2 e_5$ & $e_2 e_1= e_5$\\
        & $e_2 e_2= e_4$ &  $e_2 e_3 = -3 e_5$ & $e_3 e_1 = -3 e_4$  & $e_3 e_2 = 3 e_5$  & $e_4 e_1=e_5$\\
 \hline
\end{longtable}

%\[\T 5{21}(\alpha), \quad \T 5{22}(\alpha), \quad \T 5{23}(\beta),\quad \T 5{24}(\beta),\quad \T 5{25}(\beta)\]

\subsection{$1$-dimensional  central extensions of $\T {4}{03}$}\label{ext-T_03^4}
	Let us use the following notations:
	\begin{align*}
	\nb 1 = \Dl 14, \quad \nb 2 = \Dl 21, \quad  \nb 3 = \Dl 22 - 3\Dl 31 .
	\end{align*}
Take $\0=\sum\limits_{i=1}^3\af_i\nb i\in {\rm H_T^2}(\T {4}{03}).$
	For the element
	$$
	\phi=
	\begin{pmatrix}
x &    0  &  0 &  0\\
y &  x^2  &  0&  0\\
z &  xy  &  x^3&  0\\
t &  xz &  x^2y & x^4
\end{pmatrix}\in\aut{\T {4}{03}},
	$$
we have 
	$$
	\phi^T\begin{pmatrix}
	0      & 0     & 0  & \af_1\\
	\af_2  & \af_3 & 0  & 0   \\
    -3\af_3 & 0 & 0  & 0    \\
    0      & 0  & 0  & 0
	\end{pmatrix} \phi=
	\begin{pmatrix}
    \af^{***}      &  \af^{**}     & \af^*  & \af^*_1\\
	\af^{*}_2      & \af^*_3  & 0  & 0    \\
    -3\af^*_3  & 0  & 0  & 0    \\
    0    & 0  & 0  & 0
	\end{pmatrix},
	$$
where
\[
\begin{array}{rclcrclcrcl}
	\af^*_1 &=& x^5\af_1, &&
	\af^*_2 &=& x^3\af_2-2x^2y\af_3, &&
	\af^*_3 &=& x^4\af_3.
	\end{array} \]

Since $\af_1\ne 0,$ we have the following cases:
\begin{enumerate}
 \item if $\af_3 \neq 0,$ then choosing $x = \frac {\alpha_3} {\alpha_1},$ $y = \frac{x\af_2}{2\af_3},$ we have the representatives $\la \nb 1 + \nb 3 \ra.$
  \item if  $\af_3 = 0,$ then we have the representatives $\la \nb 1 \ra$ and $\la \nb 1 + \nb 2 \ra$
  depending on whether $\alpha_2=0$ or not.
\end{enumerate}

Hence we have the following $5$-dimensional algebras:
\begin{longtable}{llllllll}
  \hline $\T{5}{26} $ &  $e_1 e_1 = e_2$ & $e_1 e_2=e_3$ & $e_1 e_3= e_4$&  $e_1 e_4=e_5$ & $e_2 e_2= e_5$ & $e_3 e_1 = -3e_5$ \\
    \hline $\T{5}{27} $ &  $e_1 e_1 = e_2$ & $e_1 e_2=e_3$ & $e_1 e_3= e_4$ & $e_1 e_4=e_5$ & \\
    \hline $\T{5}{28} $ &  $e_1 e_1 = e_2$ & $e_1 e_2=e_3$ & $e_1 e_3= e_4$ & $e_1 e_4=e_5$ & $e_2 e_1= e_5$\\
\hline  
\end{longtable}%\[\T 5{26}, \quad \T 5{27}, \quad \T 5{28}.\]

\subsection{$1$-dimensional  central extensions of $\T {4}{04}$}\label{ext-T_04^4}
	Let us use the following notations:
	\begin{align*}
	\nb 1 = \Dl 13, \quad \nb 2 = \Dl 14+ 3\Dl 31, \quad \nb 3 = \Dl 22 - 3\Dl 31.
	\end{align*}
	Take $\0=\sum\limits_{i=1}^3\af_i\nb i\in {\rm H_T^2}(\T {4}{04}).$
	For the element
	$$
	\phi=
	\begin{pmatrix}
	1 &    0  &  0&  0\\
	x &  1  &  0&  0\\
	y &   x  &  1 &  0\\
    z &   y+z  &  x &  1
	\end{pmatrix}\in\aut{\T {4}{04}},
	$$
	we have
	$$
	\phi^T\begin{pmatrix}
	0      &  0     & \af_1  & \af_2\\
	0      & \af_3  & 0  & 0    \\
    3(\af_2-\af_3)  & 0  & 0  & 0    \\
    0      & 0  & 0  & 0
	\end{pmatrix} \phi=
	\begin{pmatrix}
    \af^{***}      &  \af^{**}     & \af^{*} + \af^*_1  & \af^*_2\\
	\af^{*}      & \af^*_3  & 0  & 0    \\
    3(\af^*_2-\af^*_3)  & 0  & 0  & 0    \\
    0     & 0  & 0  & 0
	\end{pmatrix},
	$$
	where
\[
\begin{array}{rclcrclcrcl}
	\af^*_1 &=& \af_1+2x(\af_3-\af_2),&&
	\af^*_2 &=& \af_2,&&
	\af^*_3 &=& \af_3.
	\end{array} \]

If $\af_2\ne \af_3,$ then choosing $x=\frac{\af_1}{2(\af_2-\af_3)},$ we have the representatives $\la \nb 2+\af\nb 3\ra_{\alpha\neq1}.$
If  $\af_2= \af_3,$ then we have the  representatives $\la \af\nb 1 + \nb 2+\nb 3\ra.$

Hence, we have following $5$-dimensional algebras:
\begin{longtable}{llllllllll}
  \hline $\T{5}{29}(\alpha)_{\alpha\neq 1} $ &  $e_1 e_1 = e_2$ & $e_1 e_2=e_3$ & $e_1 e_3=e_4$ & $e_1 e_4=e_5$ & \\ 
        & $e_2 e_1= e_4$& $e_2 e_2 = \alpha e_5$ & $e_3 e_1=3(1 - \alpha)e_5$ && \\
    
    \hline $\T{5}{30}(\alpha) $ &  $e_1 e_1 = e_2$ & $e_1 e_2=e_3$ & $e_1 e_3=e_4+\alpha e_5$  &&\\
        & $e_1 e_4=e_5$ & $e_2 e_1= e_4$ & $e_2 e_2 =  e_5$ && \\
        \hline
 \end{longtable}

%\[\T 5{29}(\alpha), \quad \T 5{30}(\alpha)\]

\subsection{$1$-dimensional  central extensions of $\T {4}{05}(\lambda, \alpha)$}\label{ext-T_05^4}

\begin{enumerate}
    \item $\lambda =1, \alpha=0.$ Let us use the following notations:
	\begin{align*}
	\nb 1 = \Dl 12, \quad \nb 2 = \Dl 13+ \Dl 22, \quad \nb 3 = \Dl 14+ \Dl 32, \quad \nb 4 = \Dl 34.
	\end{align*}
	Take $\0=\sum\limits_{i=1}^4\af_i\nb i\in {\rm H_T^2}(\T {4}{05}(1, 0)).$
	For the element
	$$
	\phi=
	\begin{pmatrix}
x &    0  &  0 &  0\\
0 &  x^2  &  0&  0\\
y &   0  &  x^3&  0\\
z &  xy  &   0 & x^4
\end{pmatrix}\in\aut{\T {4}{05}(1, 0)},
	$$ we have
\[\phi^T\left(
             \begin{array}{cccc}
               0 & \alpha_1 & \alpha_2 & \alpha_3 \\
               0 & \alpha_2 & 0 & 0 \\
               0 &  \alpha_3 & 0 & \alpha_4 \\
               0 & 0 & 0 & 0 \\
             \end{array}
           \right)\phi=\left(
             \begin{array}{cccc}
               \alpha^{*} & \alpha_1^* & \alpha_2^* & \alpha_3^* \\
               0 & \alpha_2^* & 0 & 0 \\
               \alpha^{**} &  \alpha_3^* & 0 & \alpha_4^* \\
               0 & 0 & 0 & 0 \\
             \end{array}
           \right),\]
where
\[
\begin{array}{rclcrclcrcl}
	\alpha^*_1 &=& x^3\alpha_1+2x^2y\alpha_3+xy^2\alpha_4,&&
	\alpha^*_2 &=& x^4 \alpha _2,\\
	\alpha^*_3 &=& x^5\alpha_3+x^4y\alpha_4, &&
    \alpha^*_4 &=& x^7 \alpha _4.
	\end{array} \]

Since we are interested only in new algebras, we have $(\af_3, \af_4) \neq (0,0)$ and the following cases:
\begin{enumerate}
\item if $\alpha_4 =0,$ then:
\begin{enumerate}
\item if $\alpha_2 =0,$ then choosing $y=-\frac{x\af_1}{2\af_3},$ we have the representatives $\la \nb 3\ra$.
\item if $\alpha_2 \neq 0,$ then choosing $x=\frac{\af_2}{\af_3},$ $y=-\frac{x\af_1}{2\af_3},$ we have the representatives $\la \nb 2 + \nb 3\ra$.

\end{enumerate}

\item if $\alpha_4 \neq 0,$ then:

\begin{enumerate}
\item if $\alpha_2 =0,$  $\alpha_1\alpha_4 -\alpha_3^2 =0,$ then choosing $y=-\frac{x\af_3}{\af_4},$ we have the representatives $\la \nb 4\ra$.
\item if $\alpha_2 =0,$  $\alpha_1\alpha_4 -\alpha_3^2 \neq0,$ then choosing $x=\sqrt[4]{\frac{\alpha_1\alpha_4 -\alpha_3^2}{\af_4^2}},$ $y=-\frac{x\af_3}{\af_4},$ we have the representatives $\la \nb 1 + \nb 4\ra$.
\item if $\alpha_2 \neq 0,$  then choosing $x=\sqrt[3]{\frac{\alpha_2}{\af_4}},$ $y=-\frac{x\af_3}{\af_4},$ we have the representatives \[\la \beta\nb 1 + \nb 2 +\nb 4\ra. \]

\end{enumerate}

\end{enumerate}

\item $\lambda =0, \alpha=\frac 1 3.$ Let us use the following notations:
\begin{align*}
	\nb 1 = \Dl 12, \quad \nb 2 = \Dl 13+ \Dl 22, \quad \nb 3 = \Dl 23, \quad \nb 4 = \Dl 24 + \Dl 33.
	\end{align*}

Take $\0=\sum\limits_{i=1}^4\af_i\nb i\in {\rm H_T^2}(\T{4}{05}(0, \frac 1 3)).$
	For the element
	$$
	\phi=
	\begin{pmatrix}
x &    0  &  0 &  0\\
0 &  x^2  &  0&  0\\
y &   0  &  x^3&  0\\
z &  \frac 4 3 xy  &   0 & x^4
\end{pmatrix}\in\aut{\T{4}{05}(0, \frac 1 3)},
	$$
we have 
 \[\phi^T\left(
             \begin{array}{cccc}
               0 & \alpha_1 & \alpha_2 & 0 \\
               0 & \alpha_2 & \alpha_3 & \alpha_4 \\
               0 & 0 & \alpha_4 & 0 \\
               0 & 0 & 0 & 0
             \end{array}
           \right)\phi=\left(
             \begin{array}{cccc}
               \alpha^{***} & \alpha_1^* & \alpha^{*} + \alpha_2^* & 0 \\
               \alpha^{**} & 2 \alpha^{*} + \alpha_2^* & \alpha_3^* & \alpha_4^* \\
               3 \alpha^{*} & 0 & \alpha_4^* & 0 \\
               0 & 0 & 0 & 0
             \end{array}
           \right)\]
where
\[
\begin{array}{rclcrclcrcl}
	\alpha^*_1 &=& x^3 \alpha _1, &&
	\alpha^*_2 &=& x^4\alpha _2 +\frac{2}{3} x^3y \alpha _4,\\
	\alpha^*_3 &=& x^5 \alpha _3, &&
    \alpha^*_4 &=& x^6 \alpha _4.
	\end{array} \]
%\begin{enumerate}
%\item $\alpha_4 =0,$ then $\af_3\ne 0$ and we have the representatives $\la \nb 3\ra,$ $\la \nb 1 +\nb 3\ra$ and $\la \beta \nb 1 +\nb 2 +\nb 3\ra$ depending on whether $\alpha_1=0,$ $\alpha_2=0$ or not.
%\item
Since $\alpha_4 \neq 0,$ choosing $y=-\frac{3x\alpha_2}{2\alpha_4},$ we have the representatives 
\[\la \nb 4\ra, \ \la \nb 1 +\nb 4\ra \mbox{ and }\la \beta \nb 1 +\nb 3 +\nb 4\ra\]
depending on whether $\alpha_1=0,$ $\alpha_3=0$ or not.

%\end{enumerate}

\item $\lambda =0, \alpha=-\frac 2 3.$ Let us use the following notations:
\begin{align*}
	\nb 1 = \Dl 12, \quad \nb 2 = \Dl 13+ \Dl 22, \quad \nb 3 = \Dl 23, \quad \nb 4 = \Dl 14+ \Dl 32 .
	\end{align*}
Take $\0=\sum\limits_{i=1}^4\af_i\nb i\in {\rm H_T^2}(\T{4}{05}(0, -\frac 2 3)).$
	For the element
	$$
	\phi=
	\begin{pmatrix}
x &    0  &  0 &  0\\
0 &  x^2  &  0&  0\\
y &   0  &  x^3&  0\\
z &  \frac 1 3 xy  &   0 & x^4
\end{pmatrix}\in\aut{\T{4}{05}(0, -\frac 2 3)},
$$
we have  \[\phi^T\left(
             \begin{array}{cccc}
               0 & \alpha_1 & \alpha_2 & \alpha_4 \\
               0 & \alpha_2 & \alpha_3 & 0 \\
               0 & \alpha_4 & 0 & 0 \\
               0 & 0 & 0 & 0
             \end{array}
           \right)\phi=\left(
             \begin{array}{cccc}
               \alpha^* & \alpha_1^* & \alpha_2^* & \alpha_4^* \\
               0 & \alpha_2^* & \alpha_3^* & 0 \\
               0 & \alpha_4^* & 0 & 0 \\
               0 & 0 & 0 & 0
             \end{array}
           \right),\]
where
\[
\begin{array}{rclcrclcrcl}
	\alpha^*_1 &=&x^3 \alpha _1 +\frac{4}{3}x^2y \alpha _4,&&
	\alpha^*_2 &=& x^4 \alpha_2 ,\\
	\alpha^*_3 &=& x^5 \alpha_3,&&
	\alpha^*_4 &=& x^5 \alpha_4.
	\end{array} \]

Since $\alpha_4\neq 0,$ choosing $y=-\frac{3x\alpha_1}{4\alpha_4},$ we have the representatives  $\la \beta \nb 3 +  \nb 4 \ra$ and $\la \nb 2 + \beta\nb 3 +  \nb 4\ra,$ depending on whether $\alpha_2=0$ or not.

\item $\lambda =0, \alpha=-\frac 1 5.$
Let us use the following notations:
\begin{align*}
	\nb 1 = \Dl 12, \quad \nb 2 = \Dl 13+ \Dl 22, \quad \nb 3 = \Dl 23, \quad \nb 4 = \Dl 14 - \frac  {2}{5}\Dl 32-\frac  {7}{2}\Dl 41.
	\end{align*}

Take $\0=\sum\limits_{i=1}^4\af_i\nb i\in {\rm H_T^2}(\T{4}{05}(0, -\frac{1} 5)).$
	For the element
	$$
	\phi=
	\begin{pmatrix}
x &    0  &  0 &  0\\
0 &  x^2  &  0&  0\\
y &   0  &  x^3&  0\\
z &  \frac 4 5 xy&0 & x^4
\end{pmatrix}\in\aut{\T{4}{05}(0, -\frac{1} 5)},$$
we have
\[\phi^T\left(
             \begin{array}{cccc}
               0 & \alpha_1 & \alpha_2 & \alpha_4 \\
               0 & \alpha_2 & \alpha_3 & 0 \\
               0 & -\frac  {2}{5}\alpha_4 & 0 & 0 \\
               -\frac  {7}{2}\alpha_4 & 0 & 0 & 0
             \end{array}
           \right)\phi=\left(
             \begin{array}{cccc}
               \alpha^{**} & \alpha_1^* & \alpha_2^* & \alpha_4^* \\
               \alpha^{*} & \alpha_2^* & \alpha_3^* & 0 \\
               0 & -\frac  {2}{5}\alpha_4^* & 0 & 0 \\
               -\frac  {7}{2}\alpha_4^* & 0 & 0 & 0
             \end{array}
           \right),\]
where
\[
\begin{array}{rclcrclcrcl}
	\alpha^*_1 &=& x^3\alpha _1 +\frac{2}{5} x^2y \alpha _4,&&
	\alpha^*_2 &=& x^4 \alpha _2,\\
	\alpha^*_3 &=& x^5 \alpha _3, &&
    \alpha^*_4 &=& x^5 \alpha _4.
	\end{array} \]

Since $\alpha_4 \neq 0,$ choosing $y=-\frac{5x\af_1}{2\af_4},$ we have the representatives $\la \beta \nb 3 + \nb 4 \ra$ or $\la \nb 2 + \beta \nb 3 + \nb 4 \ra$ depending on whether $\alpha_2=0$ or not.

\item $\lambda =-1.$
Let us use the following notations:
\begin{longtable}{ll}
	$\nb 1 = \Dl 12,$ & $\nb 3 = \Dl 14-\frac  {3\alpha+4}{3}\Dl 23+\Dl 32,$\\ 
	$\nb 2 = \Dl 13+ \Dl 22,$& $ \nb 4 = \frac  {3\alpha+1}{3}\Dl 23-\frac  {3\alpha+1}{2}\Dl 32+\Dl 41.$
	\end{longtable}

Take $\0=\sum\limits_{i=1}^4\af_i\nb i\in {\rm H_T^2}(\T{4}{05}(-1, \alpha)).$
	For the element
	$$ \phi=
	\begin{pmatrix}
x &    0  &  0 &  0\\
0 &  x^2  &  0&  0\\
z &   0  &  x^3&  0\\
t &  xz(\alpha+1)&0 & x^4
\end{pmatrix}\in\aut{\T{4}{05}(-1, \alpha)_{\alpha\neq -\frac 2 3}}, \quad
\phi=
	\begin{pmatrix}
x &    0  &  0 &  0\\
y &  x^2  &  0&  0\\
z &   0  &  x^3&  0\\
t &  \frac 1 3 xz&0 & x^4
\end{pmatrix}\in\aut{\T{4}{05}(-1, -\frac 2 3)},$$
we have
\[\phi^T\left(
             \begin{array}{cccc}
               0 & \alpha_1 & \alpha_2 & \alpha_3 \\
               0 & \alpha_2 & -\frac  {3\alpha+4}{3}\alpha_3+\frac  {3\alpha+1}{3}\alpha_4& 0 \\
               0 & \alpha_3-\frac  {3\alpha+1}{2}\alpha_4 & 0 & 0 \\
               \alpha_4 & 0 & 0 & 0
             \end{array}
           \right)\phi=\left(
             \begin{array}{cccc}
               \alpha^{**} & -\alpha^*+\alpha_1^* & \alpha_2^* & \alpha_3^* \\
               \alpha^* & \alpha_2^* & -\frac  {3\alpha+4}{3}\alpha_3^*+\frac  {3\alpha+1}{3}\alpha_4^*& 0 \\
               0 & \alpha_3^*-\frac  {3\alpha+1}{2}\alpha_4^* & 0 & 0 \\
               \alpha_4^* & 0 & 0 & 0
             \end{array}
           \right)\]
where
\[\alpha \neq -\frac 2 3:  \left\{\begin{array}{lllll}
	\alpha^*_1 =x^3\alpha _1+\frac{4 \alpha _3+(3 \alpha +5) \alpha _4}{6}x^2z ,\\
	\alpha^*_2 = x^4 \alpha _2,\\
	\alpha^*_3 = x^5 \alpha _3,\\
    \alpha^*_4 = x^5 \alpha _4.
	\end{array}\right., \quad \alpha = -\frac 2 3: \left\{ \begin{array}{lllll}
	\alpha^*_1 =x^3\alpha _1+2x^2y\alpha _2+\frac{4 \alpha _3+3\alpha _4}{6} x^2z,\\
	\alpha^*_2 = x^4 \alpha _2,\\
	\alpha^*_3 = x^5 \alpha _3,\\
    \alpha^*_4 = x^5 \alpha _4.
	\end{array}\right.\]

Since $(\af_3, \af_4) \neq (0,0),$  we have the following cases:

\begin{enumerate}
\item if $\alpha_4 =0,$ then by choosing  $z=-\frac{3x\af_1}{2\af_3}$ (or $z=-\frac{3x\af_1+6y\alpha_2}{2\af_3}$), we have the representatives $\la \nb 3\ra$ and $\la \nb 2 + \nb 3\ra$ depending on whether $\alpha_2=0$ or not.

\item if $\alpha_4 \neq 0,$ then:

\begin{enumerate}
\item if $4 \alpha _3+(3 \alpha +5) \alpha _4 \neq 0,$ then by choosing  $z=-\frac{6x\af_1}{4 \alpha _3+(3 \alpha +5) \alpha _4}$ (or $z=-\frac{6x\af_1+12y\alpha_2}{4 \alpha _3+3 \alpha _4}$), we have the representatives $\la \beta \nb 3 + \nb 4\ra$ and $\la \nb 2 + \beta \nb 3 + \nb 4\ra$ depending on whether $\alpha_2=0$ or not.

\item if $4 \alpha _3+(3 \alpha +5) \alpha _4 =0,$  then in case of $\alpha \neq -\frac 2 3,$ we have the representatives
\[\la -\frac{3\alpha+5} 4\nb 3 + \nb 4\ra, \ \la \nb 2 -\frac{3\alpha+5} 4\nb 3 + \nb 4\ra \mbox{ and }\la  \nb 1 + \beta \nb 2-\frac{3\alpha+5} 4\nb 3 + \nb 4\ra\]
depending on whether $\alpha_1=0,$ $\alpha_2=0$ or not.

If $\alpha = -\frac 2 3,$ we have the representatives
\[\la -\frac{3} 4\nb 3 + \nb 4\ra, \ \la \nb 1 -\frac{3} 4\nb 3 + \nb 4\ra\mbox{ and }\la \nb 2-\frac{3} 4\nb 3 + \nb 4\ra\] 
depending on whether $\alpha_1=0,$ $\alpha_2=0$ or not.

\end{enumerate}
\end{enumerate}

\item  $\alpha=\frac{4\lambda-1}{5-2\lambda},$ $\lambda \neq 0; -1; \frac 5 2.$
Let us use the following notations:
\begin{longtable}{ll}
	$\nb 1 = \Dl 12,$&$ \nb 3 = \Dl 14+\frac  {4\lambda^2-2\lambda+7}{3\lambda(5-2\lambda)}\Dl 23+\Dl 32,$\\
 $\nb 2 = \Dl 13+ \Dl 22, $ & $\nb 4 = \frac  {2(\lambda-1)^2}{3\lambda(5-2\lambda)}\Dl 23-\frac  {2(\lambda-1)}{5-2\lambda}\Dl 32+\Dl 41.$
	\end{longtable}
	
Take $\0=\sum\limits_{i=1}^4\af_i\nb i\in {\rm H_T^2}(\T{4}{05}(\lambda, \frac{4\lambda-1}{5-2\lambda})).$
	For the element
	$$
	\phi=
	\begin{pmatrix}
x &    0  &  0 &  0\\
y &  x^2  &  0&  0\\
z &  2xy  &  x^3&  0\\
t &  2xz+y^2  &  3x^2y & x^4
\end{pmatrix}\in\aut{\T {4}{05}(1, 1)},\quad
	\phi=
	\begin{pmatrix}
x &    0  &  0 &  0\\
y &  x^2  &  0&  0\\
z &   \frac{1}{2}xy  &  x^3&  0\\
t &  \frac 1 2 xz &\frac{1}{2}x^2y & x^4
\end{pmatrix}\in\aut{\T{4}{05}(-\frac 1 2, -\frac 1 2)},$$
$$
	\phi=
	\begin{pmatrix}
x &    0  &  0 &  0\\
0 &  x^2  &  0&  0\\
z &   0  &  x^3&  0\\
t &  \frac{2(\lambda+2)}{5-2\lambda}xz&0 & x^4
\end{pmatrix}\in\aut{\T{4}{05}(\lambda, \frac{4\lambda-1}{5-2\lambda})},\quad \lambda \neq 1; -\frac 1 2,$$
we have
\[\phi^T\left(
             \begin{array}{cccc}
               0 & \alpha_1 & \alpha_2 & \alpha_3 \\
               0 & \alpha_2 & \frac  {4\lambda^2-2\lambda+7}{3\lambda(5-2\lambda)}\alpha_3+\frac  {2(\lambda-1)^2}{3\lambda(5-2\lambda)}\alpha_4& 0 \\
               0 & \alpha_3-\frac  {2(\lambda-1)}{(5-2\lambda)}\alpha_4 & 0 & 0 \\
               \alpha_4 & 0 & 0 & 0
             \end{array}
           \right)\phi=\]
\[\left(
             \begin{array}{cccc}
               \alpha^{**} & \lambda\alpha^{*}+\alpha_1^* & \alpha_2^* & \alpha_3^* \\
               \alpha^{*} & \alpha_2^* & \frac  {4\lambda^2-2\lambda+7}{3\lambda(5-2\lambda)}\alpha_3^*+\frac  {2(\lambda-1)^2}{3\lambda(5-2\lambda)}\alpha_4^*& 0 \\
               0 & \alpha_3^*-\frac  {2(\lambda-1)}{(5-2\lambda)}\alpha_4^* & 0 & 0 \\
               \alpha_4^* & 0 & 0 & 0
             \end{array}
           \right),\]
where
\[\begin{array}{rcl}
\lambda = 1&:& \left\{\begin{array}{lll}
	\alpha^*_1 =x^3 \alpha _1 +2 x^2 y \alpha _2 +(2 x^2 z+ x y^2) (\alpha _3-\alpha _4),\\
	\alpha^*_2 = x^4 \alpha _2+3x^3 y(\alpha _3-\alpha _4),\\
	\alpha^*_3 = x^5 \alpha _3,\\
    \alpha^*_4 = x^5 \alpha _4.
	\end{array}\right.\\
\lambda = -\frac 1 2&:& \left\{\begin{array}{lll}
	\alpha^*_1 =x^3\alpha_1+2x^2y\alpha_2 + (x^2z-\frac 1 4 xy^2)(\alpha _3 + \frac 1 2\alpha_4),\\
	\alpha^*_2 = x^4 \alpha _2,\\
	\alpha^*_3 = x^5 \alpha _3,\\
    \alpha^*_4 = x^5 \alpha _4.
	\end{array}\right.\\
\lambda \neq 1; -\frac 1 2&:& \left\{\begin{array}{lll}
	\alpha^*_1 =x^3\alpha _1+x^2z\frac{2 (\lambda +2)}{3}\big(\alpha_3-\frac{ 4\lambda-1}{5-2\lambda}\alpha_4\big),\\
	\alpha^*_2 = x^4 \alpha _2,\\
	\alpha^*_3 = x^5 \alpha _3,\\
    \alpha^*_4 = x^5 \alpha _4.
	\end{array}\right.
	\end{array}\]

Since $(\af_3, \af_4) \neq (0,0),$  we have the following cases:
\begin{enumerate} \item $\lambda =1$:

\begin{enumerate}
\item if $\alpha_4 =0,$ then  choosing $y=-\frac{x\af_2}{3\af_3},$ $z=-\frac{x^2\af_1+2xy\alpha_2+y^2\af_3}{2x\af_3},$ we have the representatives $\la \nb 3\ra$.
\item if $\alpha_4 \neq 0,$ then:
\begin{enumerate}
\item if $\alpha_3 = \alpha_4,$ $\alpha_2 =0,$ then we have the representatives
$\la  \nb 3 + \nb 4\ra$ and  $\la \nb 1 + \nb 3 + \nb 4\ra$ depending on whether $\alpha_1=0$ or not.
\item if $\alpha_3 = \alpha_4,$ $\alpha_2 \neq 0,$ then choosing $y = -\frac {x\af_1}{2\af_2},$ we have the representative
$\la \nb 2 + \nb 3 + \nb 4\ra$.

\item if $\alpha_3 \neq \alpha_4,$  then choosing $y=-\frac{x\af_2}{3(\af_3-\alpha_4)},$ $z=-\frac{x^2\af_1+2xy\alpha_2+y^2(\af_3-\af_4)}{2x(\af_3-\af_4)},$ we have the representatives $\la \beta \nb 3 + \nb 4\ra.$
\end{enumerate}

\end{enumerate}

\item $\lambda =-\frac 1 2$:

\begin{enumerate}
\item if $\alpha_4 =0,$ then choosing $y=0,$ $z=-\frac{x\af_1}{\af_3},$ we have the representatives $\la \nb 3\ra$ and $\la \nb 2 + \nb 3\ra$ depending on whether $\alpha_2=0$ or not.
\item if $\alpha_4 \neq 0,$ then:
\begin{enumerate}
\item if $\alpha_3 =- \frac 1 2 \alpha_4,$ $\alpha_2 =0,$ then we have the representatives
$\la - \frac 1 2 \nb 3 + \nb 4\ra$ and  $\la \nb 1 - \frac 1 2 \nb 3 + \nb 4\ra$ depending on whether $\alpha_1=0$ or not.
\item if $\alpha_3 = - \frac 1 2 \alpha_4,$ $\alpha_2 \neq 0,$ then choosing $y = -\frac {x\af_1}{2\af_2},$ we have the representative
$\la \nb 2 - \frac 1 2 \nb 3 + \nb 4\ra$.
\item if $\alpha_3 \neq - \frac 1 2 \alpha_4,$  then choosing $z=-\frac{8x^2\af_1+16xy\alpha_2- y^2(2\af_3+\af_4)}{4x(2\af_3+\af_4)},$ we have the representative
$\la \beta \nb 3 + \nb 4\ra$ and $\la \nb 2 + \beta \nb 3 + \nb 4\ra$ depending on whether $\alpha_2=0$ or not.

\end{enumerate}

\end{enumerate}

\item $\lambda =-2$:

\begin{enumerate}
\item if $\alpha_4 =0,$ then we have the representatives 
\[\la \nb 3\ra, \ \la \nb 1 + \nb 3\ra\mbox{ and }\la  \beta \nb 1 + \nb 2 + \nb 3\ra\] depending on whether $\alpha_1=0,$ $\alpha_2=0$ or not.
\item if $\alpha_4 \neq 0,$ then we have the representatives $\la \beta \nb 3 + \nb 4\ra,$  $\la \nb 2 + \beta \nb 3 + \nb 4\ra$ and $\la \nb 1 + \beta \nb 2 +  \gamma \nb 3 + \nb 4\ra$ depending on whether $\alpha_1=0,$ $\alpha_2=0$ or not.

\end{enumerate}

\item $\lambda \neq 1; -2; - \frac 1 2$:

\begin{enumerate}
\item if $\alpha_4 =0,$ then choosing $z=-\frac{3x\af_1}{2(\lambda+2)\af_3},$
we have the representatives $\la \nb 3\ra$ and $\la \nb 2 + \nb 3\ra$  depending on whether $\alpha_2=0$ or not.
\item if $\alpha_4 \neq 0,$ then:
\begin{enumerate}
\item if $\alpha_3 =\frac{ 4\lambda-1}{5-2\lambda}\alpha_4,$ then we have the representatives
$\la \frac{ 4\lambda-1}{5-2\lambda} \nb 3 + \nb 4\ra,$ $\la \nb 2 + \frac{ 4\lambda-1}{5-2\lambda} \nb 3 +  \nb 4\ra$ and
$\la \nb 1 + \beta \nb 2 + \frac{ 4\lambda-1}{5-2\lambda} \nb 3  + \nb 4\ra$ depending on whether $\alpha_1=0,$ $\alpha_2=0$ or not.
\item if $\alpha_3 \neq \frac{ 4\lambda-1}{5-2\lambda}\alpha_4 ,$  then choosing $z=-\frac{3x\af_1(5-2\lambda)}{2(\lambda+2)\big((5-2\lambda)\alpha_4 - (4\lambda-1)\af_3\big)},$ we have the representative
$\la \beta \nb 3 + \nb 4\ra$ and $\la \nb 2 + \beta \nb 3 + \nb 4\ra$ depending on whether $\alpha_2=0$ or not.

\end{enumerate}

\end{enumerate}
  \end{enumerate}

\

\item $\lambda \neq 0; -1,$ $\alpha\neq\frac{4\lambda-1}{5-2\lambda},$ $(\lambda; \alpha) \neq (1; 0).$
	Let us use the following notations:
	\begin{align*}
	\nb 1 = \Dl 12, \quad \nb 2 = \Dl 13+ \Dl 22, \quad \nb 3 = \Dl 14+\frac  {3\alpha-2(\lambda-1)}{3\lambda}\Dl 23+\Dl 32.
	\end{align*}

	Take $\0=\sum\limits_{i=1}^3\af_i\nb i\in {\rm H_T^2}(\T {4}{05}(\lambda, \alpha)).$
	For the element
$$
	\phi=
	\begin{pmatrix}
x &    0  &  0 &  0\\
0 &  x^2  &  0&  0\\
z &   0  &  x^3&  0\\
t &  (\alpha+1)xz&0 & x^4
\end{pmatrix}\in\aut{\T{4}{05}(\lambda, \alpha)},\quad \alpha\neq \frac{2\lambda^2+5\lambda-1}{6},$$
$$
	\phi=
	\begin{pmatrix}
x &    0  &  0 &  0\\
y &  x^2  &  0&  0\\
z &   (\lambda+1)xy  &  x^3&  0\\
t &  \frac  {2\lambda^2+5\lambda+5}{6}xz+\frac  {2\lambda^2+3\lambda+1}{6}y^2&\frac  {2\lambda^2+9\lambda+7}{6}x^2y & x^4
\end{pmatrix}\in\aut{\T{4}{04}\left(\lambda, \frac{2\lambda^2+5\lambda-1}{6}\right)},$$
we have
\[\phi^T\left(
             \begin{array}{cccc}
               0 & \alpha_1 & \alpha_2 & \alpha_3 \\
               0 & \alpha_2 & \frac  {3\alpha-2(\lambda-1)}{3\lambda}\alpha_3& 0 \\
               0 & \alpha_3 & 0 & 0 \\
               0 & 0 & 0 & 0
             \end{array}
           \right)\phi=\left(
             \begin{array}{cccc}
               \alpha^{**} & \lambda \alpha^{*}+\alpha_1^* & \alpha\alpha^*+\alpha_2^* & \alpha_3^* \\
               \alpha^{*} & \big(\alpha +\frac{1-\lambda} 3 \big)\alpha^{*} +\alpha_2^* & \frac  {3\alpha-2(\lambda-1)}{3\lambda}\alpha_3^*& 0 \\
               \alpha^{*} & \alpha_3^* & 0 & 0 \\
               0 & 0 & 0 & 0
             \end{array}
           \right),\]
where

\[\begin{array}{rcl}
\alpha\neq \frac{2\lambda^2+5\lambda-1}{6}&:& \left\{\begin{array}{llll}
	\alpha^*_1 =x^3\alpha_1+\frac{2}{3}x^2y(\lambda +2)\alpha_3,\\
	\alpha^*_2 = x^4 \alpha _2,\\
	\alpha^*_3 = x^5 \alpha _3.
	\end{array}\right.\\

\alpha= \frac{2\lambda^2+5\lambda-1}{6}&:& \left\{\begin{array}{llll}
	\alpha^*_1=x^3\alpha_1+2x^2y\alpha_2+\frac{(\lambda+1)(-2 \lambda ^2+2 \lambda +3)xy^2+4\lambda(\lambda +2)x^2z}{6 \lambda }\alpha _3,\\
	\alpha^*_2= x^4\alpha_2+\frac{(2\lambda+1)(\lambda+1)}{2\lambda }x^3y\alpha _3,\\
	\alpha^*_3= x^5 \alpha _3.
	\end{array}\right.
	\end{array}\]
	
Since $\af_3 \neq 0,$  we have the following cases:
\begin{enumerate}
\item $\lambda =-2,$ $\alpha =-\frac 1 2,$ then choosing $y=\frac {4x\alpha_2}{3\alpha_3}$
we have the representatives $\la \nb 3\ra$ and  $\la \nb 1 + \nb 3\ra$ depending on whether $3\alpha_1\alpha_3 +4 \alpha_2^2=0$  or not.
\item $\lambda =-2,$ $\alpha \neq -\frac 1 2,$ then 
we have the representatives 
\[\la \nb 3\ra, \ \la \nb 1 + \nb 3\ra\mbox{ and }\la\beta \nb 1 +  \nb 2 + \nb 3\ra\] 
depending on whether $\alpha_1=0,$ $\alpha_2=0$ or not.
\item $\lambda \neq -2,$ $\alpha = \frac{2\lambda^2+5\lambda-1}{6},$ then choosing $y=-\frac {2\lambda x\alpha_2} {(2\lambda+1)(\lambda+1)\alpha_3},$
    $z=-\frac{6 x^2\lambda \alpha_1+12xy\lambda \alpha_2-(\lambda+1)(2 \lambda ^2-2 \lambda -3)xy\alpha _3}{4\lambda\alpha _3},$
we have the representatives $\la \nb 3\ra.$ 
\item $\lambda \neq -2,$ $\alpha \neq \frac{2\lambda^2+5\lambda-1}{6},$ then choosing $y=-\frac {3x\alpha_1} {2(\lambda+2)\alpha_3},$
    we have the representatives $\la \nb 3\ra$  and $\la \nb 2 + \nb 3\ra$ depending on whether $\alpha_2=0$ or not.

\end{enumerate}

\end{enumerate}

Thus, we have the following distinct orbits:
\begin{longtable}{llllllll}
$\lambda= 1, \ \alpha=0$&$:$   & $\la \beta \nb 1+ \nb 2 + \nb 4\ra, $ & $\la \nb 1 + \nb 4\ra,$ & $\la \nb 4\ra,$  \\[2mm]
\hline$\lambda= 0, \ \alpha=\frac 1 3$&$:$    & $
\la \beta \nb 1 +\nb 3+ \nb 4\ra,$  & $\la \nb 1 + \nb 4\ra, $& $\la \nb 4 \ra,$ \\[2mm]
\hline$\lambda= 0, \ \alpha=-\frac 2 3 $&$:$ & 
$\la \nb 2 + \beta \nb 3 +  \nb 4\ra,$ & $\la \beta \nb 3 + \nb 4 \ra,$ \\[2mm]
\hline$\lambda= 0, \ \alpha=-\frac 1 5 $&$:$  & 
$\la \nb 2 + \beta \nb 3+\nb 4\ra,$ & $\la \beta\nb 3 + \nb 4\ra,$ \\[2mm]
\hline$\lambda= -1$&$:$  & 
$\la \nb 2 + \beta \nb 3 + \nb 4\ra,$ & $\la \beta \nb 3 + \nb 4\ra,$ \\[2mm]
\hline&& $\la \nb 1 + \beta\nb 2 - \frac {3\alpha+5} 4 \nb 3 + \nb 4\ra_{\alpha \neq -\frac 2 3},$ &
 $\la \nb 1 - \frac {3} 4 \nb 3 + \nb 4\ra_{\alpha = -\frac 2 3},$ & $\la \nb 2 + \nb 3\ra_{\alpha = -\frac 2 3},$ \\[2mm]
$\alpha = \frac {4\lambda-1}{5-2\lambda}$&$:$ &
$\la \nb 1 + \beta\nb 2 + \frac {4\lambda-1}{5-2\lambda} \nb 3 + \nb 4\ra_{\lambda \neq 1, -2, -\frac 1 2},$ & $\la \nb 2 + \beta \nb 3+\nb 4\ra_{\lambda \neq 1},$  & $\la \beta \nb 3 + \nb 4\ra,$\\[2mm]
  && $\la \nb 1 + \nb 3 + \nb 4\ra_{\lambda =1},$ & $\la \nb 2 + \nb 3 + \nb 4\ra_{\lambda =1},$\\[2mm]
 & & $\la \nb 1 + \beta \nb 2 + \gamma \nb 3+\nb 4\ra_{\lambda =-2},$ & $\la \beta \nb 1 + \nb 2 + \nb 3\ra_{\lambda =-2},$\\
 && $\la \nb 2 -\frac 1 2 \nb 3+ \nb 4\ra_{\lambda =-\frac 1 2},$ &  $\la \nb 2 + \nb 3\ra_{\lambda =-\frac 1 2},$ \\[2mm]
\hline$\lambda \neq 0$&$:$ &  
$\la \nb 2 + \nb 3 \ra_{\lambda\neq -2; \alpha\neq \frac{2\lambda^2+5\lambda-1}{6}},$ & $\la \nb 3\ra, $\\[2mm]
 && $\la \beta \nb 1 +   \nb 2 + \nb 3 \ra_{\lambda =-2, \alpha=-\frac{1}{2}},$ & $\la \nb 1 + \nb 3\ra_{\lambda =-2}.$
 \end{longtable}

The corresponding algebras are:
%\[\T{5}{31}, \dots \T{5}{58}.\]
\begin{longtable}{llllllllll}
   \hline $\T{5}{31}(\beta) $ &  
   $e_1 e_1 = e_2$ & $e_1 e_2= e_3+\beta e_5$ & $e_1 e_3=e_5$ \\
   & $e_2 e_1=e_3$ & $e_2 e_2=e_5$ & $e_3 e_1= e_4$  & $e_3 e_4= e_5$ &&\\

    \hline $\T{5}{32} $ &  
    $e_1 e_1 = e_2$ & $e_1 e_2= e_3+ e_5$ &  $e_2 e_1=e_3$ \\
    & $e_3 e_1= e_4$& $e_3 e_4= e_5$\\

    \hline $\T{5}{33} $ &  
    $e_1 e_1 = e_2$ & $e_1 e_2= e_3$ &  $e_2 e_1=e_3$ \\
    & $e_3 e_1= e_4$& $e_3 e_4= e_5$\\
    
    \hline $\T{5}{34}(\beta) $ &  
    $e_1 e_1 = e_2$ & $e_1 e_2=\beta e_5$ & $e_1 e_3= \frac 1 3 e_4$ \\
    & $e_2 e_1=e_3$ & $e_2 e_2=\frac{2} 3 e_4$      & $e_2 e_3= e_5$ \\
    & $e_2 e_4= e_5$ & $e_3 e_1= e_4$ & $e_3 e_3= e_5$ \\

    \hline $\T{5}{35} $ &  
    $e_1 e_1 = e_2$ & $e_1 e_2=e_5$ & $e_1 e_3= \frac 1 3 e_4$ & $e_2 e_1=e_3$    \\ 
    & $e_2 e_2=\frac{2} 3 e_4$  &  $e_2 e_4= e_5$ & $e_3 e_1= e_4$ & $e_3 e_3= e_5$   \\
    
    \hline $\T{5}{36} $ &  
    $e_1 e_1 = e_2$  & $e_1 e_3= \frac 1 3 e_4$ & $e_2 e_1=e_3$ \\
    & $e_2 e_2=\frac{2} 3 e_4$    &  $e_2 e_4= e_5$  & $e_3 e_1= e_4$     & $e_3 e_3= e_5$ &&\\
    
    \hline $\T{5}{37}(\beta) $ &  
    $e_1 e_1 = e_2$  & $e_1 e_3= - \frac 2 3 e_4+e_5$ & $e_2 e_1=e_3$    & $e_1 e_4=e_5$  \\
    & $e_2 e_2=-\frac{1} 3 e_4+e_5$    & $e_2 e_3= \beta e_5$ & $e_3 e_1= e_4$ & $e_3 e_2= e_5$ & \\
    
    \hline $\T{5}{38}(\beta) $ &  
    $e_1 e_1 = e_2$  & $e_1 e_3= - \frac 2 3 e_4$ & $e_2 e_1=e_3$ & $e_1 e_4=e_5$ \\
    & $e_2 e_2=-\frac{1} 3 e_4$ &  $e_2 e_3= \beta e_5$ & $e_3 e_1= e_4$ & $e_3 e_2= e_5$ & \\
    
    \hline $\T{5}{39}(\beta) $ &  
    $e_1 e_1 = e_2$ & $e_1 e_3= -\frac 1 5 e_4+e_5$ & $e_1 e_4=e_5$ \\
    & $e_2 e_1=e_3$ & $e_2 e_2=\frac{2} {15} e_4+e_5$    & $e_2 e_3= \beta e_5$ \\ 
    & $e_3 e_1= e_4$ & $e_3 e_2= -\frac 2 5 e_5$ & $e_4 e_1= -\frac 7 2 e_5$ \\
    
    \hline $\T{5}{40}(\beta) $ &  
    $e_1 e_1 = e_2$ & $e_1 e_3= -\frac 1 5 e_4$ & $e_1 e_4=e_5$ \\
    & $e_2 e_1=e_3$ & $e_2 e_2=\frac{2} {15} e_4$  & $e_2 e_3= \beta e_5$ \\
    & $e_3 e_1= e_4$ & $e_3 e_2= -\frac 2 5 e_5$  & $e_4 e_1= -\frac 7 2 e_5$ \\

    \hline $\T{5}{41}(\alpha, \beta) $ &  
    $e_1 e_1 = e_2$ & $e_1 e_2=- e_3$ & $e_1 e_3=\alpha e_4+e_5$ \\
    & $e_1 e_4=\beta e_5$ & $e_2 e_1=e_3$  & \multicolumn{2}{l}{$e_2 e_2=\big(\alpha+\frac{2} 3\big) e_4+e_5$} \\
    & \multicolumn{2}{l}{ $e_2 e_3=\big(\frac{3\alpha+1} 3 -\frac{\beta(3\alpha+4)} 3 \big) e_5$} & $e_3 e_1= e_4$ \\
    & \multicolumn{2}{l}{$e_3 e_2=\big(\beta - \frac{3\alpha+1} 2\big) e_5$}    & $e_4e_1 =e_5$ \\
    
    \hline $\T{5}{42}(\alpha, \beta) $ &  
    $e_1 e_1 = e_2$ & $e_1 e_2=- e_3$ & $e_1 e_3=\alpha e_4$ \\
    & $e_1 e_4=\beta e_5$ & $e_2 e_1=e_3$   & $e_2 e_2=\big(\alpha+\frac{2} 3\big) e_4$ \\
    & \multicolumn{2}{l}{$e_2 e_3=\big(\frac{3\alpha+1} 3 -\frac{\beta(3\alpha+4)} 3 \big) e_5$}  & $e_3 e_1= e_4$ \\
    & \multicolumn{2}{l}{$e_3 e_2=\big(\beta - \frac{3\alpha+1} 2\big) e_5$} &  $e_4e_1 =e_5$ \\

    \hline $\T{5}{43}(\alpha, \beta) $ &  
    $e_1 e_1 = e_2$ & $e_1 e_2=- e_3 + e_5$ & $e_1 e_3=\alpha e_4+\beta e_5$ \\
     $\alpha \neq -\frac 2 3$ & $e_1 e_4=-\frac {3\alpha+5} 4 e_5$  & $e_2 e_1=e_3$& 
     \multicolumn{2}{l}{$e_2 e_2=\big(\alpha+\frac{2} 3\big) e_4+\beta e_5$ } \\ 
     &  \multicolumn{2}{l}{$e_2 e_3=\frac{3\alpha^2+13\alpha+8} 4  e_5$}& $e_3 e_1= e_4$    \\
     & $e_3 e_2= -\frac{9\alpha+7} 4 e_5$  & $e_4e_1 =e_5$ \\
     
    \hline $\T{5}{44} $ &  
    $e_1 e_1 = e_2$ & $e_1 e_2= -  e_3+e_5$ & $e_1 e_3=- \frac 2 3  e_4$ \\
    & $e_1 e_4= -\frac 3 4 e_5$ & $e_2 e_1=e_3$& $e_2 e_3= \frac 1 6 e_5$ \\
    & $e_3 e_1 =e_4$         & $e_3 e_2= - \frac 1 4 e_5$ & $e_4 e_1= e_5$    \\
    
    \hline $\T{5}{45} $ &  
    $e_1 e_1 = e_2$ & $e_1 e_2= -  e_3$ & $e_1 e_3=- \frac 2 3  e_4+e_5$ \\
    & $e_1 e_4=  e_5$ & $e_2 e_1=e_3$ & $e_2 e_2= e_5$ \\
    & $e_2 e_3= -\frac 2 3 e_5$ & $e_3 e_1 =e_4$ & $e_3 e_2= e_5$  \\

    \hline $\T{5}{46}(\lambda, \beta) $ &  
    $e_1 e_1 = e_2$ & $e_1 e_2=\lambda e_3+e_5$ & \multicolumn{2}{l}{$e_1 e_3=\frac{4\lambda-1}{5-2\lambda} e_4+\beta e_5$}\\ 
    $\lambda \neq 1, -2, -\frac 1 2, \frac  5 2$ & $e_1 e_4=\frac{4\lambda-1}{5-2\lambda} e_5$ &     $e_2 e_1=e_3$  &   \multicolumn{2}{l}{$e_2 e_2=\frac{(2\lambda+1)(\lambda+2)}{3(5-2\lambda)} e_4+\beta e_5$} \\
    &\multicolumn{2}{l}{$e_2 e_3=\frac {2(2\lambda+1)(3\lambda^2+1)}{3\lambda(5-2\lambda)^2}e_5$}& $e_3 e_1= e_4$ \\
    &   $e_3 e_2=\frac {2\lambda+1}{5-2\lambda}e_5$     &$e_4 e_1= e_5$  \\
    
    \hline $\T{5}{47}(\lambda, \beta) $ &  
    $e_1 e_1 = e_2$ & $e_1 e_2=\lambda e_3$ & $e_1 e_3=\frac{4\lambda-1}{5-2\lambda} e_4+ e_5$  \\
    $\lambda \neq 1, \frac   5 2$ & $e_1 e_4=\beta e_5$  & $e_2 e_1=e_3 $ 
          & \multicolumn{2}{l}{$e_2 e_2=\frac{(2\lambda+1)(\lambda+2)}{3(5-2\lambda)} e_4+ e_5$ }\\
    & \multicolumn{2}{l}{$e_2 e_3=\frac {\beta(4\lambda^2-2\lambda+7)+2(\lambda-1)^2}{3\lambda(5-2\lambda)}e_5$}  & $e_3 e_1= e_4$\\  
    &\multicolumn{2}{l}{$e_3 e_2=\Big(\beta - \frac {2(\lambda-1)}{5-2\lambda}\Big)e_5$ }       & $e_4 e_1= e_5$ \\
    
    \hline $\T{5}{48}(\lambda, \beta) $ &  
    $e_1 e_1 = e_2$ & $e_1 e_2=\lambda e_3$ & $e_1 e_3=\frac{4\lambda-1}{5-2\lambda} e_4$\\ 
     $  \lambda \neq 0, \frac  5 2$ & $e_1 e_4=\beta e_5$   & $e_2 e_1=e_3$   & $e_2 e_2=\frac{(2\lambda+1)(\lambda+2)}{3(5-2\lambda)} e_4$\\ 
     &\multicolumn{2}{l}{$e_2 e_3=\frac {\beta(4\lambda^2-2\lambda+7)+2(\lambda-1)^2}{3\lambda(5-2\lambda)}e_5$} &   $e_3 e_1= e_4$  \\
     & \multicolumn{2}{l}{$e_3 e_2=\Big(\beta - \frac {2(\lambda-1)}{5-2\lambda}\Big)e_5$}  & $e_4 e_1= e_5$ \\
     
    \hline $\T{5}{49} $ &  $e_1 e_1 = e_2$ & $e_1 e_2=e_3+e_5$ & $e_1 e_3= e_4$ \\
    & $e_1 e_4= e_5$ & $e_2 e_1=e_3$ & $e_2 e_2= e_4$ \\
    & $e_2 e_3=e_5$  &  $e_3 e_1= e_4$ &  $e_3 e_2= e_5$   \\ & $e_4 e_1= e_5$ \\
            
    \hline $\T{5}{50} $ &  
    $e_1 e_1 = e_2$ & $e_1 e_2=e_3$ & $e_1 e_3= e_4+e_5$ \\
    & $e_1 e_4= e_5$ &   $e_2 e_1=e_3$   & $e_2 e_2= e_4+e_5$ \\ 
    & $e_2 e_3=e_5$ &  $e_3 e_1= e_4$ & $e_3 e_2= e_5$   &$e_4 e_1= e_5$ \\
    
    \hline $\T{5}{51}(\beta, \gamma) $ &  
    $e_1 e_1 = e_2$ & $e_1 e_2=-2 e_3 + e_5$ & $e_1 e_3= - e_4+\beta e_5$ \\ 
    & $e_1 e_4= \gamma e_5$ &  $e_2 e_1=e_3$ & $e_2 e_2= \beta e_5$   \\ 
    & $e_2 e_3=-\frac{3\gamma+2} 6 e_5$    &  $e_3 e_1= e_4$    & $e_3 e_2= \frac{3\gamma+2} 3 e_5$    & $e_4 e_1= e_5$  \\
    
    \hline $\T{5}{52}(\beta) $ &  
    $e_1 e_1 = e_2$ & $e_1 e_2=-2 e_3 + \beta e_5$ & $e_1 e_3= - e_4+e_5$ \\
    & $e_1 e_4= e_5$ & $e_2 e_1=e_3$ & $e_2 e_2=  e_5$ \\
    &    $e_2 e_3=-\frac{1} 2 e_5$ &  $e_3 e_1= e_4$ &   $e_3 e_2=  e_5$    \\
    
    \hline $\T{5}{53} $ &  
    $e_1 e_1 = e_2$ & $e_1 e_2=-\frac 1 2 e_3$ & $e_1 e_3= -\frac 1 2 e_4+e_5$ & $e_1 e_4= -\frac 1 2 e_5$    \\ 
    & $e_2 e_1=e_3$  & $e_2 e_2= e_5$ & $e_3 e_1= e_4$ & $e_4 e_1= e_5$ \\
    
    \hline $\T{5}{54} $ &  
    $e_1 e_1 = e_2$ & $e_1 e_2=-\frac 1 2 e_3$ & $e_1 e_3= -\frac 1 2 e_4+e_5$ \\
    &  $e_1 e_4= e_5$ & $e_2 e_1=e_3$  &  $e_2 e_2= e_5$ \\
    & $e_2 e_3= - e_5$ & $e_3 e_1= e_4$ & $e_3 e_2= e_5$   \\
   
   	\hline $\T{5}{55}(\lambda, \alpha)_{\lambda\neq -2,0} $ &  
   	$e_1 e_1 = e_2$ & $e_1 e_2=\lambda e_3$ & $e_1 e_3=\alpha e_4+e_5$ \\
   	 $\alpha\neq \frac{2\lambda^2+5\lambda-1}{6}$   & $e_1 e_4= e_5$ &   $e_2 e_1=e_3$ & 
   	   \multicolumn{2}{l}{$e_2 e_2=\big(\alpha+\frac{1-\lambda}3\big) e_4+e_5$ }\\
   	 & $e_2e_3 =  \frac  {3\alpha-2(\lambda-1)}{3\lambda} e_5$ & $e_3 e_1= e_4$     & $e_3e_2 = e_5$& \\

    \hline $\T{5}{56}(\lambda, \alpha) $ &  
    $e_1 e_1 = e_2$ & $e_1 e_2=\lambda e_3$ & $e_1 e_3=\alpha e_4$  \\
     $ \lambda\neq  0$ & $e_1 e_4= e_5$ & $e_2 e_1=e_3$  &  $e_2 e_2=\big(\alpha+\frac{1-\lambda}3\big) e_4$ \\ 
     & $e_2e_3 =  \frac  {3\alpha-2(\lambda-1)}{3\lambda} e_5$ & $e_3 e_1= e_4$   & $e_3e_2 = e_5$& \\
        
    \hline $\T{5}{57}(\beta) $ &  
    $e_1 e_1 = e_2$ & $e_1 e_2=-2 e_3$ & $e_1 e_3=-\frac 1 2 e_4 + e_5$ \\
    & $e_1 e_4= e_5$    &  $e_2 e_1=e_3$   & $e_2 e_2=\frac 1 2 e_4 + e_5$ \\
    & $e_2e_3 = - \frac  {3}{4} e_5$&  $e_3 e_1= e_4$ & $e_3e_2 = e_5$ \\
        
    \hline $\T{5}{58}(\alpha) $ &  
    $e_1 e_1 = e_2$ & $e_1 e_2=-2 e_3+e_5$ & $e_1 e_3=\alpha e_4$ \\
    & $e_1 e_4= e_5$   & $e_2 e_1=e_3$   &  $e_2 e_2=(\alpha+1) e_4$ \\
    & $e_3 e_1= e_4$ & $e_2e_3 =  - \frac  {\alpha+2}{2} e_5$& $e_3e_2 = e_5$  \\
 \hline
\end{longtable}

\subsection{$1$-dimensional  central extensions of $\T {4}{06}(\lambda)$}\label{ext-T_06^4}
	\begin{enumerate}
    \item $\lambda =1.$ Let us use the following notations:
	\begin{align*}
	\nb 1 = \Dl 12, \quad  \nb 2 = \Dl 31, \quad \nb 3 = \Dl 14 + \Dl 23+\Dl 32, \quad \nb 4 =3\Dl 41 - \Dl 22 .
	\end{align*}
Take $\0=\sum\limits_{i=1}^4\af_i\nb i\in {\rm H_T^2}(\T {4}{05}(1)).$
    For the element	$$ \phi=
	\begin{pmatrix}
1 &    0  &  0 &  0\\
x &  1  &  0&  0\\
y &  2x  &  1&  0\\
z &  2y+x+x^2 &  3x & 1
\end{pmatrix}\in\aut{\T {4}{05}(1)},
	$$
we have 
	$$
	\phi^T\begin{pmatrix}
	0 & \alpha_1 & 0 & \alpha_3 \\
    0 & - \alpha_4 & \alpha_3 & 0 \\
    \alpha_2 & \alpha_3 & 0 & 0 \\
    3\alpha_4 & 0 & 0 & 0
	\end{pmatrix} \phi=
	\begin{pmatrix}
    \alpha^{***} & \alpha^{**} + \alpha^{*}+\alpha_1^* & \alpha^{*} & \alpha_3^* \\
               \alpha^{**} &  \alpha^{*} - \alpha_4^* & \alpha_3^* & 0 \\
               \alpha^{*}+ \alpha_2^* & \alpha_3^* & 0 & 0 \\
               3\alpha_4^* & 0 & 0 & 0
	\end{pmatrix},
	$$
where
\[	\begin{array}{rclcrcl}
	\alpha^*_1 &=& \alpha _1- x(2\alpha _2+3\alpha _3+3\alpha _4)+(x^2+2y)(\alpha _3-3\alpha _4) , &&
	\alpha^*_2 &=& \alpha_2-3x(\alpha_3 - 3\alpha_4),\\
    \alpha^*_3 &=& \alpha_3, &&
    \alpha^*_4 &=& \alpha_4.
	\end{array}\]

Since $(\af_3, \af_4) \neq (0,0),$  we have the following cases:

\begin{enumerate}
\item if $\alpha_4 =0,$ then choosing $x=\frac{\af_2}{3\af_3},$ $y=\frac{-\alpha_1 + 2x\af_2 - (x^2-3x)\af_3 }{2\af_3},$ we have the representative $\la \nb 3\ra$.

\item if $\alpha_4 \neq 0,$ then:
\begin{enumerate}
\item if $\alpha_3 \neq 3\alpha_4,$ then choosing $x=\frac{\af_2}{3(\af_3-3\af_4)},$
$y=\frac{-\alpha_1 + 2x \af_2 + (3x-x^2)\af_3  + 3(x^2+x)\af_4}{2(\af_3-3\af_4)},$ we have the representative $\la \beta \nb 3+ \nb 4\ra_{\beta\neq 3}$.

\item if $\alpha_3 = 3\alpha_4,$ then:
\begin{enumerate}
\item if $\alpha_2 \neq -6\alpha_4,$ then choosing $y=\frac{\af_1}{\af_2+6\af_4},$
 we have the representative $\la  \beta \nb 2 + 3 \nb 3+ \nb 4\ra_{\beta\neq - 6}$.

\item if $\alpha_2 = -6\alpha_4,$ then we have the representative $\la \beta \nb 1 - 6\nb 2 + 3 \nb 3 + \nb 4\ra$.
\end{enumerate}
\end{enumerate}

\end{enumerate}

%\begin{longtable}{lllllllll}
%   %\T{5}{1} &:& e_1 e_1 = e_2, & e_1 e_2= e_3+e_4, & e_1 e_3= e_4, & e_1 e_4=e_5, & e_2 e_1=e_3, \\ & & e_2 e_2= e_4, & e_2 e_3= e_5, &e_3 e_1= e_4, &e_3 e_2= e_5;\\
%   	$\T{5}{2}(\af)$ &:& $e_1 e_1 = e_2$, & $e_1 e_2= e_3+e_4$, & $e_1 e_3= e_4$, & $e_1 e_4=\af e_5$, & $e_2 e_1=e_3$, \\ & & $e_2 e_2= e_4-e_5$, & $e_2 e_3= \af e_5$, & $e_3 e_1= e_4$, &$e_3 e_2=\af e_5$,
%   &$e_4 e_1= 3e_5$;\\
%   $\T{5}{3}(\af)$ &:& $e_1 e_1 = e_2$, & $e_1 e_2= e_3+e_4$, & $e_1 e_3= e_4$, & $e_1 e_4=3e_5$, & $e_2 e_1=e_3$, \\ & & $e_3 e_1= e_4$, & $e_2 e_3=3e_5$, & $e_2 e_2= e_4-e_5$, & $e_3 e_2= 3e_5$, & $e_4 e_1= 3e_5$;\\
%       $\T{5}{4}(\af)$ &:& $e_1 e_1 = e_2$, & $e_1 e_2= e_3+e_4+\af e_5$, & $e_1 e_3= e_4$, & $e_1 e_4= 3 e_5$, & $e_2 e_1=e_3$, \\ &  & $e_2 e_2= e_4 -e_5$, &  $e_2 e_3= 3 e_5$, & $e_3 e_1= e_4-6e_5$, & $e_3 e_2=3 e_5$, & $e_4 e_1= 3 e_5$;
%\end{longtable}

\item  $\lambda =-1.$

Let us use the following notations:
	\begin{align*}
	\nb 1 = \Dl 12, \quad \nb 2 = \Dl 13 + \Dl 22, \quad \nb 3 = 3\Dl 14-2 \Dl 23+3\Dl 32, \quad \nb 4 =2 \Dl 13 -3 \Dl 14 +6\Dl 41.
	\end{align*}
Take $\0=\sum\limits_{i=1}^4\af_i\nb i\in {\rm H_T^2}(\T {4}{06}(-1)).$
  For the element
    $$
	\phi=
	\begin{pmatrix}
1 &    0  &  0 &  0\\
x &  1  &  0&  0\\
y &  0  &  1&  0\\
z &  x+\frac{1}{3}y &  0 & 1
\end{pmatrix}\in\aut{\T {4}{06}(-1)},
	$$
	we have 
	$$
	\phi^T\begin{pmatrix}0 & \alpha_1 & \alpha_2 + 2\alpha_4 & 3\alpha_3-3\alpha_4 \\
               0 & \alpha_2 & -2\alpha_3 & 0 \\
               0 & 3\alpha_3 & 0 & 0 \\
               6\alpha_4 & 0 & 0 & 0
	\end{pmatrix} \phi=
	\begin{pmatrix}
    \alpha^{***} & -\alpha^{**} + 3\alpha^{*}+\alpha_1^* & -2 \alpha^{*} + \alpha_2^* + 2\alpha_4^*& 3\alpha_3^*-3\alpha_4^* \\
               \alpha^{**} & \alpha_2^* & -2\alpha_3^* & 0 \\
               3\alpha^{*} & 3\alpha_3^*  & 0 & 0 \\
               6\alpha_4^* & 0 & 0 & 0
	\end{pmatrix},
	$$
where
\[
\begin{array}{rclcrclcrcl}
	\alpha^*_1 &=& \alpha _1+x(2\alpha _2+3\alpha _4)+y(2\alpha _3+\alpha _4), &&    \alpha^*_2 &=& \alpha_2 ,\\
    \alpha^*_3 &=& \alpha_3, &&
    \alpha^*_4 &=& \alpha_4.
	\end{array} \]

Since $(\af_3, \af_4) \neq (0,0),$  we have the following cases:
\begin{enumerate}
\item if $\alpha_4 =0,$ then 
 choosing $y=-\frac{\alpha_1+2 \alpha_2x}{2\af_3},$
  we have the representatives $\la  \beta \nb 2+\nb3 \ra$.

\item if $\alpha_4 \neq 0, $ then: 

\begin{enumerate}
\item if $\alpha_4=-2\af_3,$ $ \af_2=3 \af_3,$ then we have the representative 
\[\la \beta \nb 1 + 3 \nb 2 +\nb 3-2\nb 4\ra .\]
\item if $\alpha_4=-2\af_3,$ $ \af_2\neq 3 \af_3,$ then choosing $x=-\frac {\alpha_1} {2\alpha_2-6\alpha_3},$ we have the representative \[ \la \beta \nb 2 + \nb 3-2\nb 4\ra.\]
\item if $\alpha_4\neq -2\af_3,$ then choosing $y=-\frac {\alpha_1+x(2\alpha_2+3\alpha_4)} {2\alpha_3+\alpha_4},$ we have the representative \[ \la \beta \nb 2 + \gamma \nb 3+ \nb 4\ra.\]
    
\end{enumerate}
\end{enumerate}

\item  $\lambda =-\frac{1}{2}.$
Let us use the following notations:
	\begin{align*}
	\nb 1 = \Dl 12, \quad \nb 2 = \Dl 13 + \Dl 22, \quad \nb 3 = \Dl 14- \Dl 23+\Dl 32, \quad \nb 4 =2 \Dl 13 -3 \Dl 14 +6\Dl 41.
	\end{align*}
Take $\0=\sum\limits_{i=1}^4\af_i\nb i\in {\rm H_T^2}(\T {4}{06}(-\frac{1}{2})).$
  For the element
    $$
	\phi=
	\begin{pmatrix}
1 &    0  &  0 &  0\\
x &  1  &  0&  0\\
y &  \frac{1}{2}x  &  1&  0\\
z &  x+\frac{1}{2}y &  \frac{1}{2}x  & 1
\end{pmatrix}\in\aut{\T {4}{06}(-\frac{1}{2})},
	$$
	we have
	$$
	\phi^T\begin{pmatrix}
0 & \alpha_1 & \alpha_2 + 2\alpha_4 & \alpha_3-3\alpha_4 \\
               0 & \alpha_2 & -\alpha_3 & 0 \\
               0 & \alpha_3 & 0 & 0 \\
               6\alpha_4 & 0 & 0 & 0
	\end{pmatrix} \phi=
	\begin{pmatrix}
    \alpha^{***} & -\frac{1}{2}\alpha^{**} + \alpha^{*}+\alpha_1^* & -\frac{1}{2} \alpha^{*} + \alpha_2^* + 2\alpha_4^*& \alpha_3^*-3\alpha_4^* \\
               \alpha^{**} & \alpha_2^* & -\alpha_3^* & 0 \\
               \alpha^{*} & \alpha_3^*  & 0 & 0 \\
               6\alpha_4^* & 0 & 0 & 0
	\end{pmatrix},
	$$
where
\[
\begin{array}{rclcrcl}
	\alpha^*_1 &=& \alpha_1+2x(\alpha _2-\alpha _4)-(\frac{1}{4}x^2-y)\alpha _3, &&     \alpha^*_2 &=& \alpha_2 ,\\
    \alpha^*_3 &=& \alpha_3, &&    \alpha^*_4 &=& \alpha_4.
	\end{array} \]

Since $(\af_3, \af_4) \neq (0,0),$  we have the following cases:
\begin{enumerate}
\item if $\alpha_4 =0,$ then
 choosing $y = -\frac{4\alpha_1+8x\alpha _2-x^2\alpha_3}{4\alpha_3},$
  we have the representatives $\la  \beta \nb 2+\nb3 \ra$.

\item if $\alpha_4 \neq 0, $ then:

\begin{enumerate}
\item if $\alpha_3=0,$ $ \af_2= \af_4,$ then we have the representative $ \la \beta \nb 1 + \nb 2 +\nb 4\ra $.
\item if $\alpha_3=0,$ $ \af_2\neq \af_4,$ then choosing $x = -\frac{\alpha_1}{2(\alpha _2-\alpha _4)},$ we have the representative $ \la \beta \nb 2 + \nb 4\ra $.
\item if $\alpha_3\neq 0,$ then choosing $y = -\frac{4\alpha_1+8x(\alpha _2-\alpha _4)-x^2\alpha_3}{4\alpha_3},$ we have the representative $ \la \beta \nb 2 + \gamma \nb 3+ \nb 4\ra $.

\end{enumerate}
\end{enumerate}

\item  $\lambda \neq 0; 1; -1; -\frac{1}{2}.$
Let us use the following notations:
	\begin{align*}
	\nb 1 = \Dl 12, \quad \nb 2 = \Dl 13 + \Dl 22, \quad \nb 3 = 6\lambda \Dl 14 +(2\lambda^2+\lambda+3) \Dl 23 + 6\lambda\Dl 32.
	\end{align*}
Take $\0=\sum\limits_{i=1}^3\af_i\nb i\in {\rm H_T^2}(\T {4}{06}(\lambda)).$
For the element
	$$
	\phi=
	\begin{pmatrix}
1 &    0  &  0 &  0\\
x &  1  &  0&  0\\
y &  (\lambda+1)x  &  1&  0\\
z &  \frac{2 \lambda ^2+3 \lambda +1}{6}x^2+x+\frac{2 \lambda ^2+5 \lambda +5}{6}y &  \frac{2 \lambda ^2+9 \lambda +7}{6}x & 1
\end{pmatrix}\in\aut{\T {4}{06}(\lambda)},
	$$
	we have
	$$
	\phi^T\begin{pmatrix}
	0 & \alpha_1 & \alpha_2 & 6\lambda \alpha_3 \\
               0 & \alpha_2 & (2\lambda^2+\lambda +3)  \alpha_3 & 0 \\
               0 & 6\lambda \alpha_3 & 0 & 0 \\
               0 & 0 & 0 & 0
	\end{pmatrix} \phi=
	\begin{pmatrix}
    \alpha^{***} & \lambda \alpha^{**} + \alpha^{*}+\alpha_1^* & \frac{2\lambda^2+5\lambda-1}{6} \alpha^{*} + \alpha_2^*& 6\lambda \alpha_3^* \\
               \alpha^{**} & \frac{2\lambda^2+3\lambda+1}{6} \alpha^{*} + \alpha_2^* & (2\lambda^2+\lambda +3)  \alpha_3^* & 0 \\
               \alpha^{*} & 6\lambda \alpha_3^* & 0 & 0 \\
               0 & 0 & 0 & 0
	\end{pmatrix},
	$$
where
\begin{longtable}{rcl}
	$\alpha^*_1$ &$=$&$ \alpha _1+2x\alpha _2+\big((-2 \lambda ^3+5 \lambda +3) x^2+4 y \lambda  (\lambda +2) \big)\alpha _3,$\\ 
	$\alpha^*_2$ &$=$&$ \alpha _2+3x(2 \lambda ^2+3 \lambda +1)\alpha _3,$\\
	$\alpha^*_3$ &$=$& $\alpha_3.$
	\end{longtable}

Since $\af_3 \neq 0$  we have the following cases:

\begin{enumerate}
\item if $\lambda \neq - 2,$ then choosing 
\[x= -\frac{\af_2}{3(2\lambda^2+3\lambda+1)\af_3},\ \ y= -\frac{1} {4\lambda (\lambda+2)\alpha_3} \Big(\alpha_1 + 2x\alpha_2  - (2\lambda^3 -5\lambda -3)x^2\alpha_3 \Big),\]
  we have the representatives $\la \nb 3 \ra.$

\item  if $\lambda = -2,$ then choosing 
 $x= -\frac{\af_2}{9\af_3},$ we have the representative $\la \beta \nb 1 + \nb 3 \ra.$

\end{enumerate}

\end{enumerate}

Summarizing all cases, we have the following distinct orbits:
\begin{longtable}{llllllll}
$\lambda=1$&$:$ &  $\la \beta \nb 3+ \nb 4\ra,$ &  $\la \beta \nb 2 +  3 \nb 3+ \nb 4\ra_{\beta\neq 0},$
& $\la \beta \nb 1 -  6\nb 2 + 3 \nb 3 +  \nb 4\ra_{\beta\neq 0}$ \\[2mm]

$\lambda=- 1$&$:$ & $\la \beta \nb 2+ \nb 3\ra_{\beta\neq 0},$ &
$\la \beta\nb 1+ 3\nb 2+ \nb 3 -2\nb 4\ra_{\beta\neq 0},$ & $\la \beta\nb 2+ \gamma \nb 3 +\nb 4\ra$ & \\[2mm]

$\lambda=- \frac 1 2$&$:$ & $\la \beta \nb 2+ \nb 3\ra_{\beta\neq 0},$ &
$\la \beta\nb 1+ \nb 2+ \nb 4\ra_{\beta\neq 0},$ & 
$\la \beta\nb 2+ \gamma \nb 3 +\nb 4\ra$  \\[2mm]
$\lambda=- 2$&$:$ & $\la \beta \nb 1+ \nb 3\ra_{\beta\neq 0}$ \\[2mm]
$\lambda\neq 0$&$:$ & $\la \nb 3\ra.$ \\[2mm]
\end{longtable}

The corresponding algebras are:
%\[\T{5}{59}, \dots \T{5}{69}.\]
\begin{longtable}{llllllllllllllll}
    
    \hline $\T{5}{59}(\beta) $ &  $e_1 e_1 = e_2$ & $e_1 e_2= e_3+e_4$ & $e_1 e_3= e_4$ & $e_1 e_4= \beta e_5$  \\
    & $e_2 e_1=e_3$ & $e_2 e_2= e_4 -e_5$ &  $e_2 e_3= \beta e_5$ & $e_3 e_1= e_4$ \\
    & $e_3 e_2= \beta e_5$ & $e_4 e_1= 3 e_5$\\

    \hline $\T{5}{60}(\beta)_{\beta\neq0} $ &  
    $e_1 e_1 = e_2$ & $e_1 e_2= e_3+e_4$ & $e_1 e_3= e_4$ & $e_1 e_4= 3 e_5$ \\ 
    & $e_2 e_1=e_3$ & $e_2 e_2= e_4 -e_5$ &  $e_2 e_3= 3 e_5$ & $e_3 e_1= e_4 +  \beta e_5$ \\
    & $e_3 e_2= 3 e_5$ & $e_4 e_1= 3 e_5$\\
    
    \hline $\T{5}{61}(\beta)_{\beta\neq0} $ &  
    $e_1 e_1 = e_2$ & $e_1 e_2= e_3+e_4+\beta e_5$ & $e_1 e_3= e_4$ & $e_1 e_4= 3 e_5$ \\
    & $e_2 e_1=e_3$ & $e_2 e_2= e_4 -e_5$ &  $e_2 e_3= 3 e_5$ & $e_3 e_1= e_4 -6 e_5$ \\
    & $e_3 e_2= 3 e_5$ & $e_4 e_1= 3 e_5$\\
    
    \hline $\T{5}{62}(\beta)_{\beta\neq0} $ &  
    $e_1 e_1 = e_2$ & $e_1 e_2= - e_3+e_4$ & $e_1 e_3= -\frac 2 3 e_4+\beta e_5$ & $e_1 e_4=  3e_5$ \\
    & $e_2 e_1=e_3$    & $e_2 e_2= \beta e_5$ &  $e_2 e_3= -2 e_5$ & $e_3 e_1= e_4$ \\
    & $e_3 e_2= 3 e_5$ & \\
    
    \hline $\T{5}{63}(\beta)_{\beta\neq0} $ &  
    $e_1 e_1 = e_2$ & $e_1 e_2= - e_3+e_4 + \beta e_5$ & $e_1 e_3= -\frac 2 3 e_4 - e_5$ & $e_1 e_4=  9 e_5$ \\ 
    & $e_2 e_1=e_3$ & $e_2 e_2= 3 e_5$ &  $e_2 e_3= - 2 e_5$ & $e_3 e_1= e_4$ \\
    & $e_3 e_2= 3 e_5$ & $e_4 e_1= -12 e_5$\\
    
    \hline $\T{5}{64}(\beta, \gamma) $ &  
    $e_1 e_1 = e_2$ & $e_1 e_2= - e_3+e_4$ & $e_1 e_3= -\frac 2 3 e_4 +(\beta+2) e_5$ & $e_1 e_4=3(\gamma-1) e_5$ \\
    & $e_2 e_1=e_3$ & $e_2 e_2= \beta e_5$ &  $e_2 e_3= - 2\gamma e_5$ & $e_3 e_1= e_4$ \\
    & $e_3 e_2= 3\gamma e_5$ & $e_4 e_1= 6 e_5$\\
    
    \hline $\T{5}{65}(\beta)_{\beta\neq0} $ &  
    $e_1 e_1 = e_2$ & $e_1 e_2= -\frac 1 2  e_3+e_4$ & $e_1 e_3= -\frac 1 2 e_4+\beta e_5$ & $e_1 e_4=  e_5$ \\ 
    & $e_2 e_1=e_3$    & $e_2 e_2= \beta e_5$ &  $e_2 e_3= - e_5$ & $e_3 e_1= e_4$ \\ 
    & $e_3 e_2= e_5$&\\
    
    \hline $\T{5}{66}(\beta)_{\beta\neq0} $ &  
    $e_1 e_1 = e_2$ & $e_1 e_2= -\frac 1 2  e_3+e_4 + \beta e_5$ & $e_1 e_3= -\frac 1 2 e_4+3 e_5$ & $e_1 e_4= -3e_5$ &\\ 

    &$e_2 e_1=e_3$ & $e_2 e_2= e_5$  & $e_3 e_1= e_4$ & $e_4 e_1= 6e_5$&\\
    
    \hline $\T{5}{67}(\beta, \gamma) $ &  
    $e_1 e_1 = e_2$ & $e_1 e_2= -\frac 1 2  e_3+e_4 $ & $e_1 e_3= -\frac 1 2 e_4+ (\beta+2) e_5$ & $e_1 e_4= (\gamma-3)e_5$ \\
    & $e_2 e_1=e_3$ & $e_2 e_2= \beta e_5$ & $e_2 e_3= -\gamma e_5$ & $e_3 e_1= e_4$ \\
    & $e_3 e_2= \gamma e_5$  & $e_4 e_1= 6e_5$\\
    
    \hline $\T{5}{68}(\beta)_{\beta\neq 0} $ &  
    $e_1 e_1 = e_2$ & $e_1 e_2= -2 e_3+e_4+\beta e_5$ & $e_1 e_3= -\frac 1 2e_4$ & $e_1 e_4= 4 e_5$\\
    & $e_2 e_1=e_3$ & $e_2 e_2= \frac 1 2 e_4$ &  $e_2 e_3= -3 e_5$ & $e_3 e_1= e_4$ \\ 
    & $e_3 e_2= 4e_5$&\\
    
    \hline $\T{5}{69}(\lambda) $ &  
    $e_1 e_1 = e_2$ & $e_1 e_2=\lambda e_3+e_4$ & $e_1 e_3= \frac{2\lb^2  + 5\lb  - 1}{6}e_4$ & $e_1 e_4= 6\lambda e_5$ \\ 
    & $e_2 e_1=e_3$ & $e_2 e_2= \frac{(2\lb+1)(\lb+1)}{6}e_4$ & $e_2 e_3= (2\lb^2+\lb+3) e_5$ & $e_3 e_1= e_4$  \\ 
    & $e_3 e_2= 6\lambda e_5$&\\
\hline  \end{longtable}

\subsection{$1$-dimensional  central extensions of $\T {4}{07}(\lambda)$}\label{ext-T_07^4}
	
\begin{enumerate}

\item $\lambda\neq -1.$ Let us use the following notations:
	\begin{align*}
	\nb 1 = \Dl 12, \quad \nb 2 = (1-\lambda)\Dl 22 + \Dl 31, \quad \nb 3 = \lambda \Dl 14 + \Dl 23.
	\end{align*}
	
Take $\0=\sum\limits_{i=1}^3\af_i\nb i\in {\rm H_T^2}(\T {4}{07}(\lambda)).$
	For the element
	$$
	\phi=
	\begin{pmatrix}
x &    0  &  0 &  0\\
0 &  x^2  &  0&  0\\
y &   0  &  x^3&  0\\
z &  xy  &   0 & x^4
\end{pmatrix}\in\aut{\T {4}{07}(\lambda)},
	$$
we have
	$$
	\phi^T\begin{pmatrix}
	0      &  \af_1     &  0 & \lambda \af_3\\
	0      & (1-\lambda)\af_2  & \af_3  & 0    \\
    3\af_2 & 0 & 0  & 0    \\
    0      & 0  & 0  & 0
	\end{pmatrix} \phi=
	\begin{pmatrix}
    \af^{**}      &  \af^*_1+\af^{*}     & 0  & \af^*_3\\
	\af^{*}      & (1-\lambda)\af^*_2  & \af^*_3  & 0    \\
    3\af^*_2  & 0  & 0  & 0    \\
    0     & 0  & 0  & 0
	\end{pmatrix},
	$$
	where
\[
\begin{array}{rclcrclcrcl}
	\af^*_1 &=& x^3\af_1,&&
	\af^*_2 &=& x^4\af_2,&&
	\af^*_3 &=& x^5\af_3.
	\end{array} \]

Since $\af_3\ne 0,$ then we have the representatives 
\[\la \alpha \nb 1 + \nb 2 +  \nb 3\ra_{\lambda\neq -1}, \
\la  \nb 1 + \nb 3\ra_{\lambda\neq -1} \mbox{ and }\la  \nb 3 \ra_{\lambda\neq -1}\]
 depending on whether $\alpha_1=0,$ $\alpha_2=0$ or not.

\item $\lambda= -1.$ Let us use the following notations:
	\begin{align*}
	\nb 1 = \Dl 12, \quad \nb 2 = 2\Dl 22 + 3 \Dl 31, \quad \nb 3 = \Dl 14 - \Dl 23, \quad \nb 4 =\Dl 41 + 2\Dl 23-3\Dl 32.
	\end{align*}
Take $\0=\sum\limits_{i=1}^4\af_i\nb i\in {\rm H_T^2}(\T {4}{07}(-1)).$
	For the element
	$$
	\phi=
	\begin{pmatrix}
x &    0  &  0 &  0\\
0 &  x^2  &  0&  0\\
y &   0  &  x^3&  0\\
z &  xy  &   0 & x^4
\end{pmatrix}\in\aut{\T {4}{07}(-1)},
	$$
we have 
	$$
	\phi^T\begin{pmatrix}
	0      &  \af_1     & 0  & \af_3\\
	0      & 2\af_2  & 2\af_4-\af_3  & 0    \\
    3\af_2 & -3\af_4 & 0  & 0    \\
    \af_4      & 0  & 0  & 0
	\end{pmatrix} \phi=
	\begin{pmatrix}
    \af^{**}      &  \af^*_1+ \af^{*}     & 0  & \af^*_3\\
	\af^{*}      & 2\af^*_2  & 2\af^*_4-\af^{*}_3  & 0    \\
    3\af^{*}_2  & -3\af^{*}_4  & 0  & 0    \\
    \af^{*}_4     & 0  & 0  & 0
	\end{pmatrix},
	$$
where
\[
\begin{array}{rclcrclcrcl}
	\af^*_1 &=& x^3\af_1, &&
	\af^*_2 &=& x^4\af_2,\\
	\af^*_3 &=& x^5\af_3, &&
    \af^*_4 &=& x^5\af_4.
	\end{array} \]

Since $(\af_3, \af_4) \neq (0,0),$  we have the following cases:
\begin{enumerate}
\item if $\alpha_4 =0,$ then and we have the representatives 
\[\la \alpha \nb 1 + \nb 2 +  \nb 3\ra_{\lambda= -1}, \ 
\la  \nb 1 + \nb 3 \ra_{\lambda= -1} \mbox{ and }\la  \nb 3 \ra_{\lambda= -1}\]
 depending on whether $\alpha_1=0,$ $\alpha_2=0$ or not.

\item if $\alpha_4 \neq 0,$  then we have three distinct orbits with family of representatives
    \[\la \nb 1 +\alpha \nb 2+\beta\nb 3+ \nb 4 \ra, \   \la  \nb 2+ \alpha \nb 3+ \nb 4 \ra \mbox{ and }\la \alpha \nb 3+ \nb 4 \ra.\]
\end{enumerate}

\end{enumerate}
		
Summarizing, we have the following representatives of distinct orbits:
\[ \la \alpha \nb 1 + \nb 2 +  \nb 3\ra,\quad \la  \nb 1 + \nb 3\ra, \quad \la  \nb 2 \ra,\quad
\la \nb 1 +\alpha \nb 2+\beta\nb 3+ \nb 4 \ra, \quad \la  \nb 2+ \alpha \nb 3+ \nb 4 \ra,\quad \la \alpha \nb 3+ \nb 4 \ra.\]

Hence, we have following $5$-dimensional algebras:
\begin{longtable}{lllllllllll}
    \hline $\T{5}{70}(\lambda, \alpha) $ &  $e_1 e_1 = e_2$ & $e_1 e_2=\lambda e_3+\alpha e_5$ & $e_1 e_3= e_4$ & $e_1 e_4=\lambda e_5$&\\
        & $e_2 e_1= e_3$& $e_2 e_2 = e_4+(1-\lambda)e_5$ & $e_2 e_3 = e_5$ & $e_3 e_1=3 e_5$&\\
     \hline $\T{5}{71}(\lambda) $ &  $e_1 e_1 = e_2$ & $e_1 e_2=\lambda e_3 + e_5$ & $e_1 e_3= e_4$ & $e_1 e_4=\lambda e_5$ &\\
        & $e_2 e_1= e_3$& $e_2 e_2 = e_4$ & $e_2 e_3 = e_5$&&\\
     \hline $\T{5}{72}(\lambda) $ &  $e_1 e_1 = e_2$ & $e_1 e_2=\lambda e_3$ & $e_1 e_3= e_4$ & $e_1 e_4=\lambda e_5$ &\\
        & $e_2 e_1= e_3$& $e_2 e_2 = e_4$ & $e_2 e_3 = e_5$&&\\
     \hline $\T{5}{73}(\alpha,\beta) $ &  $e_1 e_1 = e_2$ & $e_1 e_2=- e_3+e_5$ & $e_1 e_3= e_4$ & $e_1 e_4=\alpha e_5$ & $e_2 e_1= e_3$\\
        & $e_2 e_2 = e_4+2\beta e_5$ & $e_2 e_3 = (2-\alpha)e_5$ & $e_3 e_1=3\beta e_5$ & $e_3 e_2=-3 e_5$ & $e_4 e_1= e_5$\\
     \hline $\T{5}{74}(\alpha) $ &  $e_1 e_1 = e_2$ & $e_1 e_2=- e_3$ & $e_1 e_3= e_4$ & $e_1 e_4=\alpha e_5$ & $e_2 e_1= e_3$\\
        & $e_2 e_2 = e_4+2 e_5$ & $e_2 e_3 = (2-\alpha)e_5$ & $e_3 e_1= 3 e_5$ & $e_3 e_2=-3 e_5$ & $e_4 e_1= e_5$\\
     \hline $\T{5}{75}(\alpha) $ &  $e_1 e_1 = e_2$ & $e_1 e_2=- e_3$ & $e_1 e_3= e_4$ & $e_1 e_4=\alpha e_5$ & $e_2 e_1= e_3$\\
        & $e_2 e_2 = e_4$ & $e_2 e_3 = (2-\alpha)e_5$ & $e_3 e_2=-3 e_5$ & $e_4 e_1= e_5$ &\\
        \hline
 \end{longtable}
 
%\[\T 5{70}(\lambda, \alpha), \quad \T 5{71}(\lambda), \quad \T 5{72}(\lambda), \quad \T 5{73}(\alpha,\beta), \quad 
%\T 5{74}(\alpha), \quad \T 5{75}(\alpha).\]

\subsection{$1$-dimensional  central extensions of $\T {4}{08}$}\label{ext-T_08^4}
	Let us use the following notations:
	\begin{align*}
	\nb 1 = \Dl 12, \quad \nb 2 = \Dl 13 + \Dl 22, \quad \nb 3 =\Dl 22+3\Dl 31, \quad  \nb 4 =\Dl 24 .
	\end{align*}
Take $\0=\sum\limits_{i=1}^4\af_i\nb i\in {\rm H_T^2}(\T {4}{08}).$
For the element
	$$
	\phi=
	\begin{pmatrix}
x &    0  &  0 &  0\\
0 &  x^2  &  0&  0\\
y &  0  &  x^3&  0\\
z &  0 &  x^2y & x^5
\end{pmatrix}\in\aut{\T {4}{08}},
	$$
we have 
	$$
	\phi^T\begin{pmatrix}
	0      & \af_1     & \af_2  & 0\\
	0      & \af_2 + \af_3  & 0  & \af_4    \\
    3\af_3 & 0 & 0  & 0    \\
    0      & 0  & 0  & 0
	\end{pmatrix} \phi=
	\begin{pmatrix}
    \af^{***}      &  \af^{*}_1     & \af^*_2  & 0\\
	\af^{**}      & \af^*_2+\af^*_3  & \af^{*}  & \af^*_4    \\
    3\af^*_3  & 0  & 0  & 0    \\
    0    & 0  & 0  & 0
	\end{pmatrix},
	$$
where
\[
\begin{array}{rclcrclcrcl}
	\af^*_1 &=& x^3\af_1, &&
	\af^*_2 &=& x^4\af_2,\\
	\af^*_3 &=& x^4\af_3, &&
    \af^*_4 &=& x^7\af_4.
	\end{array} \]

Since $\af_4\ne 0,$ we have the following distinct orbits with representatives \[\la \nb 1 + \alpha \nb 2 + \beta \nb 3 +  \nb 4\ra, \ 
\la  \nb 2 + \alpha \nb 3 +  \nb 4\ra,  \ 
\la  \nb 3 + \nb 4\ra, \mbox{ and }\la  \nb 4 \ra .\]

Hence, we have following $5$-dimensional algebras:

\begin{longtable}{llllllllllllllllll}
   \hline $\T{5}{76}(\alpha,\beta) $ &  $e_1 e_1 = e_2$ & $e_1 e_2=e_5$ & $e_1 e_3=\alpha e_5$ & $e_2 e_1=e_3$ & \multicolumn{2}{l}{$e_2 e_2= (\alpha+\beta)e_5$}\\& $e_2 e_3 = e_4$ & $e_2 e_4 =  e_5$ & $e_3 e_1=3\beta e_5$&\\
     \hline $\T{5}{77}(\alpha) $ &  $e_1 e_1 = e_2$ &  $e_1 e_3= e_5$ & $e_2 e_1=e_3$ & \multicolumn{2}{l}{$e_2 e_2= (\alpha+1)e_5$} &\\
        & $e_2 e_3 = e_4$ & $e_2 e_4 =  e_5$ & $e_3 e_1=3\alpha e_5$ &&\\
     \hline $\T{5}{78} $ &  $e_1 e_1 = e_2$ &  $e_2 e_1=e_3$ & $e_2 e_2= e_5$ & $e_2 e_3 = e_4$ & $e_2 e_4 =  e_5$ & $e_3 e_1=3 e_5$ &&\\
     \hline $\T{5}{79} $ &  $e_1 e_1 = e_2$ &  $e_2 e_1=e_3$ & $e_2 e_3 = e_4$ & $e_2 e_4 =  e_5$&\\ 
     \hline
     \end{longtable}

%\[\T 5{76}(\alpha,\beta), \quad \T 5{77}(\alpha), \quad \T 5{78}, \quad \T 5{79}.\]

\
\subsection{Classification theorem}
Summarizing results of the present sections,
 we have the following theorem.

\begin{theoremB}\label{main-alg}
Let $\mathbf A$ be a complex $5$-dimensional one-generated nilpotent terminal algebra. 
Then $\mathbf A$ is isomorphic to one of the algebras $\T 5{01}-\T 5{79}$  found in Table C (see Appendix).
\end{theoremB}

\section*{Appendix}

{\tiny 
\begin{longtable}{lllllllll}

\multicolumn{8}{c}{ \mbox{ {\bf \large{Table A.}}
{\it \large{The list of $4$-dimensional one-generated nilpotent terminal algebras.}}}} \\
\multicolumn{8}{c}{  } \\\hline
    $\T 4{01}$&:& $e_1 e_1 = e_2$ & $e_1 e_2=e_4$ & $e_2 e_1=e_3$\\\hline
    $\T{4}{02}(\alpha)$&:& $e_1e_1 = e_2$& $e_1e_2 = e_3$& $e_1e_3 = \alpha e_4$& $e_2e_2 = e_4$& $e_3e_1 = -3e_4$\\\hline
    $\T 4{03}$  &: &  $e_1e_1 = e_2$ & $e_1e_2 = e_3$ &  $e_1e_3 = e_4$ \\\hline
    $\T{4}{04}$&:& $e_1e_1 = e_2$ & $e_1e_2 = e_3$& $e_1e_3 = e_4$& $e_2e_1 = e_4$ \\\hline
    $\T{4}{05}(\lambda, \alpha)$ &:& $e_1 e_1 = e_2$ & $e_1 e_2=\lambda e_3$ & $e_1 e_3=\alpha e_4$ & $e_2 e_1=e_3$ & $e_2 e_2=\big(\alpha+\frac{1-\lambda}3\big) e_4$ & $e_3 e_1= e_4$\\\hline
    $\T{4}{06}(\lambda)$ &:& $e_1 e_1 = e_2$ & $e_1 e_2=\lambda e_3+e_4$ & $e_1 e_3= \frac{2\lb  + 5\lb  - 1}{6}e_4$ & $e_2 e_1=e_3$ & $e_2 e_2= \frac{(2\lb  + 1)(\lb  + 1)}{6}e_4$ & $e_3 e_1= e_4$\\\hline
   	$\T{4}{07}(\lambda)$ &:& $e_1 e_1 = e_2$ & $e_1 e_2=\lambda e_3$ & $e_1 e_3= e_4$ & $e_2 e_1=e_3$ & $e_2 e_2= e_4$\\\hline
    $\T 4{08}$   &:& $e_1 e_1 = e_2$ & $e_2 e_1=e_3$    & $e_2e_3 = e_4$\\\hline
    $\T 4{09}$   &:& $e_1 e_1 = e_2$ & $e_1 e_2=e_4$ & $e_2 e_1=e_3$ & $e_2e_3 = e_4$\\\hline
    $\T 4{10}(\af)$   &:& $e_1 e_1 = e_2$ & $e_1 e_2=\af e_4$ &  $e_2 e_1=e_3$ & $e_2e_2=\frac 1 3 e_4$ & $e_2e_3 = e_4$ & $e_3e_1=e_4$\\\hline

\end{longtable}
}

\begin{longtable}{|lcl|}
 
 \multicolumn{3}{c}{ \mbox{ {\bf {Table B.}}
{\it {Cohomology groups of $4$-dimensional one-generated nilpotent terminal algebras.}}}} \\
\multicolumn{3}{c}{  } \\
\hline
 
 \hline
${\rm B_T^2}(\T {4}{01})$ &$=$&
$\left\langle \Dt 11, \Dt 12, \Dt 21 \right\rangle$\\

${\rm Z_T^2}(\T {4}{01})$ &$=$&
$\left\langle
     \Dt 11, \Dt 12, \Dt 13 +\Dt 22, \Dt 14,     \Dt 21, \Dt 22 + 3\Dt 31, \Dt 23, \Dt 31+ \Dt 41
\right\rangle$\\

${\rm H_T^2}(\T {4}{01})$ &$=$&
$\left\langle
     \Dl 13+\Dl 22, \Dl 14,     \Dl 22 + 3\Dl 31, \Dl 23, \Dl 31 + \Dl 41
\right\rangle$ \\
\hline

${\rm B_T^2}(\T {4}{02}(\alpha)_{\alpha\neq6})$&$=$&$
\left\langle  \Dt 11, \Dt 12,  \alpha \Dt 13+ \Dt 22 - 3 \Dt 13 \right\rangle$\\

${\rm Z_T^2}(\T {4}{02}(\alpha)_{\alpha\neq6})$&$=
$&$\left\langle
     \Dt 11, \Dt 12, \Dt 13, \Dt 14 + 2\Dt 23-3\Dt 32,
     \Dt 21, \Dt 22 - 3\Dt 31
\right\rangle
$\\

${\rm H_T^2}(\T {4}{02}(\alpha)_{\alpha\neq6})$&$=
$&$\left\langle
      \Dl 13, \Dl 14 + 2\Dl 23-3\Dl 32, \Dl 21
\right\rangle$\\
\hline

${\rm B_T^2}(\T {4}{02}(6))$&$=$&
$\left\langle  \Dt 11, \Dt 12,  6 \Dt 13+ \Dt 22 - 3 \Dt 13 \right\rangle$\\

${\rm Z_T^2}(\T {4}{02}(6))$&$=$&
$\left\langle
     \Dt 11, \Dt 12, \Dt 13, \Dt 14 + 2\Dt 23-3\Dt 32,
     \Dt 21, \Dt 22 - 3\Dt 31, \Dt 41 + \Dt 23-3\Dt 32
\right\rangle$\\

${\rm H_T^2}(\T {4}{02}(6))$&$=$&
$\left\langle
      \Dl 13, \Dl 14 + 2\Dl 23-3\Dl 32, 
      \Dl 21,  \Dl 41 + \Dl 23-3\Dl 32
\right\rangle$
\\
\hline

${\rm B_T^2}(\T {4}{03})$&$=$
&
$\left\langle   \Dt 11,  \Dt 12,  \Dt 13   \right\rangle$\\

${\rm Z_T^2}(\T {4}{03})$&$=$
&
$\left\langle \Dt 11, \Dt 12 , \Dt 13, \Dt 14, \Dt 21, \Dt 22-3\Dt 31   \right\rangle$\\

${\rm H_T^2}(\T {4}{03})$&$=$
&
$\left\langle    \Dl 14,  \Dl 21, \Dl 22-3\Dl 31 \right\rangle$
\\
\hline

${\rm B_T^2}(\T {4}{04})$&$=$
&
$\left\langle \Dt 11, \Dt 12, \Dt 13+ \Dt 21 \right\rangle$\\

${\rm Z_T^2}(\T {4}{04})$&$=$
&
$\left\langle
     \Dt 11, \Dt 12, \Dt 13, \Dt 14 + 3\Dt 31,
     \Dt 21, \Dt 22 - 3\Dt 31
\right\rangle$\\

${\rm H_T^2}(\T {4}{04})$&$=$
&
$\left\langle
     \Dl 13, \Dl 14 + 3 \Dl 31, \Dl 22 - 3\Dl 31
\right\rangle$\\ \hline

${\rm B_T^2}(\T {4}{05}(1, 0))$&$=$
&
$\left\langle  \Dt 11, \Dt 12 + \Dt 21, \Dt 31  \right\rangle$\\

${\rm Z_T^2}(\T {4}{05}(1, 0))$&$=$
&
$\left\langle
     \Dt 11, \Dt 12, \Dt 13+\Dt 22, \Dt 21, \Dt 14 + \Dt 32, \Dt 31, \Dt 34
\right\rangle$\\

${\rm H_T^2}(\T {4}{05}(1, 0))$&$=$&
$\left\langle
      \Dl 12, \Dl 13+\Dl 22,  \Dl 14 + \Dl 32,
      \Dl 34
\right\rangle$\\
\hline

${\rm B_T^2}( \T {4}{05}(0, \frac 1 3)) $&$=$
&
$\left\langle 
\Dt 11, \Dt 21,   \Dt 13+ 2 \Dt 22 + 3 \Dt 31 \right\rangle$\\

${\rm Z_T^2}( \T {4}{05}(0, \frac 1 3)) $&$=$&
$\left\langle
     \Dt 11, \Dt 12, \Dt 13+\Dt 22, \Dt 21, 
        \Dt 22 + 3\Dt 31, \Dt 23, \Dt 24 + \Dt 33
\right\rangle$\\

${\rm H_T^2}( \T {4}{05}(0, \frac 1 3)) $&$=$
&
$\left\langle
\Dl 12, \Dl 13 + \Dl 22,  \Dl 23, \Dl 24 + \Dl 33
\right\rangle$\\
\hline

$ {\rm B_T^2}(\T {4}{05}\big(0, -\frac{2} 3\big)) $&$=$
&
$\left\langle  
\Dt 11, \Dt 21,  2 \Dt 13 + \Dt 22 - 3\Dt 31  
\right\rangle$\\

$ {\rm Z_T^2}(\T {4}{05}\big(0, -\frac{2} 3\big)) $&$=$
&
$\left\langle
     \Dt 11, \Dt 12, \Dt 13+\Dt 22, \Dt 21, \Dt 22 + 3 \Dt 31, \Dt 23,  \Dt 14 + \Dt 32
\right\rangle$\\

$ {\rm H_T^2}(\T {4}{05}\big(0, -\frac{2} 3\big)) $&$=$&
$\left\langle
      \Dl 12,  \Dl 13 +  \Dl 22,   \Dl 23, \Dl 14 + \Dl 32
\right\rangle$\\
\hline

$ {\rm B_T^2}( \T {4}{05}\big(0, -\frac{1} 5\big)) $&$=$
&
$\left\langle 
\Dt 11, \Dt 21,  -\frac{1} 5 \Dt 13 + \frac{2} {15}\Dt 22 + \Dt 31 
\right\rangle$\\

$ {\rm Z_T^2}( \T {4}{05}\big(0, -\frac{1} 5\big)) $&$=$
&
$\left\langle
     \Dt 11, \Dt 12, \Dt 13+\Dt 22, \Dt 21, \Dt 23, \Dt 22 + 3 \Dt 31,   \Dt 14 -\frac{2} 5 \Dt 32 -\frac{7} 2 \Dt 41
\right\rangle$\\

$ {\rm H_T^2}( \T {4}{05}\big(0, -\frac{1} 5\big)) $&$=$
&
$\left\langle
      \Dl 12,  \Dl 13 +  \Dl 22, \Dl 23,   \Dl 14 -\frac{2} 5 \Dl 32 -\frac{7} 2 \Dl 41
\right\rangle$\\
\hline

$\alpha \neq \frac 1 3; - \frac 23; - \frac 1 5$&&\\
$ {\rm B_T^2}(\T {4}{05}\big(0, \alpha))$&$=$
&
$\left\langle 
\Dt 11, \Dt 21,  \alpha \Dt 13+ (\alpha+\frac 2 3) \Dt 22 + \Dt 31 \right\rangle$\\

$ {\rm Z_T^2}(\T {4}{05}\big(0, \alpha))$&$=$
&
$\left\langle
     \Dt 11, \Dt 12, \Dt 13+\Dt 22, \Dt 21, \Dt 22 + 3 \Dt 31,  \Dt 23
\right\rangle$\\

$ {\rm H_T^2}(\T {4}{05}\big(0, \alpha))$&$=$
&
$\left\langle
\Dl 12,  \Dl 13 + \Dl 22,  \Dl 23 \right\rangle$ \\

\hline

$  {\rm B_T^2}(\T {4}{05}(-1, \alpha)) $&$=$
&
$\left\langle 
\Dt 11, -\Dt 12 + \Dt 21,  \alpha \Dt 13+ (\alpha+\frac 2 3) \Dt 22 + \Dt 31 \right\rangle$\\
$  {\rm Z_T^2}(\T {4}{05}(-1, \alpha)) $&$=$
&
$\left\langle
\begin{array}{l}     
\Dt 11, \Dt 12, \Dt 13+\Dt 22, \Dt 21, 2 \Dt 22 + 3\Dt 31, \\ \Dt 14 - \frac {3\alpha+4} {3} \Dt 23 + \Dt 32, \frac{3\alpha +1} { 3} \Dt 23 - \frac{3\alpha +1} { 2} \Dt 32 + \Dt 41
\end{array}
\right\rangle$\\

$  {\rm H_T^2}(\T {4}{05}(-1, \alpha)) $&$=$
&
$\left\langle
\begin{array}{l}     
      \Dl 12,  \Dl 13 + \Dl 22, \Dl 14 - \frac {3\alpha+4} {3} \Dl 23 + \Dl 32, \\
      \frac{3\alpha +1} { 3} \Dl 23 - \frac{3\alpha +1} { 2} \Dl 32 + \Dl 41
\end{array}     
\right\rangle$\\
\hline

\multicolumn{3}{|l|}{ $\lambda \neq 0; -1; \frac 5 2, \alpha=\frac{4\lambda-1}{5-2\lambda}$}\\

$ {\rm B_T^2}(\T {4}{05}(\lambda, \alpha))$&$=$&
$\left\langle 
\Dt 11, \lambda \Dt 12 + \Dt 21,   \alpha \Dt 13+ \big(\alpha + \frac{1-\lambda} 3\big) \Dt 22 + \Dt 31  \right\rangle$\\

$ {\rm Z_T^2}(\T {4}{05}(\lambda, \alpha))$&$=$
&
$\left\langle
\begin{array}{l}
     \Dt 11, \Dt 12, \Dt 13+\Dt 22, \Dt 21, 
          (1-\lambda) \Dt 22 + 3\Dt 31, \\ \Dt 14 + \frac {4\lambda^2-2\lambda+7} {3\lambda(5-2\lambda)} \Dt 23 + \Dt 32,        \frac {2(\lambda-1)^2} {3\lambda(5-2\lambda)} \Dt 23 -\frac {2(\lambda-1)} {5-2\lambda} \Dt 32 + \Dt 41
\end{array}\right\rangle$\\

$ {\rm H_T^2}(\T {4}{05}(\lambda, \alpha))$&$=$
&
$\left\langle
\begin{array}{l}
      \Dl 12,  \Dl 13 + \Dl 22, 
        \Dl 14 + \frac {4\lambda^2-2\lambda+7} {3\lambda(5-2\lambda)} \Dl 23 + \Dl 32
          \\  \frac {2(\lambda-1)^2} {3\lambda(5-2\lambda)} \Dl 23 -\frac {2(\lambda-1)} {5-2\lambda} \Dl 32 + \Dl 41
\end{array}\right\rangle$
\\
\hline

\multicolumn{3}{|l|}{$\lambda \neq 0; -1; (\lambda; \alpha) \neq (1; 0); \alpha\neq\frac{4\lambda-1}{5-2\lambda}$}\\

$  {\rm B_T^2}(\T {4}{05}(\lambda, \alpha))$&$=$&
$\left\langle
\Dt 11, \lambda \Dt 12 + \Dt 21,  \alpha \Dt 13+ (\alpha +\frac{1-\lambda} 3) \Dt 22 + \Dt 31  \right\rangle$\\

$  {\rm Z_T^2}(\T {4}{05}(\lambda, \alpha))$&$=$&
$\left\langle
     \Dt 11, \Dt 12, \Dt 13+\Dt 22, \Dt 21, 
          (1-\lambda) \Dt 22 + 3\Dt 31,  \Dt 14 + \frac {3\alpha - 2(\lambda-1)} {3\lambda} \Dt 23 + \Dt 32
\right\rangle$\\

$  {\rm H_T^2}(\T {4}{05}(\lambda, \alpha))$&$=$&
$\left\langle
      \Dl 12,  \Dl 13 + \Dl 22, 
       \Dl 14 + \frac  {3\alpha - 2(\lambda-1)}{3\lambda}\Dl 23 + \Dl 32
\right\rangle$

\\
\hline

$ {\rm B_T^2}( \T {4}{06}(0)) $&$=$
&
$\left\langle 
 \Dt 11,  \Dt 21,  6 \Dt 12-\Dt 13 + \Dt 22 + 6\Dt 31 
 \right\rangle$\\
 
$ {\rm Z_T^2}(\T {4}{06}(0)) $&$=$
&
$\left\langle
     \Dt 11, \Dt 12, \Dt 13+\Dt 22, \Dt 21, \Dt 22 + 3\Dt 31,  \Dt 23
\right\rangle$\\

$ {\rm H_T^2}(\T {4}{06}(0)) $&$=$
&
$\left\langle
      \Dl 12, \Dl 23 + \Dl 22,  \Dl 23
\right\rangle$\\
\hline

$ {\rm B_T^2}(  \T {4}{06}(1)) $&$=$&
$\left\langle  
\Dt 11,  \Dt 12 + \Dt 21,  \Dt 12+ \Dt 13+ \Dt 22 + \Dt 31 \right\rangle$\\
$ {\rm Z_T^2}(  \T {4}{06}(1)) $&$=$&
$\left\langle
     \Dt 11, \Dt 12, \Dt 13+\Dt 22, \Dt 21,  \Dt 31, 
           \Dt 14 + \Dt 23+\Dt 32,   3\Dt 41-\Dt 22
\right\rangle$\\

$ {\rm H_T^2}(  \T {4}{06}(1)) $&$=$&
$\left\langle
 \Dl 12,   \Dl 14 + \Dl 23+\Dl 32, 
          \Dl 31,    3\Dl 41-\Dl 22
\right\rangle$\\
\hline

$ {\rm B_T^2}( \T {4}{06}(-1)) $&$=$
&
$\left\langle 
\Dt 11,  -\Dt 12 + \Dt 21,  3\Dt 12 -2 \Dt 13 + 3\Dt 31 \right\rangle$\\

$ {\rm Z_T^2}( \T {4}{06}(-1)) $&$=$
&
$\left\langle
\begin{array}{l}     
\Dt 11, \Dt 12, \Dt 13+\Dt 22, \Dt 21,  3 \Dt 14-2 \Dt 23+3 \Dt 32, \\ 2\Dt 22 + 3\Dt 31,
        2 \Dt 13 -3 \Dt 14 + 6\Dt 41
        \end{array}
\right\rangle$\\

$ {\rm H_T^2}( \T {4}{06}(-1)) $&$=$&
$\left\langle
 \Dl 12, \Dl 13+\Dl 22,   3\Dl 14-2 \Dl 23+3\Dl 32,
2 \Dl 13 -3 \Dl 14 +6\Dl 41
\right\rangle$\\
\hline

$ {\rm B_T^2}(  \T {4}{06}(-\frac{1}{2}))  $&$=$
&
$\left\langle 
\Dt 11, - \Dt 12 + 2\Dt 21,  2\Dt 12-\Dt 13+2\Dt 31  \right\rangle$\\

$ {\rm Z_T^2}(  \T {4}{06}(-\frac{1}{2}))  $&$=$
&
$\left\langle
\begin{array}{l}
     \Dt 11, \Dt 12, \Dt 13+\Dt 22, \Dt 21,  \Dt 14-\Dt 23 + \Dt 32, \\
     \Dt 22 + 2\Dt 31,
            2   \Dt 13 -3 \Dt 14 + 6 \Dt 41
            \end{array}
\right\rangle$\\

$ {\rm H_T^2}(  \T {4}{06}(-\frac{1}{2}))  $&$=$
&
$\left\langle
\Dl 12, \Dl 13+\Dl 22,   \Dl 14-\Dl 23 + \Dl 32,
 2 \Dl 13-3 \Dl 14 + 6\Dl 41
\right\rangle$
\\ \hline

\multicolumn{3}{|l|}{$ \lambda \neq 0; 1; -1; -\frac 1 2$}\\

$ {\rm B_T^2}(\T {4}{06}(\lambda)) $&$=$
&
$\left\langle 
 \Dt 11, \lambda \Dt 12 + \Dt 21,  \Dt 12+ \frac{2\lambda^2+5\lambda-1}{6}\Dt 13+   \frac{(2\lambda+1)(\lambda+1)}{6}\Dt 22 + \Dt 31  \right\rangle$\\

$ {\rm Z_T^2}(\T {4}{06}(\lambda)) $&$=$
&
$\left\langle
\begin{array}{l}     
\Dt 11, \Dt 12, \Dt 13+\Dt 22, \Dt 21,            (1-\lambda)\Dt 22 + 3\Dt 31,  \\
6\lambda \Dt 14 + (2\lambda^2+\lambda+3)\Dt 23 + 6\lambda\Dt 32
\end{array}
\right\rangle$\\

$ {\rm H_T^2}(\T {4}{06}(\lambda)) $&$=$
&
$\left\langle
\Dl 12, \Dt 13+\Dt 22, 6\lambda \Dl 14 + (2\lambda^2+\lambda+3) \Dl 23 + 6\lambda \Dl 32
\right\rangle$\\
\hline

$ {\rm B_T^2}( \T {4}{07}(-1)) $&$=$
&
$\left\langle
     \Dt 11, -\Dt 12+\Dt 21, \Dt 13+\Dt 22
\right\rangle$\\

$ {\rm Z_T^2}( \T {4}{07}(-1)) $&$=$
&
$\left\langle
     \Dt 11, \Dt 12, \Dt 13+\Dt 22, \Dt 21,  2\Dt 22 + 3\Dt 31,
            \Dt 14 - \Dt 23,  2 \Dt 23 - 3\Dt 32 + \Dt 41
\right\rangle$\\

$ {\rm H_T^2}( \T {4}{07}(-1)) $&$=$
&
$\left\langle
      \Dl 12, 2\Dl 22 + 3\Dl 31, 
            \Dl 14 - \Dl 23, 2 \Dl 23 - 3\Dl 32 + \Dl 41
\right\rangle$\\
\hline

$  {\rm B_T^2}( \T {4}{07}(\lambda)_{\lambda \neq -1}) $&$=$
&
$\left\langle
     \Dt 11, \lambda\Dt 12+\Dt 21, \Dt 13+\Dt 22
\right\rangle$\\

$  {\rm Z_T^2}( \T {4}{07}(\lambda)_{\lambda \neq -1}) $&$=$
&
$\left\langle
     \Dt 11, \Dt 12, \Dt 13+\Dt 22, \Dt 21, (1-\lambda)\Dt 22 + 3\Dt 31,
           \lambda \Dt 14 + \Dt 23
 \right\rangle$\\

$  {\rm H_T^2}( \T {4}{07}(\lambda)_{\lambda \neq -1}) $&$=$
&
$\left\langle
      \Dl 12, (1-\lambda)\Dl 22 + 3\Dl 31, 
           \lambda \Dl 14 + \Dl 23
\right\rangle$ \\
\hline

$  {\rm B_T^2}( \T {4}{08} )$&$=$
&
$\left\langle  \Dt 11, \Dt 21,  \Dt 23  \right\rangle$\\

$  {\rm Z_T^2}( \T {4}{08} )$&$=$
&
$\left\langle  \Dt 11, \Dt 12 , \Dt 13 + \Dt 22, \Dt 21, \Dt 22+ 3\Dt 31, \Dt 23, \Dt 24 \right\rangle$\\

$  {\rm H_T^2}( \T {4}{08} )$&$=$
&
$\left\langle    \Dl 13 + \Dl 22,  \Dl 21,  \Dl 22+ 3\Dl 31,  \Dl 24   \right\rangle$
\\
\hline

$  {\rm B_T^2}( \T {4}{09}) $&$=$
&
$\left\langle 
\Dt 11, \Dt 12 + \Dt 23,  \Dt 21   \right\rangle$\\

$  {\rm Z_T^2}( \T {4}{09}) $&$=$
&
$\left\langle 
 \Dt 11, \Dt 12 , \Dt 13 - 3 \Dt 31, \Dt 21,  \Dt 22+ 3\Dt 31,  \Dt 23 \right\rangle$\\

$  {\rm H_T^2}( \T {4}{09}) $&$=$
&
$\left\langle    \Dl 12 , \Dl 13 - 3 \Dl 31,   \Dl 22+ 3\Dl 31  \right\rangle$\\

\hline
$  {\rm B_T^2}(\T {4}{10}(\alpha)) $&$=$
&
$\left\langle 
 \Dt 11,  \Dt 21,   \alpha \Dt 12 + \frac 1 3 \Dt 22+ \Dt 23 +\Dt 31  \right\rangle$\\

$  {\rm Z_T^2}(\T {4}{10}(\alpha)) $&$=$
&
$\left\langle
     \Dt 11, \Dt 12, \Dt 13+\Dt 22,  \Dt 21, \Dt 22 + 3\Dt 31,
           \Dt 23 \right\rangle$\\

$  {\rm H_T^2}(\T {4}{10}(\alpha)) $&$=$
&
$\left\langle
      \Dl 12, \Dl 22 + 3\Dl 31,           \Dl 23
\right\rangle$\\
\hline
\end{longtable}

{\tiny
\begin{longtable}{l llll}
\multicolumn{5}{c}{ \mbox{ {\bf \large{Table C.}}
{\it \large{The list of $5$-dimensional one-generated nilpotent terminal algebras.}}}} \\
\multicolumn{5}{c}{  } \\

    \hline $\T{5}{01} $ &  $e_1 e_1 = e_2$ & $e_1 e_2=e_3$ & $e_1 e_3=e_4$ & $e_2 e_1=e_5$\\

    \hline $\T{5}{02}(\alpha) $ & 
        $e_1 e_1 = e_2$ & $e_1 e_2=e_3$ & $e_1 e_3=\alpha e_5$ &\\ 
        & $e_2 e_1=e_4$ & $e_2 e_2=e_5$  & $e_3 e_1=-3e_5$ &\\

    \hline $\T{5}{03} $ &  
    $e_1 e_1 = e_2$ & $e_1 e_2=e_3$ & $e_1 e_3=e_4$ \\ 
    &  $e_2 e_2=e_5$ & $e_3 e_1=-3e_5$\\

    \hline $\T{5}{04} $ &  
    $e_1 e_1 = e_2$ & $e_1 e_2=e_3$ & $e_1 e_3=e_4$  \\
        & $e_2 e_1=e_4$ & $e_2 e_2=e_5$   & $e_3 e_1=-3e_5$ \\

    \hline $\T{5}{05}(\lambda) $ &  
    $e_1 e_1 = e_2$ & $e_1 e_2=  \lambda e_3+e_4$ & $e_1 e_3= e_5$ \\ 
    & $e_2 e_1=e_3$  & $e_2 e_2= e_5$\\
    
    \hline $\T{5}{06}(\lambda,\alpha) $ &  
    $e_1 e_1 = e_2$& $e_1 e_2= \lambda e_3+e_4$ & $e_1 e_3= \alpha e_5$ \\
        & $e_2 e_1=e_3$ & $e_2 e_2= \big(\alpha+\frac{1-\lambda} 3 \Big) e_5$& $e_3 e_1=  e_5$ \\
        
    \hline $\T{5}{07}(\lambda) $ &  
    $e_1 e_1 = e_2$ & $e_1 e_2= \lambda e_3$ & $e_1 e_3=  e_4$ \\
        & $e_2 e_1=e_3$ & $e_2 e_2=  e_4 + \frac{1-\lambda} 3e_5$& $e_3 e_1=  e_5$ \\
        
    \hline $\T{5}{08}(\lambda) $ &  
    $e_1 e_1 = e_2$& $e_1 e_2= \lambda e_3+e_5$ & $e_1 e_3=  e_4$ \\
        & $e_2 e_1=e_3$ & $e_2 e_2=  e_4 + \frac{1-\lambda} {3} e_5$& $e_3 e_1=  e_5$ \\

    \hline $\T{5}{09} $ &  
    $e_1 e_1 = e_2$ & $e_1 e_2= e_4$ & $e_2 e_1=e_3$ & $e_2 e_3=e_5$ \\
    
    \hline $\T{5}{10} $ &  
    $e_1 e_1 = e_2$ & $e_1 e_2= e_4$ & $e_2 e_1=e_3$  \\ 
        & $e_2 e_2= e_5$ & $e_2 e_3=e_5$ & $e_3 e_1=3e_5$ \\
        
    \hline $\T{5}{11} $ &  
    $e_1 e_1 = e_2$ & $e_1 e_3= e_4$ & $e_2 e_1=e_3$ \\
    & $e_2 e_2=e_4$ & $e_2 e_3=e_5$\\
    
    \hline $\T{5}{12} $ &  
    $e_1 e_1 = e_2$ & $e_1 e_2= e_5$ & $e_1 e_3=e_4$  \\ 
        & $e_2 e_1=e_3$ & $e_2 e_2=e_4$  & $e_2 e_3=e_5$  \\
        
    \hline $\T{5}{13}(\alpha) $ &  
    $e_1 e_1 = e_2$ & $e_1 e_2=  \alpha e_5$ & $e_1 e_3= e_4$ & $e_2 e_1=e_3$   \\ 
        & $e_2 e_2=e_4 + e_5$& $e_2 e_3=e_5$& $e_3 e_1=3e_4$   \\
        
    \hline $\T{5}{14}(\alpha) $ &  
    $e_1 e_1 = e_2$ & $e_1 e_3= \alpha e_4$ & $e_2 e_1=e_3$  \\
        & $e_2 e_2= (\alpha +1)e_4$ & $e_2 e_3=e_5$& $e_3 e_1=3e_4$  \\
        
    \hline $\T{5}{15}(\alpha) $ &  
    $e_1 e_1 = e_2$ & $e_1 e_2=  e_4$  & $e_1 e_3= \alpha e_4$  & $e_2 e_1=e_3$   \\
        & $e_2 e_2= (\alpha+1)e_4$ & $e_2 e_3=e_5$& $e_3 e_1=3 e_4$  \\
        
    \hline $\T{5}{16}(\alpha, \beta) $ &  
    $e_1 e_1 = e_2$& $e_1 e_2= \beta e_4 + e_5$ & $e_1 e_3= \alpha e_4$  & $e_2 e_1=e_3$  \\ 
        &      $e_2 e_2= (\alpha+1)e_4$ & $e_2 e_3=e_5$& $e_3 e_1=3e_4$   \\
        
    \hline $\T{5}{17}(\alpha, \beta) $ &  
    $e_1 e_1 = e_2$ & $e_1 e_2=e_4$ & $e_1 e_4=e_5$ & $e_2 e_1=e_3$ \\
        & $e_2 e_2= \alpha e_5$ & $e_2 e_3 = e_5$ & $e_3 e_1=(3\alpha+\beta) e_5$ & $e_4 e_1= \beta e_5$\\
        
    \hline $\T{5}{18}(\alpha) $ &  
    $e_1 e_1 = e_2$ & $e_1 e_2=e_4$  & $e_2 e_1=e_3$ & $e_2 e_2= \alpha e_5$   \\
        & $e_2 e_3 = e_5$  &  $e_3 e_1=(3\alpha+1) e_5$ & $e_4 e_1=  e_5$ \\
        
    \hline $\T{5}{19}(\alpha, \beta, \gamma) $ &  
    $e_1 e_1 = e_2$ & $e_1 e_2=e_4$ & $e_1 e_3=\alpha e_5$ & $e_1 e_4=\beta e_5$   \\
    $(\alpha; \gamma)\neq(0 ; - \frac 1 3 )$ & $e_2 e_1=e_3$ & $e_2 e_2= (\alpha+\gamma) e_5$ & $e_3 e_1=(3\gamma+1) e_5$ & $e_4 e_1= e_5$  \\
    
    \hline $\T{5}{20}(\alpha, \beta) $ &  
    $e_1 e_1 = e_2$ & $e_1 e_2=e_4$ & $e_1 e_3=\alpha e_5$ & $e_1 e_4= e_5$    \\
    $(\alpha; \beta)\neq(0 ; 0 )$ &  $e_2 e_1=e_3$ & $e_2 e_2= (\alpha+\beta) e_5$ & $e_3 e_1=3\beta e_5$  \\
    
    \hline $\T{5}{21}(\alpha) $ &  
    $e_1 e_1 = e_2$ & $e_1 e_2=e_3$ & $e_1 e_3=\alpha e_4$ & $e_1 e_4=e_5$  \\
        & $e_2 e_2= e_4$  & $e_2 e_3 = 2 e_5$ & $e_3 e_1 = -3 e_4$ & $e_3 e_2=-3e_5$ \\
        
    \hline $\T{5}{22}(\alpha) $ &  
    $e_1 e_1 = e_2$ & $e_1 e_2=e_3$ & $e_1 e_3=\alpha e_4 +e_5$ & $e_1 e_4=e_5$  \\
        & $e_2 e_2= e_4$  & $e_2 e_3 = 2 e_5$ & $e_3 e_1 = -3 e_4$ & $e_3 e_2=-3e_5$  \\
        
    \hline $\T{5}{23}(\beta) $ &  $e_1 e_1 = e_2$ & $e_1 e_2=e_3$ & $e_1 e_3=6 e_4$ \\& $e_1 e_4=\beta e_5$ 
    & $e_2 e_2= e_4$ & $e_2 e_3 = (2\beta+1) e_5$ \\& $e_3 e_1 = -3 e_4$  
    & $e_3 e_2 = -3(\beta+1) e_5$  & $e_4 e_1=e_5$&\\
    
    \hline $\T{5}{24}(\beta) $ &  
    $e_1 e_1 = e_2$ & $e_1 e_2=e_3$ & $e_1 e_3=6 e_4 + e_5$ \\
    & $e_1 e_4=\beta e_5$ & $e_2 e_2= e_4$     &  $e_2 e_3 = (2\beta+1) e_5$ \\& $e_3 e_1 = -3 e_4$  
    & $e_3 e_2 = -3(\beta+1) e_5$  & $e_4 e_1=e_5$&\\
    
    \hline $\T{5}{25}(\beta) $ &  
    $e_1 e_1 = e_2$ & $e_1 e_2=e_3$ & $e_1 e_3=6 e_4 + \beta e_5$ & $e_1 e_4=-2 e_5$ \\
    & $e_2 e_1= e_5$ & $e_2 e_2= e_4$ &  $e_2 e_3 = -3 e_5$ \\
    & $e_3 e_1 = -3 e_4$  & $e_3 e_2 = 3 e_5$  & $e_4 e_1=e_5$\\
    
    \hline $\T{5}{26} $ &  $e_1 e_1 = e_2$ & $e_1 e_2=e_3$ & $e_1 e_3= e_4$  \\ 
        & $e_1 e_4=e_5$ & $e_2 e_2= e_5$ & $e_3 e_1 = -3e_5$  \\
        
    \hline $\T{5}{27} $ &  $e_1 e_1 = e_2$ & $e_1 e_2=e_3$ & $e_1 e_3= e_4$ & $e_1 e_4=e_5$   \\
    
    \hline $\T{5}{28} $ &  $e_1 e_1 = e_2$ & $e_1 e_2=e_3$ & $e_1 e_3= e_4$ \\ 
    & $e_1 e_4=e_5$ & $e_2 e_1= e_5$\\
    
    \hline $\T{5}{29}(\alpha) $ &  
    $e_1 e_1 = e_2$ & $e_1 e_2=e_3$ & $e_1 e_3=e_4$ & $e_1 e_4=e_5$   \\ 
    $\alpha\neq 1$    & $e_2 e_1= e_4$& $e_2 e_2 = \alpha e_5$ & $e_3 e_1=3(1 - \alpha)e_5$   \\
    
    \hline $\T{5}{30}(\alpha) $ &  
    $e_1 e_1 = e_2$ & $e_1 e_2=e_3$ & $e_1 e_3=e_4+\alpha e_5$   \\
        & $e_1 e_4=e_5$ & $e_2 e_1= e_4$ & $e_2 e_2 =  e_5$   \\
        
    \hline $\T{5}{31}(\beta) $ &  $e_1 e_1 = e_2$ & $e_1 e_2= e_3+\beta e_5$ & $e_1 e_3=e_5$ & $e_2 e_1=e_3$  \\
        & $e_2 e_2=e_5$ & $e_3 e_1= e_4$& $e_3 e_4= e_5$  \\
        
    \hline $\T{5}{32} $ &  $e_1 e_1 = e_2$ & $e_1 e_2= e_3+ e_5$ &  $e_2 e_1=e_3$ \\
    & $e_3 e_1= e_4$ & $e_3 e_4= e_5$\\
    
    \hline $\T{5}{33} $ &  $e_1 e_1 = e_2$ & $e_1 e_2= e_3$ &  $e_2 e_1=e_3$ \\ 
    & $e_3 e_1= e_4$& $e_3 e_4= e_5$\\
    
    \hline $\T{5}{34}(\beta) $ &  $e_1 e_1 = e_2$ & $e_1 e_2=\beta e_5$ & $e_1 e_3= \frac 1 3 e_4$ & $e_2 e_1=e_3$ \\ & $e_2 e_2=\frac{2} 3 e_4$  & $e_2 e_3= e_5$ & $e_2 e_4= e_5$\\ 
    & $e_3 e_1= e_4$ & $e_3 e_3= e_5$ &\\
    
    \hline $\T{5}{35} $ &  
    $e_1 e_1 = e_2$ & $e_1 e_2=e_5$ & $e_1 e_3= \frac 1 3 e_4$ & $e_2 e_1=e_3$  \\
        & $e_2 e_2=\frac{2} 3 e_4$  &  $e_2 e_4= e_5$ & $e_3 e_1= e_4$ & $e_3 e_3= e_5$   \\
        
    \hline $\T{5}{36} $ &  
    $e_1 e_1 = e_2$  & $e_1 e_3= \frac 1 3 e_4$ & $e_2 e_1=e_3$ & $e_2 e_2=\frac{2} 3 e_4$  \\
        &  $e_2 e_4= e_5$  & $e_3 e_1= e_4$ & $e_3 e_3= e_5$  \\
        
     \hline $\T{5}{37}(\beta) $ &  $e_1 e_1 = e_2$  & $e_1 e_3= - \frac 2 3 e_4+e_5$ & $e_2 e_1=e_3$ & $e_1 e_4=e_5$   \\
        & $e_2 e_2=-\frac{1} 3 e_4+e_5$    & $e_2 e_3= \beta e_5$ & $e_3 e_1= e_4$ & $e_3 e_2= e_5$  \\
        
    \hline $\T{5}{38}(\beta) $ &  
    $e_1 e_1 = e_2$  & $e_1 e_3= - \frac 2 3 e_4$ & $e_2 e_1=e_3$ & $e_1 e_4=e_5$   \\
        & $e_2 e_2=-\frac{1} 3 e_4$ &  $e_2 e_3= \beta e_5$ & $e_3 e_1= e_4$ & $e_3 e_2= e_5$   \\
        
    \hline $\T{5}{39}(\beta) $ &  $e_1 e_1 = e_2$ & $e_1 e_3= -\frac 1 5 e_4+e_5$ & $e_1 e_4=e_5$ & $e_2 e_1=e_3$\\
    & $e_2 e_2=\frac{2} {15} e_4+e_5$      & $e_2 e_3= \beta e_5$  & $e_3 e_1= e_4$ \\
    & $e_3 e_2= -\frac 2 5 e_5$  & $e_4 e_1= -\frac 7 2 e_5$ & \\
    
    \hline $\T{5}{40}(\beta) $ &  
    $e_1 e_1 = e_2$ & $e_1 e_3= -\frac 1 5 e_4$ & $e_1 e_4=e_5$ \\& $e_2 e_1=e_3$ 
    & $e_2 e_2=\frac{2} {15} e_4$ & $e_2 e_3= \beta e_5$ & \\$e_3 e_1= e_4$ 
    & $e_3 e_2= -\frac 2 5 e_5$  & $e_4 e_1= -\frac 7 2 e_5$ & \\

    \hline $\T{5}{41}(\alpha, \beta) $ &  
    $e_1 e_1 = e_2$ & $e_1 e_2=- e_3$ & $e_1 e_3=\alpha e_4+e_5$ & $e_1 e_4=\beta e_5$ \\
    & $e_2 e_1=e_3$ & $e_2 e_2=\big(\alpha+\frac{2} 3\big) e_4+e_5$ &  \multicolumn{2}{l}{ $e_2 e_3=\big(\frac{3\alpha+1} 3 -\frac{\beta(3\alpha+4)} 3 \big) e_5$} \\
    & $e_3 e_1= e_4$ & $e_3 e_2=\big(\beta - \frac{3\alpha+1} 2\big) e_5$  & $e_4e_1 =e_5$ \\
    
    \hline $\T{5}{42}(\alpha, \beta) $ &  
    $e_1 e_1 = e_2$ & $e_1 e_2=- e_3$ & $e_1 e_3=\alpha e_4$ &  $e_1 e_4=\beta e_5$ \\
    & $e_2 e_1=e_3$ & $e_2 e_2=\big(\alpha+\frac{2} 3\big) e_4$ &  \multicolumn{2}{l}{$e_2 e_3=\big(\frac{3\alpha+1} 3 -\frac{\beta(3\alpha+4)} 3 \big) e_5$}\\ 
    & $e_3 e_1= e_4$ & $e_3 e_2=\big(\beta - \frac{3\alpha+1} 2\big) e_5$ &  $e_4e_1 =e_5$ \\
    
    \hline $\T{5}{43}(\alpha, \beta) $ &  
    $e_1 e_1 = e_2$ & $e_1 e_2=- e_3 + e_5$ & $e_1 e_3=\alpha e_4+\beta e_5$ & $e_1 e_4=-\frac {3\alpha+5} 4 e_5$\\
    $\alpha \neq -\frac 2 3$ & $e_2 e_1=e_3$  &     $e_2 e_2=\big(\alpha+\frac{2} 3\big) e_4+\beta e_5$    &  \multicolumn{2}{l}{$e_2 e_3=\frac{3\alpha^2+13\alpha+8} 4  e_5$}\\ & $e_3 e_1= e_4$ 
        & $e_3 e_2= -\frac{9\alpha+7} 4 e_5$ & $e_4e_1 =e_5$ \\

    \hline $\T{5}{44} $ &  
    $e_1 e_1 = e_2$ & $e_1 e_2= -  e_3+e_5$ & $e_1 e_3=- \frac 2 3  e_4$ \\& $e_1 e_4= -\frac 3 4 e_5$ 
    & $e_2 e_1=e_3$ & $e_2 e_3= \frac 1 6 e_5$ \\& $e_3 e_1 =e_4$ 
    & $e_3 e_2= - \frac 1 4 e_5$ & $e_4 e_1= e_5$ &    \\
    
    \hline $\T{5}{45} $ &  
    $e_1 e_1 = e_2$ & $e_1 e_2= -  e_3$ & $e_1 e_3=- \frac 2 3  e_4+e_5$ \\& $e_1 e_4=  e_5$ 
    & $e_2 e_1=e_3$ & $e_2 e_2= e_5$ \\& $e_2 e_3= -\frac 2 3 e_5$ 
    & $e_3 e_1 =e_4$ & $e_3 e_2= e_5$    \\
    
    \hline $\T{5}{46}(\lambda, \beta) $ &  
    $e_1 e_1 = e_2$ & $e_1 e_2=\lambda e_3+e_5$ & $e_1 e_3=\frac{4\lambda-1}{5-2\lambda} e_4+\beta e_5$ & $e_1 e_4=\frac{4\lambda-1}{5-2\lambda} e_5$ \\
    $\lambda \neq 1, -2, -\frac 1 2, \frac  5 2$ &     $e_2 e_1=e_3$ &   $e_2 e_2=\frac{(2\lambda+1)(\lambda+2)}{3(5-2\lambda)} e_4+\beta e_5$ &$e_2 e_3=\frac {2(2\lambda+1)(3\lambda^2+1)}{3\lambda(5-2\lambda)^2}e_5$\\
    & $e_3 e_1= e_4$ &   $e_3 e_2=\frac {2\lambda+1}{5-2\lambda}e_5$     &$e_4 e_1= e_5$  \\
    
    \hline $\T{5}{47}(\lambda, \beta) $ &  
    $e_1 e_1 = e_2$ & $e_1 e_2=\lambda e_3$ & $e_1 e_3=\frac{4\lambda-1}{5-2\lambda} e_4+ e_5$    & $e_1 e_4=\beta e_5$  \\
    
    $\lambda \neq 1, \frac   5 2$ & $e_2 e_1=e_3 $ & $e_2 e_2=\frac{(2\lambda+1)(\lambda+2)}{3(5-2\lambda)} e_4+ e_5$ & \multicolumn{2}{l}{$e_2 e_3=\frac {\beta(4\lambda^2-2\lambda+7)+2(\lambda-1)^2}{3\lambda(5-2\lambda)}e_5$}\\ 
    & $e_3 e_1= e_4$  
    &$e_3 e_2=\Big(\beta - \frac {2(\lambda-1)}{5-2\lambda}\Big)e_5$        & $e_4 e_1= e_5$ \\
    
    \hline $\T{5}{48}(\lambda, \beta) $ &  
    $e_1 e_1 = e_2$ & $e_1 e_2=\lambda e_3$ & $e_1 e_3=\frac{4\lambda-1}{5-2\lambda} e_4$   & $e_1 e_4=\beta e_5$   \\ 
    $  \lambda \neq 0, \frac  5 2$ & $e_2 e_1=e_3$  & $e_2 e_2=\frac{(2\lambda+1)(\lambda+2)}{3(5-2\lambda)} e_4$ &\multicolumn{2}{l}{$e_2 e_3=\frac {\beta(4\lambda^2-2\lambda+7)+2(\lambda-1)^2}{3\lambda(5-2\lambda)}e_5$} \\
    &   $e_3 e_1= e_4$  & $e_3 e_2=\Big(\beta - \frac {2(\lambda-1)}{5-2\lambda}\Big)e_5$  & $e_4 e_1= e_5$ \\
    
    \hline $\T{5}{49} $ &  
    $e_1 e_1 = e_2$ & $e_1 e_2=e_3+e_5$ & $e_1 e_3= e_4$ & $e_1 e_4= e_5$ \\
    & $e_2 e_1=e_3$ & $e_2 e_2= e_4$ & $e_2 e_3=e_5$  \\
    &  $e_3 e_1= e_4$ &  $e_3 e_2= e_5$   & $e_4 e_1= e_5$ \\
            
    \hline $\T{5}{50} $ &  
    $e_1 e_1 = e_2$ & $e_1 e_2=e_3$ & $e_1 e_3= e_4+e_5$ & $e_1 e_4= e_5$ \\
    &   $e_2 e_1=e_3$ & $e_2 e_2= e_4+e_5$ & $e_2 e_3=e_5$ \\
    &  $e_3 e_1= e_4$ & $e_3 e_2= e_5$   &$e_4 e_1= e_5$ \\
    
    \hline $\T{5}{51}(\beta, \gamma) $ &  
    $e_1 e_1 = e_2$ & $e_1 e_2=-2 e_3 + e_5$ & $e_1 e_3= - e_4+\beta e_5$ & $e_1 e_4= \gamma e_5$ \\
    &  $e_2 e_1=e_3$ & $e_2 e_2= \beta e_5$   & $e_2 e_3=-\frac{3\gamma+2} 6 e_5$    &  $e_3 e_1= e_4$ \\ 
        & $e_3 e_2= \frac{3\gamma+2} 3 e_5$    & $e_4 e_1= e_5$  \\
        
    \hline $\T{5}{52}(\beta) $ &  
    $e_1 e_1 = e_2$ & $e_1 e_2=-2 e_3 + \beta e_5$ & $e_1 e_3= - e_4+e_5$ \\& $e_1 e_4= e_5$ 
     & $e_2 e_1=e_3$ & $e_2 e_2=  e_5$\\ &    $e_2 e_3=-\frac{1} 2 e_5$ 
     &  $e_3 e_1= e_4$ &   $e_3 e_2=  e_5$ &   \\

    \hline $\T{5}{53} $ &  
    $e_1 e_1 = e_2$ & $e_1 e_2=-\frac 1 2 e_3$ & $e_1 e_3= -\frac 1 2 e_4+e_5$ & $e_1 e_4= -\frac 1 2 e_5$  \\
        & $e_2 e_1=e_3$  & $e_2 e_2= e_5$& $e_3 e_1= e_4$ & $e_4 e_1= e_5$  \\
        
    \hline $\T{5}{54} $ &  
    $e_1 e_1 = e_2$ & $e_1 e_2=-\frac 1 2 e_3$ & $e_1 e_3= -\frac 1 2 e_4+e_5$ \\& $e_1 e_4= e_5$ 
    & $e_2 e_1=e_3$ &  $e_2 e_2= e_5$ \\& $e_2 e_3= - e_5$ 
    & $e_3 e_1= e_4$ & $e_3 e_2= e_5$ & \\
   	
   	\hline $\T{5}{55}(\lambda, \alpha)_{\lambda\neq -2,0} $ &  
   	$e_1 e_1 = e_2$ & $e_1 e_2=\lambda e_3$ & $e_1 e_3=\alpha e_4+e_5$ \\
   	$\alpha\neq \frac{2\lambda^2+5\lambda-1}{6}$ & $e_1 e_4= e_5$ 
   	&   $e_2 e_1=e_3$ & $e_2 e_2=\big(\alpha+\frac{1-\lambda}3\big) e_4+e_5$ \\ & $e_2e_3 =  \frac  {3\alpha-2(\lambda-1)}{3\lambda} e_5$ 
   	& $e_3 e_1= e_4$     & $e_3e_2 = e_5$& \\

    \hline $\T{5}{56}(\lambda, \alpha) $ &  $e_1 e_1 = e_2$ & $e_1 e_2=\lambda e_3$ & $e_1 e_3=\alpha e_4$  \\  $ \lambda\neq  0$& $e_1 e_4= e_5$ 
     & $e_2 e_1=e_3$  &  $e_2 e_2=\big(\alpha+\frac{1-\lambda}3\big) e_4$  \\& $e_2e_3 =  \frac  {3\alpha-2(\lambda-1)}{3\lambda} e_5$ 
     & $e_3 e_1= e_4$   & $e_3e_2 = e_5$& \\
        
    \hline $\T{5}{57}(\beta) $ &  
    $e_1 e_1 = e_2$ & $e_1 e_2=-2 e_3$ & $e_1 e_3=-\frac 1 2 e_4 + e_5$ \\
    & $e_1 e_4= e_5$
        &  $e_2 e_1=e_3$   & $e_2 e_2=\frac 1 2 e_4 + e_5$ \\
        & $e_2e_3 = - \frac  {3}{4} e_5$
        &  $e_3 e_1= e_4$ & $e_3e_2 = e_5$& \\

    \hline $\T{5}{58}(\alpha) $ &  
    $e_1 e_1 = e_2$ & $e_1 e_2=-2 e_3+e_5$ & $e_1 e_3=\alpha e_4$ \\& $e_1 e_4= e_5$    
    & $e_2 e_1=e_3$  &  $e_2 e_2=(\alpha+1) e_4$ \\& $e_3 e_1= e_4$ 
    & $e_2e_3 =  - \frac  {\alpha+2}{2} e_5$ & $e_3e_2 = e_5$ & \\
    
    \hline $\T{5}{59}(\beta) $ &  
    $e_1 e_1 = e_2$ & $e_1 e_2= e_3+e_4$ & $e_1 e_3= e_4$ & $e_1 e_4= \beta e_5$ \\  
    & $e_2 e_1=e_3$ & $e_2 e_2= e_4 -e_5$ &  $e_2 e_3= \beta e_5$ \\
    & $e_3 e_1= e_4$ & $e_3 e_2= \beta e_5$ & $e_4 e_1= 3 e_5$\\

    \hline $\T{5}{60}(\beta)_{\beta\neq0} $ &  
    $e_1 e_1 = e_2$ & $e_1 e_2= e_3+e_4$ & $e_1 e_3= e_4$ & $e_1 e_4= 3 e_5$ \\ 
    & $e_2 e_1=e_3$ & $e_2 e_2= e_4 -e_5$ &  $e_2 e_3= 3 e_5$ & $e_3 e_1= e_4 +  \beta e_5$ \\
    & $e_3 e_2= 3 e_5$ & $e_4 e_1= 3 e_5$\\
    
    \hline $\T{5}{61}(\beta)_{\beta\neq0} $ &  
    $e_1 e_1 = e_2$ & $e_1 e_2= e_3+e_4+\beta e_5$ & $e_1 e_3= e_4$ & $e_1 e_4= 3 e_5$ \\
    & $e_2 e_1=e_3$ & $e_2 e_2= e_4 -e_5$ &  $e_2 e_3= 3 e_5$ \\ 
    & $e_3 e_1= e_4 -6 e_5$ & $e_3 e_2= 3 e_5$ & $e_4 e_1= 3 e_5$\\
    
    \hline $\T{5}{62}(\beta)_{\beta\neq0} $ &  
    $e_1 e_1 = e_2$ & $e_1 e_2= - e_3+e_4$ & $e_1 e_3= -\frac 2 3 e_4+\beta e_5$ \\& $e_1 e_4=  3e_5$   
    &$e_2 e_1=e_3$     & $e_2 e_2= \beta e_5$\\ &  $e_2 e_3= -2 e_5$ 
    & $e_3 e_1= e_4$ & $e_3 e_2= 3 e_5$ & \\
    
    \hline $\T{5}{63}(\beta)_{\beta\neq0} $ &  
    $e_1 e_1 = e_2$ & $e_1 e_2= - e_3+e_4 + \beta e_5$ & $e_1 e_3= -\frac 2 3 e_4 - e_5$ & $e_1 e_4=  9 e_5$ \\ 
    & $e_2 e_1=e_3$     & $e_2 e_2= 3 e_5$ &  $e_2 e_3= - 2 e_5$ \\
    & $e_3 e_1= e_4$ & $e_3 e_2= 3 e_5$ & $e_4 e_1= -12 e_5$\\
    
    \hline $\T{5}{64}(\beta, \gamma) $ &  
    $e_1 e_1 = e_2$ & $e_1 e_2= - e_3+e_4$ & $e_1 e_3= -\frac 2 3 e_4 +(\beta+2) e_5$ & $e_1 e_4=3(\gamma-1) e_5$ \\
    & $e_2 e_1=e_3$  & $e_2 e_2= \beta e_5$ &  $e_2 e_3= - 2\gamma e_5$ \\
    & $e_3 e_1= e_4$ & $e_3 e_2= 3\gamma e_5$ & $e_4 e_1= 6 e_5$\\
    
    \hline $\T{5}{65}(\beta)_{\beta\neq0} $ &  
    $e_1 e_1 = e_2$ & $e_1 e_2= -\frac 1 2  e_3+e_4$ & $e_1 e_3= -\frac 1 2 e_4+\beta e_5$ \\
    & $e_1 e_4=  e_5$ & $e_2 e_1=e_3$  & $e_2 e_2= \beta e_5$ \\
    &  $e_2 e_3= - e_5$ & $e_3 e_1= e_4$ & $e_3 e_2= e_5$&\\
    
    \hline $\T{5}{66}(\beta)_{\beta\neq0} $ &  
    $e_1 e_1 = e_2$ & $e_1 e_2= -\frac 1 2  e_3+e_4 + \beta e_5$ & $e_1 e_3= -\frac 1 2 e_4+3 e_5$ & $e_1 e_4= -3e_5$ \\ 
    &$e_2 e_1=e_3$ & $e_2 e_2= e_5$  & $e_3 e_1= e_4$ & $e_4 e_1= 6e_5$\\
    
    \hline $\T{5}{67}(\beta, \gamma) $ &  
    $e_1 e_1 = e_2$ & $e_1 e_2= -\frac 1 2  e_3+e_4 $ & $e_1 e_3= -\frac 1 2 e_4+ (\beta+2) e_5$ & $e_1 e_4= (\gamma-3)e_5$ \\
    & $e_2 e_1=e_3$ & $e_2 e_2= \beta e_5$ & $e_2 e_3= -\gamma e_5$ & $e_3 e_1= e_4$ \\
    & $e_3 e_2= \gamma e_5$ & $e_4 e_1= 6e_5$\\
    
    \hline $\T{5}{68}(\beta)_{\beta\neq 0} $ &  
    $e_1 e_1 = e_2$ & $e_1 e_2= -2 e_3+e_4+\beta e_5$ & $e_1 e_3= -\frac 1 2e_4$ \\ & $e_1 e_4= 4 e_5$ 
    & $e_2 e_1=e_3$ & $e_2 e_2= \frac 1 2 e_4$ \\&  $e_2 e_3= -3 e_5$ 
    & $e_3 e_1= e_4$ & $e_3 e_2= 4e_5$&\\
    
    \hline $\T{5}{69}(\lambda) $ &  
    $e_1 e_1 = e_2$ & $e_1 e_2=\lambda e_3+e_4$ & $e_1 e_3= \frac{2\lb^2  + 5\lb  - 1}{6}e_4$ \\
    & $e_1 e_4= 6\lambda e_5$ 
    & $e_2 e_1=e_3$ & $e_2 e_2= \frac{(2\lb+1)(\lb+1)}{6}e_4$ \\& $e_2 e_3= (2\lb^2+\lb+3) e_5$ 
    & $e_3 e_1= e_4$  & $e_3 e_2= 6\lambda e_5$&\\
        
     \hline $\T{5}{70}(\lambda, \alpha) $ &  
     $e_1 e_1 = e_2$ & $e_1 e_2=\lambda e_3+\alpha e_5$ & $e_1 e_3= e_4$ & $e_1 e_4=\lambda e_5$\\
        & $e_2 e_1= e_3$& $e_2 e_2 = e_4+(1-\lambda)e_5$ & $e_2 e_3 = e_5$ & $e_3 e_1=3 e_5$\\
     
     \hline $\T{5}{71}(\lambda) $ &  
     $e_1 e_1 = e_2$ & $e_1 e_2=\lambda e_3 + e_5$ & $e_1 e_3= e_4$ & $e_1 e_4=\lambda e_5$ \\
        & $e_2 e_1= e_3$& $e_2 e_2 = e_4$ & $e_2 e_3 = e_5$\\
        
     \hline $\T{5}{72}(\lambda) $ &  
     $e_1 e_1 = e_2$ & $e_1 e_2=\lambda e_3$ & $e_1 e_3= e_4$ & $e_1 e_4=\lambda e_5$ \\
        & $e_2 e_1= e_3$& $e_2 e_2 = e_4$ & $e_2 e_3 = e_5$\\
        
     \hline $\T{5}{73}(\alpha,\beta) $ &  $e_1 e_1 = e_2$ & $e_1 e_2=- e_3+e_5$ & $e_1 e_3= e_4$ & $e_1 e_4=\alpha e_5$ \\
     & $e_2 e_1= e_3$ & $e_2 e_2 = e_4+2\beta e_5$ & $e_2 e_3 = (2-\alpha)e_5$ & $e_3 e_1=3\beta e_5$ \\ 
     & $e_3 e_2=-3 e_5$ & $e_4 e_1= e_5$\\
     
     \hline $\T{5}{74}(\alpha) $ &  
     $e_1 e_1 = e_2$ & $e_1 e_2=- e_3$ & $e_1 e_3= e_4$ & $e_1 e_4=\alpha e_5$ \\
     & $e_2 e_1= e_3$ & $e_2 e_2 = e_4+2 e_5$ & $e_2 e_3 = (2-\alpha)e_5$ \\
     & $e_3 e_1= 3 e_5$ & $e_3 e_2=-3 e_5$ & $e_4 e_1= e_5$\\
     
     \hline $\T{5}{75}(\alpha) $ &  
     $e_1 e_1 = e_2$ & $e_1 e_2=- e_3$ & $e_1 e_3= e_4$ \\
     & $e_1 e_4=\alpha e_5$ 
     & $e_2 e_1= e_3$ & $e_2 e_2 = e_4$ \\& $e_2 e_3 = (2-\alpha)e_5$ 
     & $e_3 e_2=-3 e_5$ & $e_4 e_1= e_5$ &\\
     
     \hline $\T{5}{76}(\alpha,\beta) $ &  
     $e_1 e_1 = e_2$ & $e_1 e_2=e_5$ & $e_1 e_3=\alpha e_5$ & $e_2 e_1=e_3$ \\
        & $e_2 e_2= (\alpha+\beta)e_5$& $e_2 e_3 = e_4$ & $e_2 e_4 =  e_5$ & $e_3 e_1=3\beta e_5$\\
        
     \hline $\T{5}{77}(\alpha) $ &  
     $e_1 e_1 = e_2$ &  $e_1 e_3= e_5$ & $e_2 e_1=e_3$ & $e_2 e_2= (\alpha+1)e_5$ \\
        & $e_2 e_3 = e_4$ & $e_2 e_4 =  e_5$ & $e_3 e_1=3\alpha e_5$ \\
     
     \hline $\T{5}{78} $ &  
     $e_1 e_1 = e_2$ &  $e_2 e_1=e_3$ & $e_2 e_2= e_5$  \\
        & $e_2 e_3 = e_4$ & $e_2 e_4 =  e_5$ & $e_3 e_1=3 e_5$ \\
        
     \hline $\T{5}{79} $ &  $e_1 e_1 = e_2$ &  $e_2 e_1=e_3$ & $e_2 e_3 = e_4$ & $e_2 e_4 =  e_5$\\ 
     \hline
\end{longtable}}

\end{document}